\documentclass{amsart}

\usepackage{amsthm,amsfonts,amsmath,amssymb,latexsym,epsfig}

\newtheorem{theorem}{Theorem}
\newtheorem{lemma}[theorem]{Lemma}
\newtheorem{proposition}[theorem]{Proposition}
\newtheorem{corollary}[theorem]{Corollary}
\newenvironment{example}{\medskip \refstepcounter{theorem}
\noindent  {\bf Example \thetheorem}.\rm}{\,}
\newenvironment{remark}{\medskip \refstepcounter{theorem}
\noindent  {\bf Remark \thetheorem}.\rm}{\,}
\newenvironment{definition}{\medskip \refstepcounter{theorem}
\noindent  {\bf Definition \thetheorem}.\rm}{\,}

\renewcommand{\thetheorem}{\thesection.\arabic{theorem}}

\theoremstyle{definition}

\theoremstyle{remark}

\def \cK{{\mathcal K}}

\def\wh{\widehat{\phantom{p}}}
\def \phip{\phi_p}
\def \di{p}
\def \ddi{q}
\def \sh{\sharp}
\def \la{\leftarrow}
\def \ua{\uparrow}
\def \tf{\tilde{f}}

\def \hU{\hat{U}}
\def \rs{\ }

\def \da{\downarrow}

\def \bZ{{\mathbb Z}}
\def \bG{{\mathbb G}}
\def \bmu{{\mathbb \mu}}
\def \D{\{\d,\dd\}}
\def \dd{\delta_q}

\def \<{\langle}
\def \>{\rangle}

\def \cU{\mathcal U}
\def \cF{\mathcal F}

\def \cK{\mathcal K}

\def \cF{\mathcal F}

\def \h{\widehat{\phantom{p}}}

\def \d{\delta_p}

\def \bZ{{\mathbb Z}}

\def \bF{{\mathbb F}}

\def \bC{{\mathbb C}}

\def \bR{{\mathbb R}}
\def \bA{{\mathbb A}}

\def \cO{\mathcal O}
\def \ra{\rightarrow}

\def \bX{{\bf X}}

\def \bF{{\mathbb F}}

\def \bQ{{\mathbb Q}}

\def \bQ{{\mathbb Q}}

\begin{document}

\title[Arithmetic PDES]{Arithmetic  partial differential
equations}
\author{Alexandru Buium and Santiago R. Simanca}
\address{University of New Mexico \\ Albuquerque, NM 87131}
\email{buium@math.unm.edu, santiago@math.unm.edu}

\begin{abstract}
We develop an arithmetic analogue of linear partial  differential
equations in two independent  ``space-time'' variables. The
spatial derivative is a Fermat quotient operator, while the time
derivative is the usual derivation. This allows us to ``flow''
integers or, more generally, points on algebraic groups with
coordinates in rings with arithmetic flavor. In particular, we
show that elliptic curves have certain canonical ``flows'' on them
that are the arithmetic analogues of the heat and wave equations.
The same is true for the additive and the multiplicative group.
\end{abstract}

\maketitle

\section{Introduction}
In this paper, we consider arithmetic partial differential
equations in two ``space-time'' variables, a higher dimensional
analogue of the theory of arithmetic ordinary differential
equations developed in \cite{char,difmod,book}. In the ordinary
case, the r\^ole of functions of one variable is played by
integers, and that of the derivative operator is played by a
``Fermat quotient operator'' with respect to a fixed prime $p$.
Instead, we now take power series in a variable $q$ with integer
coefficients as the analogues of functions of two variables, and
while maintaining the idea that a Fermat quotient type operator
with respect to $p$ is the analogue of the derivative in the
``arithmetic direction,'' we now add the usual derivative with
respect to $q$ to play the r\^ole of a derivative in the
``geometric direction.'' This leads to the study of some
``arithmetic flows'' of remarkable interest.

In the ordinary case \cite{char, difmod, book} the ``arithmetic
direction'' was viewed as a ``temporal direction.'' In the present
paper, the ``arithmetic direction'' is viewed as a ``spatial
direction,'' and the ``geometric direction'' is the ``time.''
Under this interpretation, we will be able to think of points on
algebraic varieties with coordinates in number theoretic rings as
``functions of space,'' and we will be able to flow these points
using the geometric parameter $q$, parameter that morally speaking
plays the r\^{o}le of (the exponential of) time.

We proceed to explain our idea in some detail, and begin by discussing some
of the basic aspects of evolution partial differential equations as they
appear in classical analysis, discussion that out of necessity will be carried
out in a non-rigorous fashion. In particular, the word ``function'' will be
used to refer to functions (or even distributions) belonging to
unspecified classes, and we will ignore all questions on convergence, as well
as those concerning the proper definition  of certain products or convolutions.
Instead, we will concentrate exclusively on the formal aspects of the story,
and examine only those concepts whose arithmetic analogue will later play a
r\^ole in our study. We will then describe qualitatively what these
arithmetic analogues are, and will end the introduction by a presentation
of the basic problems and results of our theory.

\subsection{Evolution equations in analysis}
We denote by ${\mathbb R}_x$ the real line with ``space''
coordinate $x$ and ${\mathbb R}_t$ the real line with
``time'' coordinate $t$. We let ${\mathcal F}({\mathbb R}_x)$ be
the ring of complex valued functions $f(x)$ on ${\mathbb R}_x$,
and we let $\cF({\mathbb R}_x \times {\mathbb R}_t)$
 be the ring of complex valued functions $v(x,t)$ on ${\mathbb R}_x \times
{\mathbb R}_t$. Both of these rings are equipped with pointwise addition
and multiplication. We will sometimes evaluate functions at
complex values of $t$ by ``analytic continuation.''

\subsubsection{Linear partial differential operators}
We consider first a general $r$-th order partial differential
operator ``in $1+1$ variables.'' This is just an operator of the
form
\begin{equation}
\label{shapteshpe}
\begin{array}{rcl}
\cF({\mathbb R}_x \times {\mathbb R}_t) & \stackrel{P}{\ra} &
\cF({\mathbb
R}_x \times {\mathbb R}_t)\\
  Pu & := & P\left(x,t,u,Du, D^2 u, \ldots, D^r u \right)\, ,
  \end{array}
\end{equation}
acting on  functions $u=u(x,t)$. In this expression, $D^n u$
stands for the $n+1$ functions $\partial^{n-k}_x
\partial_t^k u$, $0\leq k \leq n$, where $\partial_x$ and
$\partial_t$ are the corresponding partial derivative in the $x$
and $t$ directions, and $P(x,t,z)$ is a complex valued function of
$\frac{(r+1)(r+2)}{2}+2$ complex variables. The operator is said
to be {\it linear} if $P(x,t,z)$ is a linear function in the
vector variable $z$. In that case, we define the {\it full symbol}
$\sigma(P)(x,t,\xi,\tau)$ of $P$, a polynomial in $(\xi,\tau)$
with coefficients that are functions of $(x,t)$, by replacing
$\partial^{n-k}_x \partial^k_t u$ in $P$ with the monomials $i^n
\xi^{n-k}\tau^k$.  Modulo terms of degree $r-1$ or less,
$\sigma(P)$ is an invariantly defined function on the cotangent
bundle
 $T^{*}( {\mathbb R}_x \times {\mathbb R}_t)$. If $P(x,t,z)$ is
independent of $(x,t)$, we say that the operator has constant
coefficients. Standard examples of partial differential operators
that are linear and have constant coefficients are
\begin{equation} \label{ceipat}
\begin{array}{llll}
Pu & = & \partial_t u -c \partial_x u\, , &
  \textit{the convection operator}.\\
Pu & = & \partial_t u- c \partial_x^2 u \, , &  \textit{the heat operator},\\
Pu & = & \partial^2_t u-c \partial^2_x u\, , &  \textit{the wave
operator},
\end{array}
\end{equation}
respectively. Here, $c$ is a constant that in the last two cases is assumed
to be positive.
These operators are of particular importance to us as their arithmetic
analogues play a significant r\^{o}le in this article.

In regard to these examples, some remarks are in order:
\begin{enumerate}
\item If in the heat operator, we replace the real parameter $c$ by a purely
imaginary constant, then we obtain the Schr\"{o}dinger operator. As the
results of this paper will suggest, our arithmetic analogue of the heat
operator may also deserve attention as an analogue of the Schr\"{o}dinger
operator.
\item If in the heat operator we interchange $t$ and $x$, then we obtain the
{\it sideways heat operator}. This operator will also have an analogue in our
arithmetic theory.
\item If in the wave operator we replace the positive constant $c$ by a
negative one, then we obtain the Laplace operator. Similarly, if in the
convection operator $c$ is replaced by a purely imaginary constant, we
obtain the Cauchy-Riemann operator. These are the typical examples of
elliptic operators, and they do not have analogues among the
arithmetic partial differential operators discussed here.
\end{enumerate}

Given a linear partial differential operator $P$, we may consider
the {\it linear partial differential equation} $Pu=0$, and its
{\it space of solutions},
$$\cU=\cU_P:=\{u \in \cF(\bR_x \times \bR_t): \; Pu=0\}\, ,$$
that is to say, the kernel of $P$. One of the main problems of the theory is
to describe the structure of this space, and to describe the image of the
linear mapping on functions defined by $P$. Typically, a function $u$
in $\cU$ will depend  on $\rho$ arbitrary functions $C_1(\xi),\ldots,
C_{\rho}(\xi)$, $\rho$ some integer not exceeding the order $r$ of $P$,
while the description of the image of $P$ deals with the
{\it inhomogeneous equation} $Pu=\varphi$ for a given $\varphi$,
and the set of conditions on it under which the said equation
has a solution. By duality, the latter problem is usually
reduced to the consideration of the space of solutions for the dual
linear operator $P^{*}$.

\subsubsection{Exponential solutions and characteristic roots}
 Assume $P$ has constant coefficients and that $u$ is in its space of
solutions. The dependence of $u$ on the functions $C_j(\xi)$ above is easily
obtained by using Fourier transform, as follows. Let $\hat{u}(\xi,t)$ denote
the Fourier transform of $u$ in the $x$ variable. Then the equation $Pu=0$
yields
\begin{equation}
\sigma(P)(-\xi,-i \partial_t)\hat{u}(\xi,t)=0\, , \label{q0}
\end{equation}
an ordinary differential equation in the variable $t$ with parameter $\xi$.
By solving this equation and applying the inverse Fourier transform in the
parameter $\xi$, we obtain that
\begin{equation}
\label{furh} u(x,t)=\sum_{j=1}^{\rho} \int C_j(\xi) e^{-i \xi
x-i\tau_j(\xi)t} d \xi\, ,
\end{equation}
 for some functions
$C_j(\xi)$, where, for each $\xi$, the numbers $\tau_1(\xi),\ldots
,\tau_{\rho}(\xi)$ are the {\it characteristic roots} of $P$, that
is to say, the (complex) roots of the {\it characteristic
polynomial} $\sigma(P)(-\xi,-\tau)\in \bR[\tau]$, chosen to depend
continuously on $\xi$. (We ignore here the problems arising from
the possible presence of  multiple roots.)

The ``kernels'' $e^{-i \xi x-i\tau_j(\xi)t}$ in (\ref{furh}) are
the {\it exponential solutions} in the space of solutions for $P$,
and the formula exhibits the general element of this space as a
sum of $\rho$ functions. There is one exponential solution per
characteristic root, and the general solution $u$ depends
$\bC$-linearly on one arbitrary function per root. (An important
analytic aspect ignored here is that, in order to produce suitable
distributional solutions through the formal manipulations above,
we may need to choose some of the functions $C_j(\xi)$ to vanish
identically. This is dictated by the behavior of the
characteristic roots, and could make $u$ depend on fewer than
$\rho$ arbitrary functions. For instance, think of the case on an
elliptic operator, where some of the exponential solutions grow
exponentially fast.)

We observe that for any $\xi$ and $\tau$, the exponentials
$u_{\xi,\tau}(x,t):=e^{-i \xi x-i\tau t}$ ``diagonalize'' $P$.
Indeed, we have
\begin{equation} \label{pirz}
Pu_{\xi,\tau}=\sigma(P)(-\xi,- \tau) \cdot u_{\xi,\tau}\, .
\end{equation}
This fact leads naturally to the study of the inhomogeneous
equation $Pu=\varphi$ by way of Fourier inversion.

\subsubsection{Boundary value problem}
The classical approach to pinning down the functions $C_j(\xi)$
in (\ref{furh}) is by imposing ``boundary
conditions'' on the solution of $u$. For suppose we have
a $\rho$-tuple $B=(B_1, \ldots ,B_{\rho})$ of linear partial
differential operators $\cF(\bR_x \times \bR_t) \stackrel{B_j}{\ra}
 \cF(\bR_x \times \bR_t)$. We consider the {\it restriction operator}
$$\begin{array}{rcl}
\cF(\bR_x \times \bR_t) & \stackrel{\gamma}{\ra} &  \cF(\bR_x)\\
\gamma v & := & v_{|t=0}\, ,
\end{array}$$
and if we let $B^0_j$ stand for the composition $\gamma \circ B_j$, we obtain
the {\it boundary value operator}
$$\begin{array}{rcl}
\cF(\bR_x \times \bR_t) & \stackrel{B^0}{\ra} & \cF(\bR_x)^{\rho}\\
B^0 u & := & (B^0_1 u,\ldots ,B^0_{\rho} u) \, .
\end{array}$$
We then say that the {\it boundary value problem for $(P,B)$ is well
posed} if for any $g \in \cF(\bR_x)^{\rho}$ there exists a unique element
$u$ in the space of solutions ${\mathcal U}_P$ whose boundary value $B^0 u$ is
equal to $g$. In other words, the mapping
$$B^0_P:\cU_P \ra \cF(\bR_x)^{\rho} $$
given by the restriction of $B^0$ to $\cU_P$ is a $\bC$-linear
isomorphism. (Classically, the domain and range are endowed with
some topology, and the continuity of both, the mapping and its inverse, are
also required; we ignore that consideration here.) In the case where $P$ and
$B_j$ have constant coefficients, a formal computation shows that the
functions $C_j(\xi)$ and $g_j(x):=B^0 u$ are related by the
equalities
\begin{equation}
\sum_{k=1}^{\rho} \sigma(B_j)(-\xi,- \tau_k(\xi))
C_k(\xi)=\hat{g}_j(\xi),\ \ \ 1 \leq j \leq \rho\, .
\end{equation}
The determinant of the matrix of this system is the {\it
Lopatinski determinant}. Its non-vanishing is ``morally''
equivalent to the well posedness condition; cf. \cite{egorov}, pp
321-322, or \cite{hormander}.

A classical choice for the  operators $B_j$ (corresponding to the
{\it Cauchy problem}) is
\begin{equation} \label{inival}
B_j u =\partial^{j-1}_t u,\ \ \ 1 \leq j \leq \rho\, .
\end{equation}

For the classical operators $P$ listed in (\ref{ceipat}) and
the operators $B_j$ in (\ref{inival}), the corresponding boundary
value problems are all well posed.

\subsubsection{Propagator and Huygens principle}
Let us suppose that we have given operators
$P,B_1,\ldots,B_{\rho}$ such that the boundary value problem for
$(P,B)$ is well posed. Assume further that $P,B_1,\ldots
,B_{\rho}$ commute with the time translation operators
$$\begin{array}{c}
\cF(\bR_x \times \bR_t) \stackrel{L_{t_0}}{\ra}
 \cF(\bR_x \times \bR_t) \\ L_{t_0}(g(x,t))=g(x,t+t_0)
\end{array}
$$
for all $t_0$. (This is the case, for instance, if the operators
$P,B_1,\ldots ,B_{\rho}$ have constant coefficients.) Then $\cU_P$
is stable under all $L_t$, and  we have the $\bC$-linear
isomorphisms
$$B^t_{P}:=B^0_P \circ L_{t}\, : \,
\cU_P \ra \cF(\bR_x)^{\rho} \, ,$$ explicitly given by
$$B^{t_0}_{P}u=(B_1 u,\ldots ,B_{\rho} u)_{|t=t_0}\, .$$
For any pair $t_1$ and $t_2$, we obtain the {\it
evolution} or {\it propagator} operator, defined as the
$\bC$-linear isomorphism
$$S_{t_1,t_2}:=B_{P}^{t_2} \circ (B_{P}^{t_1})^{-1}=
B_P^0 \circ L_{t_2-t_1} \circ (B_P^0)^{-1}: \cF(\bR_x)^{\rho} \ra
\cF(\bR_x)^{\rho}\, .$$ This family of operators satisfies the
$1$-parameter group property
$$S_{0, t_1+t_2}=S_{0,t_1} \circ S_{0,t_2} \, ,$$
a weak form of the ``Huygens principle.''

\subsubsection{Fundamental solutions}
The idea of evolution operator above is closely
related to the concept of {\it fundamental solution}. In order to
review this concept, let $\cF_{\star}(\bR_x)$ be the Abelian group
$\cF(\bR_x)$, viewed as a ring with respect to convolution
$$(f \star g)(x):=\int f(y)g(x-y)dy\, .$$
We also let $\cF_{\star}(\bR_x \times \bR_t)$ be the Abelian group
$\cF(\bR_x \times \bR_t)$, viewed as a module over
$\cF_{\star}(\bR_x)$ with respect to convolution in the variable
$x$. Then, if $P$ commutes with the translation operators in the
variable $x$ (for instance, if $P$ has constant coefficients),
then the space of solutions $\cU_P$ is a
$\cF_{\star}(\bR_x)$-submodule of $\cF_{\star}(\bR_x \times
\bR_t)$. In this situation, we will assume further that $B^0$ is a
$\cF_{\star}(\bR_x)$-module homomorphism. This is the case for the
 operators $B_j$ in (\ref{inival}).

Under these circumstances, if the boundary value problem for
$(P,B)$  is well posed, $B_P^0$ is an
$\cF_*(\bR_x)$-module isomorphism so the space of solutions
$\cU_P$ is a free $\cF_{\star}(\bR_x)$-module of rank $\rho$ with
a unique basis
\begin{equation}
\label{fundulet} u_{1}(x,t), \ldots ,u_{\rho}(x,t)
\end{equation}
(that generally speaking, will consist of distributions) such that
the $\rho \times \rho$ matrix $(B_j^0 u_i)$ is diagonal, with
diagonal entries the Dirac delta function $\delta_0=\delta_0(x)$
centered at $0$. This basis is the {\it system of fundamental
solutions} of $(P,B)$. Clearly, for any $u \in \cU_P$, we have that
\begin{equation}
\label{uaah} u(x,t)=\sum_{i=1}^{\rho} (B_i^0 u)(x) \star u_i(x,t)\, ,
\end{equation}
and so the fundamental solutions $u_i$ appear as kernels in this
integral representation for $u$. Conversely, let us assume
that the operators $B_j$ are as in (\ref{inival}), and that we can
find $K_1,\ldots ,K_{\rho} \in \cU_P$ such that any $u \in \cU_P$
can be written uniquely as
$$u(x,t)=\sum_{i=1}^{\rho} (\partial_t^{i-1} u)(x,0) \star K_i(x,t) \, .$$
(Cf. \cite{rauch}, p. 138, for the case of the operators listed in
(\ref{ceipat}).) Applying $\partial_t^{j-1}$ to this identity, $1
\leq j \leq \rho$, and letting $t\ra 0$, we get that the matrix
$K=K(x,t)$, whose entries are given by $K_{ij}=\partial^{j-1}_t
K_i$, $1 \leq i,j \leq \rho$,  has the property that $K(x,0)$ is
equal to ${\rm diag}(\delta_0,\ldots ,\delta_0)$. Thus, we
conclude that $K_1,\ldots ,K_{\rho}$ is the system of fundamental
solutions of $(P,B)$. The matrix $K$ is the {\it fundamental
solution matrix}.

For simplicity, let us assume further that $Pu=\partial_t^{\rho} u-c
\partial_x^s u$. As $K_{\rho}$ is in the space of solutions of $P$, we
see that $\partial_t^{\rho-1} K_{\rho},
\partial_t^{\rho-2}K_{\rho}, \ldots ,\partial_t^0K_{\rho}$
 is also a system of
fundamental solutions for $(P,B)$. The uniqueness implies that
$$ K_1=\partial_t^{\rho-1}K_{\rho},\; K_2 =\partial_t^{\rho-2}K_{\rho},
\;  \ldots ,\;  K_{\rho-1}=
\partial_t K_{\rho}\, .$$
Notice that
$$S_{0,t}(g(x))=g(x) \star K(x,t)$$
for any $g(x) \in \cF(\bR_x)^{\rho}$, that is, the propagator operator
is given by convolution with the fundamental solution matrix.

 It is worth recalling that $K_{\rho}$ above may be used to solve
the {\it inhomogeneous equation} $Pu=\varphi$, as in the following discussion
where for simplicity we take once again $Pu=\partial_t^{\rho} u-c
\partial_x^s u$. Cf.
\cite{hormander}, pp. 80, 109, or \cite{egorov}, pp. 142, 235. Let
$H \in \cF(\bR_t)$ be the characteristic function of the interval
$[0, \infty)$ (the {\it Heaviside function}). If we set
$K_+(x,t):=K_{\rho}(x,t) \cdot H(t)$ where, of course, we are
implicitly assuming that the product of the distributions in the
right hand side is well defined, then we see that
$$P K_+(x,t)=\delta_0(x) \delta_0(t)\, .$$
A function of $(x,t)$ satisfying this equation is said to be a
{\it fundamental solution of the inhomogeneous equation}. A formal
computation shows  that for any $\varphi(x,t) \in \cF(\bR_x \times
\bR_t)$, the  function $K_+ \star \varphi$ (where $\star$ denotes
now the convolution with respect to both variables $x$ and $t$,)
is a solution to the equation $Pu=\varphi$, that is to say,
$$P(K_+ \star \varphi)=\varphi\, .$$

\subsection{Evolution equations in arithmetic}
The main purpose of this paper is to propose an arithmetic
analogue of the ``$1+1$ evolution picture'' above. In the remaining
portion of this introduction, we informally present our main concepts,
problems, and results.

\subsubsection{Main concepts}
In our arithmetic theory, the ring $\cF({\mathbb R}_x)$ of
functions in $x$ is replaced by a ``ring of numbers'' $R$. A
natural choice \cite{char,difmod,book} for this $R$ is given by
the completion of the maximum unramified extension of the ring
$\bZ_p$ of $p$-adic integers. We will think of $p$ as a ``space
variable,'' the analogue of $x$. And the analogue of the ring
$\cF({\mathbb R}_{t}\times {\mathbb R}_x)$ of functions of
space-time is the ring of formal power series $A=R[[q]]$, whose
elements are viewed as ``superpositions'' of the ``plane waves''
$aq^n$, $a \in R$, analogues of the plane waves $a(x)e^{-2 \pi i n
t}$ of frequency $n$. (We will use other rings also, for instance,
$R[[q^{-1}]]$. Series in $R[[q]]$ will be viewed as superpositions
of plane waves ``involving non-negative frequencies only,''
whereas series in $R[[q^{-1}]]$ will be superpositions of plane
waves ``involving non-positive frequencies only.'' It will be
interesting to further enlarge these by considering the rings of
Laurent power series and their $p$-adic completions, $R((q))\wh$
and $R((q^{-1}))\wh$, respectively.)

The r\^ole of the partial derivative $-(2 \pi i)^{-1}\partial_t$
is to be played by the derivation
\begin{equation}
\begin{array}{rcl}
A & \stackrel{\dd}{\ra} & A \\
\dd u & := & {\displaystyle q \partial_q u}
\end{array}\, ,
\end{equation}
where $\partial_q$ is the usual derivative with respect to $q$.
On the other hand, the analogue of the partial derivative $(2 \pi
i)^{-1}\partial_x$, which should be interpreted as a derivative
with respect to the prime $p$, is obtained by following the idea in
\cite{char,difmod,book}. Indeed, we propose to define this
derivative with respect to $p$ as the ``Fermat quotient operator''
given by
\begin{equation} \label{defde}
\begin{array}{rcl}
A & \stackrel{\d}{\ra} & A \\
\d u & := & {\displaystyle  \frac{u^{(\phi)}(q^p)-u(q)^p}{p}}
\end{array}\, ,
\end{equation}
where the upper index $(\phi)$ stands for the operation of twisting the
coefficients of a series by the unique automorphism $\phi: R \ra R$
that lifts the $p$-th power Frobenius automorphism of $R/pR$. Note
that the restriction of $\d$ to $R$ is the mapping
\begin{equation}
\label{defcu}
\begin{array}{rcl}
R & \stackrel{\d}{\ra} & R \\
\d a & = & {\displaystyle  \frac{\phi(a)-a^p}{p}}
\end{array}\, ,
\end{equation}
which is the arithmetic analogue of a derivation, as discussed in
\cite{char, difmod, book}. Note that the set of its {\it
constants}, $R^{\d}:=\{a \in R;\d a=0\}$ consists of $0$ and the
roots of unity in $R$. So $R^{\d}$ is a multiplicative monoid but
not a subring of $R$.

In order to proceed, we need to describe the analogues of linear partial
differential equations. We start, more generally, with maps of the
form
\begin{equation}
\label{baiie} \begin{array}{rcl} A & \stackrel{P}{\ra} &  A\\
Pu & := &P(u,Du,\ldots ,D^r u),
\end{array}
\end{equation}
where $P=P(z)$ is a  $p$-adic limit of polynomials with
coefficients in $A$ in $\frac{(r+1)(r+2)}{2}$ variables, and $D^n
u$ stands for the $n+1$ series $\d^i \dd^{n-i} u$, $0\leq i \leq
n$.
 These maps are the  ``partial differential
operators'' of order $r$ in this article.
What remains to be done is to define the notion of linearity for them.

The naive requirement that  $P(z)$ be a linear form in $z$ is not
appropriate. Indeed, the property that linearity should really
capture is that differences of solutions be again solutions, and
this is not going to happen since $\d$ itself is not additive. We
could insist upon asking that the map $u \mapsto Pu$ be additive,
but as such, this condition would lead to a rather restricted
class of examples. In order to find what we suggest is the right
concept (which, in particular, will provide sufficiently many
interesting examples), we proceed by generalizing our setting as
follows (cf. \cite{char, difmod, book} for the ordinary
differential case).

Firstly, we consider mappings $A^N \ra A$ as in (\ref{baiie}),
where now $u$ is an $N$-tuple. Secondly, we restrict such maps to
subsets $X(A) \subset \bA^N(A)=A^N$, where $X\subset {\mathbb
A}^N$ is a closed subscheme of the affine $N$-space $\bA^N$ over
$A$, and where $X(A)$ denotes the set of points of $X$ with
coordinates in $A$. If $X$ has relative dimension $n$ over $A$,
the induced maps $X(A) \ra A$ will be viewed as ``partial
differential operators'' on $X$ in $1+1$ ``independent variables''
and $n$ ``dependent variables.'' Using a gluing procedure, we then
derive the concept of a ``partial differential operator,'' $X(A)
\ra A$, on (the set of $A$-points of) an arbitrary scheme $X$ of
finite type over $A$ (which need not be affine). Finally, when we
take $X$ in this general set-up to be a commutative group scheme
$G$ over $A$, we define a {\it linear partial differential
operator} on $G$ to be a ``partial differential operator'' $G(A)
\ra A$ on $G$ that is also a group homomorphism, where $A$ is
viewed with its additive group structure. By making this
``set-theoretical'' definition one that is ``scheme-theoretical''
(varying $A$), we arrive at the notion of {\it linear partial
differential operator} that we propose in here. And once again, we
are able to associate to any linear partial differential
operator so defined a polynomial $\sigma(\xi_p,\xi_q)$ in two
variables with $A$-coefficients, which we refer to as the
({\it Picard-Fuchs}) {\it symbol} of the operator.

The construction above is motivated by points of view adopted in
analysis and mathematical physics. Indeed, we view ${\rm Spec}\, A$
as the analogue of $\bR_x \times \bR_t$, and we view schemes $X$
as the analogues of manifolds $M$, so the set $X(A)={\rm Hom}({\rm
Spec}\, A,X)$ is the analogue of the set $\cF(\bR_x \times
\bR_t,M)$ of maps $\bR_x \times \bR_t \ra M$ (which we require here
to be at least continuous). The commutative group schemes $G$ of
relative dimension $1$ (which will be the main concern of this
paper) are the analogues of Lie groups of the form $\bC/\Gamma$
where $\Gamma$ is a discrete subgroup of $\bC$. (There are $3$
cases, those corresponding to a subgroup $\Gamma$ of rank $0$, $1$, or $2$
respectively. The corresponding groups $\bC/\Gamma$ are the
additive group $\bC$, the multiplicative group $\bC^{\times}$, and
the elliptic curves $E$ over $\bC$) On the other
hand, any linear partial differential operator $ \cF(\bR_x \times
\bR_t) \ra \cF(\bR_x \times \bR_t)$ with symbol that vanishes at
$(0,0)$ induces a homomorphism
$$
 \cF(\bR_x \times \bR_t)/\Gamma \ra \cF(\bR_x \times \bR_t)\, .
$$
We can consider then the composition
\begin{equation}
\label{nuinca} \cF(\bR_x \times \bR_t,\bC/\Gamma) \simeq \cF(\bR_x
\times \bR_t)/\Gamma \ra \cF(\bR_x \times \bR_t).
\end{equation}
Our groups $G(A)$ are the analogues of the groups $\cF(\bR_x \times
\bR_t,\bC/\Gamma)$, and our linear partial differential
operators $G(A) \ra A$ are the analogues of the operators in
(\ref{nuinca}).

Given a  linear partial differential operator $G(A) \stackrel{P}{
\ra} A$ in the arithmetic setting,
 we may consider the {\it group of solutions} consisting
of those $u \in G(A)$ such that $Pu=0$. This is a subgroup of
$G(A)$. Those elements in this space that ``do not vary with
time'' will define the concept of {\it stationary solutions}.
Notice that if $G$ descends to $R$, {\it stationary} will simply
mean ``belonging to $G(R)$.''

The linear partial differential operators in the sense above will
be called $\D$-{\it characters}, and we will denote them by $\psi$
rather than $P$. This will put our terminology and notation in
line with that used in \cite{char,difmod,book}, where the case of
arithmetic ordinary differential equations was treated.

As a matter of fact, some comments on the ordinary case are in
order here. If all throughout the theory sketched above we were to
insist that the operators $u \mapsto Pu$ in (\ref{baiie}) be given
by polynomials in $u, \dd u, \ldots , \dd^r u$ only, we would then
be led to the Ritt-Kolchin theory of ``ordinary differential
equations'' with respect to $\dd$, cf.
\cite{ritt,kolchin,cassidy}. In particular, this would lead to the
notion of $\dd$-{\it character} of an algebraic group, which in
turn should be viewed as the analogue of a linear ordinary
differential operator on an algebraic group (cf. to the Kolchin
logarithmic derivative of algebraic groups defined over $R$,
\cite{kolchin, cassidy}, and the Manin maps of Abelian varieties
defined over $R[[q]]$, \cite{man,annals}).

If on the other hand we were to insist throughout the theory that
the operators $u \mapsto Pu$ in (\ref{baiie}) be given by $p$-adic
limits of polynomials in $u, \d u,\ldots ,\d^r u$  only, we would
then be led to the arithmetic analogue of the ordinary
differential equations in \cite{char,difmod,book}. In particular,
we would then arrive at the notion of a $\d$-character of a group
scheme, which is the arithmetic analogue of a linear ordinary
differential equation on a group scheme.

\subsubsection{Main problems}
 We present, in what follows, a sample of the main problems to be treated in
this paper:

\begin{enumerate}
\item[1.] Find all $\D$-characters on a given group scheme $G$.
\item[2.] For any $\D$-character $\psi$, describe the kernel of
$\psi$, that is to say, the group of solutions of $\psi u=0$, and study the
behavior of the solutions in terms of ``convolution,'' ``boundary
value problems,'' ``characteristic polynomial,'' ``propagators,''
``Huygens' principle,'' etc.
\item[3.] For any $\D$-character $\psi$, describe the image of $\psi$,
that is to say, the group of all the $\varphi$ such that the inhomogeneous
equation $\psi u=\varphi$ has a solution.
\end{enumerate}

These will be discussed in detail for the case of one dependent variable,
that is to say, for groups of dimension one. Indeed, we shall thoroughly
analyze
the additive group $G={\mathbb G}_a$, the multiplicative group
$G={\mathbb G}_m$, and elliptic curves $G=E$ over $A$, cases where $G(A)$
corresponds to the additive group $(A,+)$, the
multiplicative group $(A^{\times},\cdot)$ of invertible elements
of $A$, and the group $(E(A),+)$ of points with coordinates in $A$
of a projective non-singular cubic, with addition operation given by the
chord-tangent construction, respectively.

\subsubsection{Main results}
In regard to Problem 1, we will start by proving  that the
$\D$-characters of fixed order on a fixed group scheme form a
finitely generated $A$-module. We will then provide a rather
complete picture of the space of $\D$-characters in the cases
where $G$ is either ${\mathbb G}_a$, or ${\mathbb G}_m$, or an
elliptic curve $E$ over $A$. In particular, it will turn out that
these groups possess certain remarkable $\D$-characters  that are,
roughly speaking, the analogues of the classical operators listed
in (\ref{ceipat}) above.

In the cases where $G$ is either ${\mathbb G}_a$, or ${\mathbb
G}_m$, or an elliptic curve $E$ defined over $R$, the
$\D$-characters $\psi$ of $G$ are ``essentially'' built from
$\dd$-characters $\psi_{\ddi}$, and $\d$-characters $\psi_{\di}$
of $G$. This situation is analogous to the one in classical
analysis in  $\bR_x \times \bR_t$, where linear differential
operators are ``built'' from $\partial_x$ and $\partial_t$,
respectively. Remarkably, however, if $G$ is an elliptic curve $E$
over $A$ that is sufficiently ``general,'' then $E$ possesses a
$\D$-character $\psi^1_{\di \ddi}$ that cannot be built from
$\dd$- and $\d$-characters alone. Thus, from a global point of
view, this $\D$-character $\psi^1_{\di \ddi}$ is a ``pure partial
differential object,'' in the sense that it cannot be built from
``global ordinary linear differential objects.''

This will all unravel in the following manner. We will first prove
that for any elliptic curve $E$ over $A$, there always exists a
non-zero $\D$-character $\psi^1_{\di \ddi}$ of order
one\footnotemark. \footnotetext{This is in deep contrast with the
``ordinary case,'' where for a general $E$ over $A$, there are no
non-zero $\dd$-characters of order one \cite{man,ajm}, or where
for a general $E$ over $A$, or a general $E$ over $R$, there are
no non-zero $\d$-characters of order one \cite{char, book}. Here,
an elliptic curve $E$ over a ring is said to be {\it general} if
the coefficients of the defining cubic belong to the ring and do
not satisfy a certain ``(arithmetic) differential equation.''} We
will then show that for a general elliptic curve $E$ over $A$, the
$A$-module of $\D$-characters of order one has rank one, and so
$\psi^1_{\di \ddi}$ is essentially unique. We view it as a {\it
canonical convection equation} on the elliptic curve. Let us note
that $\psi^1_{\di \ddi}$ cannot be decomposed as a linear
combination of $\dd$- and $\d$-characters alone because, on these
elliptic curves, all such characters of order one are trivial.

The $\D$-character $\psi^1_{\di \ddi}$ will turn out to be a factor
of a canonical order two $\D$-character that can be expressed as
$\psi^2_{\ddi}+\lambda \psi^2_{\di}$, the sum of a $\dd$-character
$\psi^2_{\ddi}$ of order two, and of $\lambda$ times a $\d$-character
$\psi^2_{\di}$ of order two also, $\lambda \in A$. We will view this
$\D$-character of order two as a {\it canonical wave
equation} on the elliptic curve.

Once again, for a general elliptic curve $E$ over
$R$, we will encounter {\it heat equations} on $E$ also. These
will be sums of the form $\psi^1_{\ddi}+\lambda \psi^2_{\di}$,
where $\psi^1_{\ddi}$ is a $\dd$-character of order one,
$\psi^2_{\di}$ is a $\d$-character of order two, and $\lambda \in
R$.

In regard to Problem 2, a first remark is that the ``generic''
$\D$-characters do not admit non-stationary solutions hence a
natural question is to characterize those that admit such
solutions.  We succeed in giving such a characterization under a
mild non-degeneracy condition on the symbol $\sigma(\xi_p,\xi_q)$,
a condition that is satisfied ``generically.'' Roughly speaking,
we will show that a non-degenerate $\D$-character $\psi$ of a
group $G$ over $R$ has non-stationary solutions if, and only if,
the polynomial $\sigma(0,\xi_q)$  has an integer root. The effect
of this criterion is best explained if we consider families
$\psi_{\lambda}$ of $\D$-characters of low order (usually $1$ or
$2$) that depend linearly on a parameter $\lambda \in R$, and ask
for the values of this parameter for which $\psi_{\lambda}$
possesses non-stationary solutions. We then discover a
``quantization'' phenomenon according to which, the only values of
$\lambda$ for which this is so form a ``discrete'' set
parameterized by integers $\kappa \in \bZ$. This singles out
certain $\D$-characters as ``canonical equations'' on our groups,
and produces, for instance, a {\it canonical heat equation} on a
general elliptic curve  over $R$.

A different kind of ``quantization'' will be encountered in our
study of Tate curves with parameter $\beta q$, where $\beta \in
R$. (These curves are defined over $R((q))$ but not over $R$.) In
that case, we obtain that the canonical convection equation has
non-stationary solutions if and only if  the values of $\beta$ are
themselves ``quantized,'' that is to say, parameterized by
integers $\kappa \in \bZ$.

The criterion above on the existence of non-stationary solutions
will be a consequence of a more detailed analysis of spaces of
solutions of $\psi u=0$. In order to explain this, we take a
non-degenerate $\D$-character $\psi$ of a group $G$ over $R$, and
consider first the groups $\cU_{\pm 1}$ of solutions of $\psi u=0$
in $G(A)$ ($A=R[[q^{\pm 1}]]$) vanishing at $q^{\pm 1}=0$,
respectively. Intuitively, these are the analogues of those
solutions  in analysis that ``decay to zero'' as $t \ra \mp i
\infty$, respectively. We will prove that $\cU_{\pm 1}$ are
finitely generated free $R$-modules under a natural convolution
operation denoted by $\star$. The ranks of these $R$-modules are
given by the cardinalities $\rho_{+}$ and $\rho_{-}$ of the sets
${\mathcal K}_{\pm}$ of positive and negative integer roots of the
polynomial $\sigma(0,\xi_q)$, where $\sigma(\xi_p,\xi_q)$ is the
symbol of $\psi$.  The integers in $\cK_{\pm}$ are the {\it
characteristic integers} of $\psi$.

For instance, if $\psi$ is the arithmetic analogue of the
convection or heat equation then one of the spaces $\cU_{\pm 1}$
is zero and the other has rank one over $R$ under convolution. If
$\psi$ is the arithmetic analogue of the wave equation then both
spaces $\cU_{\pm 1}$ have rank one over $R$, an analogue of the
picture in d'Alembert's formula where one has $2$ waves traveling
in opposite directions.

Going back to the general situation of a non-degenerate
$\D$-character $\psi$, we will consider a {\it boundary value
problem at $q^{\pm 1}=0$} as follows. Firstly, we will consider
the operators
$$\begin{array}{rcl}
R[[q^{\pm 1}]] & \stackrel{\Gamma_{\kappa}}{\ra} & R \vspace{1mm} \\
\Gamma_{\kappa}(\sum a_n q^n) & :=  & a_{\kappa}=\frac{1}{\kappa
!} (\partial_{q^{\pm 1}}^{\kappa} u)_{|q^{\pm 1}=0}
\end{array}\, ,
$$
where, of course, $\partial_{q^{-1}}:=-q^2 \partial_q$. Secondly,
we note that, up to an invertible element in $R$, there is a
unique non-zero $\dd$-character $\psi_q$ of $G$ of minimal order.
The $\D$-character $\psi_q$ has order $0$ if $G=\bG_a$, and order
$1$
 if $G$ is either $\bG_m$ or an elliptic curve
over $R$; in the latter case $\psi_q$
 is the Kolchin logarithmic derivative.
 For
$\kappa \in \cK_{\pm}$ we denote by $B_{\kappa}^0$ the composition
$\Gamma_{\kappa} \circ \psi_q$. We will then prove that the {\it
boundary value operator at $q^{\pm 1}=0$},
$$\begin{array}{rcl}
\cU_{\pm 1}  & \stackrel{B_{\pm}^0}{\ra} & R^{\rho_{\pm}}\\
B_{\pm}^0 u & = & (B_{\kappa}^0 u)_{\kappa \in \cK_{\pm}}
\end{array}\, ,$$
is an $R$-module isomorphism, and, furthermore, that there exist solutions
$u_{\kappa} \in \cU_{\pm}$ such that for any $u \in \cU_{\pm 1}$
we have the formula
\begin{equation} \label{macbet} u=\sum_{\kappa \in {\mathcal
K}_{\pm}} (B_{\kappa}^0 u) \star u_{\kappa}\, .
\end{equation}
This can be viewed as analogue of the expression (\ref{furh}) because it
exhibits $u$ as a sum of $\rho_{\pm}$ terms indexed by the
characteristic integers, with  each term depending $\bZ$-linearly
on one arbitrary ``function of space.'' It can also be viewed as an
analogue of (\ref{uaah}) because of its formulation in terms of
convolution. Then we interpret the bijectivity of
$B_{\pm}^0$ by saying that the ``boundary value problem at $q^{\pm
1}=0$'' is {\it well posed}. (The language chosen here is a bit
lax since no direct analogue of this boundary value problem at $q^{\pm 1}=0$
seems to be available in real analysis; indeed, such an analogue would
prescribe boundary values at (complex) infinity, which does not
appear as a natural condition to be imposed on solutions of linear partial
differential equations in analysis.)

The solutions $u_{\kappa}$ ($\kappa \in \cK$) will be referred to
as {\it basic solutions} of the $\D$-character $\psi $. In some
sense, these elements of the kernel of $\psi$ are the analogues of
both the exponential solutions and the fundamental solutions of
the homogeneous equations in real analysis. Of course, these
analogies have significant limitations.

An interesting feature of the solutions $u$ in the kernel of
$\psi$ is the following ``algebraicity modulo $p$''
 property. Let us denote by $k$ the residue field of $R$,
hence $k=R/pR$. Then, for any such $u$,  the reduction modulo $p$,
$\overline{\psi_q u} \in k[[q^{\pm 1}]]$,
 of the series $\psi_q u$ is integral over the
polynomial ring $k[q^{\pm 1}]$, and the field extension
$$k(q) \subset k(q,\overline{\psi_q u})$$
is Abelian with Galois group killed by $p$. On the other hand,
under certain general assumptions that are satisfied, in particular, by
our ``canonical'' equations, we will prove that the
solutions $u \neq 0$ of $\psi u=0$ are transcendental over $R[q]$.
This transcendence result can be viewed as a (weak) incarnation
of Manin's Theorem of the Kernel \cite{man}.

Some of the results above about groups $G$ over $R$ have analogues
for groups not defined over $R$, such as the Tate curves. For the
latter the $\dd$-character $\psi_q$ will now be the Manin map,
which has order $2$.

If we instead consider solutions that do not necessarily decay to
$0$ as $q^{\pm 1} \ra 0$, we are able to show that, in case
$\rho_+=1$,  the appropriate boundary value  problem at $q \neq
0$ is well posed. (The condition $\rho_+=1$ is usually satisfied
by our ``canonical equations''.) Roughly speaking, for groups $G$
over $R$ and $\D$-characters $\psi$ with $\rho_+=1$, this boundary
value problem has the following meaning. We consider $q_0 \in
pR^{\times}$ and $g \in G(R)$. If $A=R[[q]]$, we ask if there exists a
(possibly unique) solution $u \in G(A)$ of the equation $\psi
u=0$ that satisfies the condition $B^{q_0}u =g$, where
$$\begin{array}{rcl}
G(A) & \stackrel{B^{q_0}}{\ra} &  G(R) \\
B^{q_0} u & := & u(q_0) \end{array}
$$
is the group homomorphism induced by the ring homomorphism
$A=R[[q]] \ra R$ that sends $q$ into $q_0$. We view $B^{q_0}$ as
a {\it boundary value operator at $q_0$}, and we view the
corresponding embedding ${\rm Spec}\, R \ra {\rm Spec}\, A$  as
the ``curve'' $q=q_0$ in the $(p,q)$-plane, along which we are
given our boundary values.

Our point of view here is reminiscent of that in real analysis,
where we replace real time by a complex number whose real part is
small in order to avoid singularities on the real axis. Indeed the
condition ``$q_0$ invertible'' (that is to say, $v_p(q_0)=0$,
where $v_p$ is the $p$-adic valuation)
 is an
analogue of the condition ``real time;'' the condition ``$q_0$
non-invertible with small $v_p(q_0)$'' is an analogue of ``complex
time close to the real axis;'' ``$q_0$ close to $0$'' (that is,
$v_p(q_0)$ big)  is an analogue of ``time close to $-i \infty$.''

For  equations with $\rho_+=1$, we will investigate also the
``limit of solutions as $q \rightarrow 0$'' (or intuitively, as $t
\rightarrow -i \infty$). This limit will sometimes exist, and if
so, the limit will usually be a torsion point of $G(R)$. In this
sense, torsion points tend to play the role of ``equilibrium
states at (complex) infinity.''

These principles do not hold uniformly in all examples. For
instance, in the case of elliptic curves $E$ over $R$, we will
need  to replace $E(R)$ by a suitable subgroup of it, $E'(R)$, in
order to avoid points whose orders are powers of $p$. And for
elliptic curves over $A$, rather than those over $R$ (such as the
Tate curves), the boundary value problem at $q \neq 0$ will take a
slightly different form.

Once the boundary value problem at $q \neq 0$ has been solved, we
can construct  {\it propagators} as follows. Again, we fix $\psi$
with $\rho_+=1$, choose ``complex times'' $q_1, q_2 \in
pR^{\times}$, and consider the endomorphism $S_{q_1,q_2}$ of the
group $G(R)$ that sends any $g_1 \in G(R)$ into $g_2:=u(q_2)$,
where $u \in G(A)$ is the unique solution to the boundary value
problem
$$\begin{array}{ccl} \psi u & = & 0\\
u(q_1) & = & g_1
\end{array}.$$ This construction does not work
uniformly in all examples the way it is described here. In order
to make it so, we need an appropriate modification of the given
recipe. Nevertheless, in all situations,  we find that the
propagator $S$ satisfies a $1$-parameter group property. For it
turns out that given ``complex times'' $q_i=\zeta_i q_0$ with
$\zeta_i \in R$ a root of unity, $i=1,2$, we have that
$$
S_{q_0,\zeta_1 \zeta_2 q_0}=S_{q_0,\zeta_2 q_0} \circ
S_{q_0,\zeta_1 q_0}\, ,
$$
identity that can be viewed as a weak incarnation
of Huygens' principle.

In regard to Problem 3, we will consider non-degenerate
$\D$-characters $\psi$ of $G$, and we will find sufficient
conditions on a given series $\varphi  \in
A=R[[q]]$ ensuring that the inhomogeneous equation $\psi
u=\varphi$ have a solution $u \in G(A)$.
Specifically, let us define the support of the series $\varphi=\sum c_nq^n$
as the set $\{ n\, : \; c_n \neq 0\}$. We will then prove that if $\varphi$
has support contained in the set $\cK'$ of {\it totally
non-characteristic} integers, then the inhomogeneous equation
above, with $\varphi$ as right hand side, has a unique solution
$u \in G(A)$ for which the support of $\psi_q u$ is disjoint from the set
$\cK$ of characteristic integers. Here, $\cK'$ is defined by
 an easy congruence involving the symbol, and, as suggested by the terminology,
 $\cK'$ is disjoint from $\cK$. Furthermore, if $\bar{\varphi}$, the reduction
mod $p$ of $\varphi$, is a polynomial, then we will prove a
corresponding ``algebraicity mod $p$'' property for $u$ stating
that the series $\overline{\psi_q u}$ is integral over $k[q^{\pm
1}]$, and that the field extension $k(q)\subset k(q,\overline{\psi_q
u})$ is Abelian with Galois group killed by $p$. On the other
hand, for a ``canonical'' $\psi$ and for a $\varphi$ with {\it short}
support, we will show that the solutions of $\psi u=\varphi$ are
transcendental over $R[q]$.

The idea behind the results above is to construct, for all integers $\kappa$
coprime to $p$, a solution $u_{\kappa}$ of the equation
\begin{equation}
\label{obbossesc} \psi u_{\kappa}=\frac{\sigma(0,\kappa)}{p} \cdot
q^{\kappa} \, ,
\end{equation}
which we shall call a {\it basic} solution of the inhomogeneous
equation. For $\kappa$ a characteristic integer, the right hand
side of (\ref{obbossesc}) vanishes and our basic solutions are the
previously mentioned basic solutions of the homogeneous equation.
For $\kappa$ a totally non-characteristic integer, the right hand
side of (\ref{obbossesc}) is a unit in $R$ times $q^{\kappa}$, and
that leads to the desired result about the inhomogeneous equation.
Notice that (\ref{obbossesc}) should be viewed as an analogue of
(\ref{pirz}), and the $u_{\kappa}$'s (for $\kappa \in \bZ
\backslash p\bZ$) should be viewed as (partially) diagonalizing
$\psi$.

\subsection{Concluding remarks}
It is natural to ask for an extension of the theory in the present
paper to the case of $d+e$ independent variables and $n$ dependent
variables, that is to say, to the case of $d$ time variables $q_1,
\ldots ,q_d$, $e$ space variables $p_1, \ldots ,p_e$, where $p_i$
are prime numbers, and groups $G$ of (relative) dimension $n$. It
is not hard to perform such an extension to the  case $e=1$, $d
\geq 1$, $n \geq 1$ (that is to say, one prime $p$ as space
variable, $d \geq 1$ indeterminates $q_i$ as time variables, and
groups of dimension $n \geq 1$). In fact, all the difficulties of
this more general case are already present in ours here, where
$d=e=n=1$. On the other hand, there seems to be no obvious way to
extend the theory in a non-trivial way even to the case $d=1$,
$e > 1$, $n=1$ (that is to say, two or more primes $p_i$ as space
variables, one time variable $q$, and groups of dimension $1$).
The main obstruction lies in the fact that when at least two primes are
made to interact in the same equation, the solutions exhibit a rather
fundamental divergent form.

We end our discussion here by summarizing (cf. the tables below)
some of the similarities and differences between the set-ups of
classical real analysis \cite{hormander, folland, rauch},
classical $p$-adic analysis \cite{dwork, koblitz}, and arithmetic
(in the spirit of  \cite{char, difmod, book} for the ordinary differential
case, and the present paper for the partial differential case).
 For the ordinary differential case we have:

\bigskip

\begin{center}
\begin{tabular}{||c|c|c|c||} \hline \hline
$\ $  & real analysis  &  $p$-adic analysis & arithmetic \\ \hline
$1$-dimensional manifold & ${\mathbb R}_{x}$ & $R$ & $R^{\d}$ \\
\hline
ring of functions & $\cF({\mathbb R}_{x})$ & $R[[x]]$ & $R$ \\
\hline
operator on functions & $\partial_x$ & $\partial_x$ & $\d$ \\
\hline
\end{tabular}
\end{center}

\bigskip

\noindent For the partial differential case we have:

\bigskip

\begin{center}
\begin{tabular}{||c|c|c|c||} \hline \hline
  & real analysis  &  $p$-adic analysis & arithmetic \\ \hline
  $2$-dimensional
manifold & ${\mathbb R}_{t}\times {\mathbb R}_{x}$ & $R\times R$ &
$R^{\d} \times R$
\\ \hline
ring of functions & $\cF({\mathbb R}_{x}\times {\mathbb R}_{t})$ &
$R[[x,t]]$ & $R[[q]]$
\\ \hline
operators on functions & $\partial_x,\partial_t$ &
$\partial_x,\partial_t$ & $\d, \dd$ \\ \hline
\end{tabular}
\end{center}
\bigskip

\noindent In the last column,  the set $R^{\d}$  plays a role
similar to that of the set of ``geometric points of the spectrum
of the field with one element'' in the sense of Deninger,
Kurokawa, Manin, Soul\'{e}, and others \cite{kurokawa, manin,
soule}. (For comments on differences between our approach here and
the ideology of the ``field with one element,'' we refer to the
Introduction of \cite{book}.) In particular $R^{\d}$ can be viewed
as an object of dimension zero. Notice that the third column
appears to be obtained from the second one by ``decreasing
dimensions by one;'' this reflects the fact that, in contrast to
the second column, the third one treats numbers as functions. Note
also that the interpretation of $R^{\d}$ as a space of dimension
zero is morally consistent with the (otherwise) curious fact that
the groups of solutions $\cU_{1}$ and $\cU_{-1}$ are modules over
$R$ with respect to a ``convolution'' operation $\star$. Indeed,
this suggests that ``pointwise'' multiplication and
``convolution'' of functions on $R^{\d}$ coincide, which in turn,
is compatible with viewing
 $R^{\d}$
as having dimension $0$.

\subsection{Plan of the paper} We begin in \S 2 by
introducing our main concepts, and where, in particular, we define
$\D$-characters, their various spaces of solutions, and the
convolution module structure on these spaces. In \S 3 and \S 4,
we construct and study partial differential jet spaces of schemes
and formal groups, respectively. These are arithmetic-geometric
analogues of the standard jet spaces of differential geometry,
arithmetic in the space direction, and geometric in the
time direction. The geometry of these jet spaces controls the
structure of the spaces of $\D$-characters. Among several other
constructions, we define in \S 4 the Picard-Fuchs
symbol of a $\D$-character. In \S 5, we develop an analogue of
the Fr\'{e}chet derivative, and use it to define the
Fr\'{e}chet symbol of a $\D$-character. We then establish a link
between the Fr\'{e}chet symbol and the Picard-Fuchs symbol that will
be useful in applications. We also develop in this section a brief analogue
of the Euler-Lagrange formalism. In \S 6, \S 7 and \S 8, we
present our main results about $\D$-characters and their solution spaces in
the cases of $\bG_a$, $\bG_m$, and elliptic curves, respectively.

\bigskip

\noindent {\bf Acknowledgements.} We would like to
 acknowledge useful conversations
with F. J. Voloch and D. Thakur on the subject of transcendence.

\section{Main concepts}
\setcounter{theorem}{0} In this section, we begin by introducing
the main algebraic concepts in our study, especially the
$\D$-rings and $\D$-prolongation sequences. These are then used to
define $\D$-morphisms of schemes, and, eventually, $\D$-characters
of group schemes. We introduce various solution spaces of
$\D$-characters, and discuss the convolution module structure on
them.

 Let $p$ be a prime integer that we fix
throughout the entire paper. For technical purposes related to the
use of logarithms of formal groups, we need to assume that $p \neq
2$. Later on, in our applications to elliptic curves, we will need
to assume that $p\neq 3$ also. All throughout, $A$
shall be a ring, and $B$ an algebra over $A$. If $x$ is an element of $A$,
we shall denote by $x$ its image in $B$ also.

We let $C_p(X,Y)$ stand for
$$C_p(X,Y):=\frac{X^p+Y^p-(X+Y)^p}{p}\, ,$$
an element of the polynomial ring ${\mathbb Z}[X,Y]$.

After \cite{char}, we say that a map
$\delta:  A \rightarrow B$ is a $p$-{\it derivation} if $\d(1)=0$, and
$$\begin{array}{rcl}
\d(x+y) & = &  \d x + \d y +C_p(x,y)\, , \\
\d(xy)  & = & x^p  \d y +y^p  \d x  +p  \d x  \d y\, ,
\end{array}$$
for all $x,y \in A$, respectively. We will always write
$\d x$ instead of $\d(x)$.

Given a $\d$-derivation, we define $\phi:=\phip:=
\phi_{\d}: A \rightarrow B$ by
\begin{equation}
\phip(x)= x^p +p\d x \, , \label{me0}
\end{equation}
a map that turns out to be a homomorphism of rings. Sometimes, we will write
$x^{\phi}$ instead of $\phip(x)$, and when omitted from the notation, the
context will indicate the $p$-derivation $\d$ that is being used.

We recall that a map $\dd : A\rightarrow B$ is said to be
a {\it derivation} if
$$\begin{array}{rcl}
\dd(x+y) & =  & \dd x + \dd y\, , \\
\dd(xy) & = & x\dd y+y\dd x\, ,
\end{array}$$
for all $x,y \in A$, respectively.
For the time being, the index $q$ will not be given any interpretation.
Later on, we will encounter the situation where $q$ is an element
of $A$, and in that case, we will think of
$\dd$ as a derivation in the ``direction'' $q$.

\begin{definition}
Let $\d:A \ra B$ be a $p$-derivation. A derivation $\dd : A
\rightarrow B$ is said to be a  $\d$-{\it derivation} if
\begin{equation}
\label{didi} \dd \d x=p\d \dd x + (\dd x)^p-x^{p-1}\dd x
\end{equation}
for all $x \in A$.
\end{definition}

In particular, for a $\d$-derivation  $\dd : A \rightarrow
B$, we have that
\begin{equation}
\label{cucu} \dd \circ \phip =p \cdot \phip \circ \dd\, .
\end{equation}
Conversely, if $p$ is a non-zero divisor in $A$, then (\ref{cucu}) implies
(\ref{didi}).

\begin{definition}
A $\{\d,\dd\}$-{\it ring} $A$ is one equipped with a
$p$-derivation $\d:A \ra A$ and a $\d$-derivation $\dd : A \ra A$.
A {\it morphisms of $\{\d,\dd\}$-rings} is a ring homomorphism
that commutes with $\d$ and $\dd$. A $\{\d, \dd\}$-ring $B$ is
said to be a $\{\d,\dd\}$-{\it ring over the $\{\d,\dd\}$-ring}
$A$ if it comes equipped with a $\{\d,\dd\}$-ring homomorphism
$A\rightarrow B$. We say that a $\{\d,\dd\}$-ring $A$ is a
$\D$-{\it subring} of the $\{\d,\dd\}$-ring $B$ if $A$ is a
subring of $B$ such that $\d A \subset A$ and $\dd A \subset A$,
respectively.
\end{definition}

We describe next the basic examples of $\{\d,\dd\}$-rings that will play
crucial r\^oles in our paper. All throughout, we shall let $R$ stand for
$R:=R_p:=\widehat{\mathbb Z}_p^{ur}$, the completion of the maximum
unramified extension of ${\mathbb Z}_p$, $k$ for the residue field $k:=R/pR$,
and $K$ for the fraction field $K:=R[1/p]$. Furthermore, we will let
$\mu(R)$ be the multiplicative group of roots of unity in $R$, and recall
that the reduction mod $p$ mapping
$$\mu (R) \rightarrow k^{\times}$$ defines
an isomorphism whose inverse is the {\it Teichm\"uller lift}. Any element of
the ring $R$ can be represented uniquely as a series $\sum_{i=0}^{\infty}
\zeta_i p^i $, where $\zeta_i \in \mu(R)\cup \{ 0\}$.
There is a unique ring isomorphism
\begin{equation}
\phi : R \rightarrow R \label{fr}
\end{equation}
that lifts the $p$-th power Frobenius isomorphism on $k$, and for
$\zeta  \in \mu(R)$, we have that $\phi(\zeta)=\zeta^p$.

The ring $R$ is isomorphic to the Witt ring on the algebraic closure
${\mathbb F}_p^a$ of ${\mathbb F}_p$, and for each $s\geq 1$, the ring
$$R^{\phi^s}:=\{ x\in R \, \mid \; x^{\phi^s}=x\}$$
is isomorphic to the Witt ring on the field ${\mathbb F}_{p^s}$ with $p^s$
elements. Notice that the ring $R^{\phi}={\mathbb Z}_p$ is simply the
ring of $p$-adic integers. As usual, we denote by ${\mathbb Z}_{(p)}$ the
ring of all fractions $a/b\in {\mathbb Q}$, where $a\in {\mathbb Z}$ and
$b\in {\mathbb Z}\setminus (p)$. We have the inclusions
$${\mathbb Z}_{(p)} \subset {\mathbb Z}_{p} \subset
R^{\phi^s}\subset R \, .$$

\begin{example} \label{unex}
The ring $R$ carries a
unique $p$-derivation $\d: R \ra R$ given by (see (\ref{me0}))
$$\d x=\frac{\phip(x)-x^p}{p}\, .$$
The {\it constants of} $\d$ are defined to be
$$R^{\d}:=\{x\in R\, : \; \d x =0\}\, ,$$
set that coincides with $\mu(R)\cup \{0\}$, the roots of unity in $R$
together with $0$. Notice that we have the trivial $\d$-derivation
$\dd =0$ on $R$, and that the pair $(\d, 0)$ equips $R$ with a
$\D$-ring structure.

We denote by $R[[q]]$ and $R[[q^{-1}]]$ the power series rings in
$q$ and $q^{-1}$, respectively, and we embed them into the rings
$$
R((q))^{\wh } =R[[q]][q^{-1}]^{\wh }=
\left\{ \sum_{n=-\infty}^{\infty} a_nq^n
\, | \; \lim_{n\rightarrow -\infty} a_n =0 \right\}
$$
and
$$
R((q^{-1}))^{\wh}=R[[q^{-1}]][q]^{\wh}=
\left\{ \sum_{n=-\infty}^{\infty} a_nq^n \, | \;
\lim_{n\rightarrow \infty} a_n =0 \right\}\, ,
$$
respectively.

The rings $R((q))^{\wh}$ and $R((q^{-1}))^{\wh}$ have
unique structures of $\D$-rings that extend that of $R$, such that
$$
\d q=0\, ,\; \dd q=q\, , \; \d(q^{-1})=0\, , \; \dd(q^{-1})=-q^{-1}\, .
$$
The rings $R[[q]]$ and $R[[q^{-1}]]$ are $\D$-subrings of
$R((q))^{\wh}$ and $R((q^{-1}))^{\wh}$, respectively.
Also the
ring
$$
R[q,q^{-1}]^{\wh}= \left\{ \sum_{n=-\infty}^{\infty} a_nq^n \, | \;
\lim_{n\rightarrow \pm \infty} a_n =0 \right\}\,
$$
is a $\D$-subring of both
$R((q))^{\wh}$ and $R((q^{-1}))^{\wh}$.

We describe these $\D$-structures in further detail.

The automorphism (\ref{fr}) extends to a unique homomorphism
$\phi: R[[q]]\ra R[[q]]$ such that $\phi(q)=q^p$. Similarly, it
extends to a unique homomorphism $\phi: R[[q^{-1}]]\ra
R[[q^{-1}]]$ such that $\phi(q^{-1})=q^{-p}$, to which for obvious
reasons we give the same name. Then the expression
$$
\d F:=\frac{F^{\phi}-F^p}{p}\, .
$$
defines a $p$-derivation $\d$ on both, $R[[q]]$ and $R[[q^{-1}]]$,
respectively. On the other hand, the
expression
$$\dd F:=q \frac{dF}{dq}=q\partial_q F$$
defines a $\d$-derivation on $R((q))^{\wh}$ and
$R((q^{-1}))^{\wh}$. Here,
$$\partial_q \left( \sum a_n q^n \right)=\frac{d}{dq}
\left( \sum a_n q^n \right):= \sum na_n q^n \, .$$
The operators $\d$ and $\dd$ so defined
provide the various rings above with their respective $\D$-ring structures.

Observe that the $R$-algebra mapping
$$
\begin{array}{ccc}
R((q))^{\wh} & \rightarrow & R((q^{-1}))^{\wh} \\
q & \mapsto & q^{-1}
\end{array}
$$
is a ring isomorphism that fails to be a $\D$-ring isomorphism.
\qed
\end{example}
\medskip

In the sequel, we will repeatedly use the following ``Dwork's Lemma,''
which we record here for convenience.

\begin{lemma}\label{dwor}
Let $v \in 1+qK[[q]]$ be such that $v^{\phi}/v^p \in 1+pqR[[q]]$. Then
$v\in 1+qR[[q]]$.
\end{lemma}

{\it Proof}. For $v \in 1+q{\mathbb Q}_p[[q]]$, this is proved
in \cite{koblitz},
p. 93. The general case follows by a similar argument. \qed

We will also need the following case of Hazewinkel's
Functional Equation Lemma. This special version of this result implies
Dwork's Lemma.

\begin{lemma}\label{haze}
\cite{haze} Let $\mu_0 \in R^{\times}$, and $\mu_1,\ldots ,\mu_s \in
R$. Consider two series $f,g \in K[[q]]$ such that $f \equiv
\lambda q$ mod $q^2$ in $K[[q]]$, for some $\lambda \in
R^{\times}$. Let $f^{-1} \in qK[[q]]$ be the compositional inverse
of $f$. For
$$
\Lambda:=\mu_0+\frac{\mu_1}{p} \phip+\cdots
+\frac{\mu_s}{p}\phip^s \in K[\phip]\, ,
$$
assume that both, $\Lambda f$ and $\Lambda g$, belong to $R[[q]]$.
Then $f^{-1} \circ g \in R[[q]]$.
\end{lemma}

Closely related to these results is the following.

\begin{lemma}
\label{floricica} Let $\mu_0 \in R^{\times}$, $\mu_1,\ldots ,\mu_s
\in R$, and $f \in K[[q]]$. If for
$$\Lambda:=\mu_0+\mu_1\phip+\cdots +\mu_s \phip^s \in R[\phip]$$
we have that $\Lambda f \in qR[[q]]$, then $f \in qR[[q]]$. If on the
other hand we have that $\Lambda f \in pqR[[q]]$, then
$f \in pqR[[q]]$.
\end{lemma}

{\it Proof}. We may assume $\mu_0=1$. Set $\Lambda_1:=1-\Lambda \in \phip
R[\phip]$, and note that $\Lambda_1 q R[[q]] \subset q^pR[[q]]$.
Set $g=\Lambda f$. Then
$$f=g+\Lambda_1 g +\Lambda_1^2 g +\cdots \, ,$$
and the result follows. \qed

Formally speaking, the $\D$-rings $R((q))\wh$ and
$R((q^{-1}))\wh$, and their various $\D$-subrings introduced
above, can be viewed as rings of ``Fourier series.'' These rings
allow us to take limits of solutions to arithmetic partial
differential equations in the ``time'' direction, as $t\rightarrow
\pm i\infty$. On the other hand, the examples of $\D$-rings of
``Iwasawa series,'' which we discuss next, control limits ``as
time goes to $0$.''

\begin{example}\label{ninna}
Let $R[[q-1]]$ be the completion of the polynomial ring $R[q]$ with respect
to the ideal $(q-1)$. Hence, $R[[q-1]]$ can be identified with the power
series ring $R[[\tau]]$ in an indeterminate $\tau$, and $R[q]\wh$ embeds into
$R[[q-1]]=R[[\tau]]$ via the map $q \mapsto 1+\tau$.

There is a unique $\D$-ring structure on $R[[q-1]]$ that extends that of
$R[q]\wh$. Indeed, in terms of $\tau$ we have
$$
\begin{array}{rcl}
\d \tau & = & {\displaystyle \frac{(1 + \tau)^p -1 -\tau^p}{p}}\, ,
\vspace{1mm} \\
\dd \tau & = &  1+\tau\, .
\end{array}
$$
Note that although $R[[q]]$ and $R[[q-1]]$ are isomorphic as rings, they are
not isomorphic as $\D$-rings. In fact, they are quite ``different'' in this
latter context.

Similarly, let $R[[q^{-1}-1]]$ be the completion of $R[q^{-1}]$
with respect to the ideal $(q^{-1}-1)$. Then there is a natural
embedding $R[q^{-1}]\wh \subset R[[q^{-1}-1]]$, and a unique
structure of $\D$-ring on $R[[q^{-1}-1]]$ that extends that of
 $R[q^{-1}]\wh$. \qed
\end{example}
\medskip

From here on, we consider the ring $A$ provided with a fixed
structure of $\D$- ring. For simplicity, we always assume that $A$
is a $p$-adically complete Noetherian integral domain of
characteristic zero, and we shall let $L$ be its fraction field.
Example \ref{unex} above illustrates the cases where $A$ is equal
to $R$, $R[[q]]$, $R((q))^{\wh}$, $R[[q^{-1}]]$,
$R((q^{-1}))^{\wh}$, and
 $R[q,q^{-1}]^{\wh}$, respectively, while Example \ref{ninna} illustrates the
cases where $A$ is equal to $R[[q-1]]$ and $R[[q^{-1}-1]]$, respectively.
Whenever any of these rings is considered below, it will be given
the $\D$-ring structure defined in these two examples.

In Examples \ref{unneg} and \ref{ex2} below we explain some basic
general constructions that can be performed using a fixed
$\D$-ring $A$.

\begin{example}
\label{unneg} We will denote by $(A[\phip,\dd],+, \, \cdot \, )$
the ring generated by $A$ and the symbols $\phip$ and $\dd$
subject to the relations
$$\begin{array}{rcl}
\dd \cdot \phip & =  & p\cdot  \phip \cdot \dd\, , \\
 \mbox{} [\dd,a]=\dd \cdot a -a\cdot \dd & = &  \dd a\, , \\
\phip \cdot a & = & a^{\phi}\cdot \phip \, ,
\end{array}
$$
for all $a\in A$. Let $\xi_p$, $\xi_q$ be two variables, which we view
as ``duals to $p$ and $q$,'' respectively. If $\mu=\mu(\xi_p,\xi_q)=\sum
\mu_{ij} \xi_p^i \xi_q^j \in A[\xi_p,\xi_q]$ is a polynomial,
then we define
$$\mu(\phip,\dd):=\sum \mu_{ij} \phip^i \dd^j \in
A[\phip,\dd ] \, .
$$
The map
$$\begin{array}{rcl}
A[\xi_p,\xi_q] & \ra &  A[\phip,\dd] \\
\mu(\xi_p,\xi_q) & \mapsto & \mu(\phip,\dd)
\end{array}
$$ is a left $A$-module isomorphism.

Given a $\mu(\xi_p,\xi_q)$ as above,
we define the polynomial
$$\mu^{(p)}(\xi_p,\xi_q):=\mu(p\xi_p,\xi_q)\, .$$
We have the following useful formula:
$$\dd \mu(\phip,\dd)=\mu^{(p)}(\phip,\dd) \dd\, .$$ \qed
\end{example}

\begin{example} \label{ex2}
Let $y$ be an $N$-tuple of indeterminates over $A$, and let $y^{(i,j)}$ be
an $N$-tuple of indeterminates over $A$ parameterized by non-negative
integers $i$, $j$, such that $y^{(0,0)}=y$. We set
$$A\{ y\}:=A[y^{(i,j)}\mid_{i\geq 0, j\geq 0}]$$
for the polynomial ring in the indeterminates $y^{(i,j)}$. This ring has
a natural $\D$-structure that we now discuss.

For let $\phip : A\{ y\} \rightarrow A\{ y\}$ be the unique ring homomorphism
extending $\phip : A \rightarrow A$, and satisfying the relation
$$\phip (y^{(i,j)})= (y^{(i,j)})^p + p y^{(i+1,j)} \, .$$
Then, we may define a $p$-derivation $\d: A\{ y\} \rightarrow A\{ y\}$ by the
expression (see (\ref{me0}) above)
 $$\d F: = \frac{\phip(F)-F^{p}}{p} \, .$$
In particular, we have that $\d y^{(i,j)}=y^{(i+1,j)}$.

We let $L\{y\}$ be $L\{y\}:=L\otimes_A A\{ y\}$, $L$ the field of
fractions of $A$. Notice that
$$L\{y\}=L[y^{(i,j)}\mid_{i\geq 0, j\geq 0}]=L[\phip^{i}y^{(0,j)}
\mid_{i\geq 0,
j\geq 0}]\, .$$
The endomorphism $\phip$ extends uniquely to an endomorphism $\phip$ of
$L\{y\}$.
Since the polynomials $\phi^{i}y^{(0,j)}$ are algebraically independent over
$L$, there exists a unique derivation $\dd: L\{ y\} \rightarrow L\{ y\}$
that
extends the derivation $\dd : A \rightarrow A$, and satisfies the relation
$$\dd (\phip^{i}y^{(0,j)})=p^i \phip^{i}y^{(0,j+1)}\, .$$
We claim that $D:=\dd \circ \phip -p\cdot \phip \circ \dd$ vanishes
on $L\{ y\}$.

Indeed, $D: L\{ y\} \rightarrow L\{ y\}$ is a derivation if we view $L\{ y\}$
as an algebra over itself via $\phip$. Thus, it suffices to observe that $D$
vanishes on $L$ and $\phip^{i}y^{(0,j)}$, and both of these are clear.

We now claim also that $\dd y^{(i,j)}\in A\{y\}$, and this will imply that
$\dd$ induces a derivation of the ring $A\{ y\}$. For, proceeding by
induction on the index $i$, with the case $i=0$ being clear, we assume it for
$i$, and prove the desired statement for $i+1$. We have
$$\begin{array}{rcl}
\dd y^{(i+1,j)} & = & {\displaystyle \dd \left( \frac{\phip(y^{(i,j)})-
(y^{(i,j)})^p}{p}\right)}\vspace{1mm} \\ & = &
\phip(\dd y^{(i,j)})-(y^{(i,j)})^{p-1}
\dd y^{(i,j)} \, ,
\end{array}
$$
and this is clearly in $A\{ y\}$ by the induction hypothesis.

The morphisms $\d$, $\dd $ endow $A\{ y\}$ with a $\D$-ring structure.
Clearly, $y^{(i,j)}=\d^i \dd^j y$, and so
$$A\{ y\} = A[y, Dy, D^2y, \ldots  ]\, ,$$
where for any whole number $n$, $D^ny$ stands for the $(n+1)$-tuple with
components $\d^i \dd^{n-1}y$, $0\leq i\leq n$. \qed
\end{example}

\begin{definition}
Let $S^{*}=\{S^n\}_{n\geq 0}$ be a sequence of rings. Suppose we
have ring homomorphisms $\varphi: S^n \ra S^{n+1}$,
$p$-derivations $\d:S^n \ra S^{n+1}$, and $\d$-derivations $\dd
:S^n \ra S^{n+1}$ such that $\d \circ \varphi=\varphi \circ \d$,
and $\dd \circ \varphi=\varphi \circ \dd$. We then say that
$(S^{*}, \varphi,\d,\dd)$, or simply $S^{*}$, is a
$\{\d,\dd\}$-{\it prolongation sequence}. A {\it morphism}
$u^*:S^* \ra \tilde{S}^*$ of $\D$-prolongation sequences is a
sequence $u^n:S^n \ra \tilde{S}^n$ of ring homomorphisms such that
$\d \circ u^n=u^{n+1} \circ \d$, $\dd \circ u^n=u^{n+1} \circ
\dd$, and $\varphi \circ u^n=u^{n+1} \circ \varphi$.
\end{definition}

Let us consider the $\D$-ring structure on the ring $A$. We obtain a natural
$\D$-prolongation sequence $A^{*}$ by setting $A^n=A$ for all $n$, and
taking the ring homomorphisms $\varphi$ to be all equal to the identity.
This leads
to the natural concept of morphisms of $\D$-prolongation
sequences over $A$.

\begin{definition}
We say that a $\D$-prolongation sequence $S^*$ is a
{\it $\D$-prolongation sequence over} $A$ if we have
a morphism of prolongation sequences
$A^* \ra S^*$.
\end{definition}

There is a natural notion of morphisms of $\D$-prolongation sequences over
$A$ that we do not explicitly state. In the sequel, all
$\D$-prolongation sequences, and morphisms of such, will be prolongation
sequences over $A$.

\begin{example}
Let $A\{ y\}$ the ring discussed in Example \ref{ex2}. We consider the
subrings
$$S^n:=A[y, Dy, \ldots, D^n y ]\, .$$
We view $S^{n+1}$ as an $S^n$-algebra via the inclusion homomorphism, and
observe that $\d S^n \subset S^{n+1}$, and $\dd S^n \subset S^{n+1}$,
respectively. Therefore, $S^{*}=\{ S^n\}$ defines a $\D$-prolongation
sequence. We then obtain the $p$-adic completion prolongation sequence
$$A[y, Dy, \ldots, D^ny ]^{\widehat{\mbox{\phantom{p}}}}
\, ,$$
and the prolongation sequence
$$A[[Dy, \ldots, D^n y ]]$$
of formal power series ring, with their corresponding $\D$-structures. \qed
\end{example}

\begin{definition}
Let $X$ and $Y$ be smooth schemes  over the fixed $\D$-ring $A$.
By a {\it $\D$-morphism of order $r$} we mean a rule $f:X \ra Y$
that attaches to any $\D$-prolongation sequence $S^*$ of
$p$-adically complete rings, a map of sets $X(S^0) \ra Y(S^r)$
that is ``functorial'' in $S^*$ in the obvious sense.
\end{definition}

For any
$\D$-prolongation sequence $S^*$, the shifted sequence $S^{*}[i]$,
$S[i]^n:=S^{n+i}$, is a new $\D$-prolongation sequence. Thus, any morphism
$f:X \ra Y$ of order $r$ induces maps of sets $X(S^i) \ra Y(S^{r+i})$
that are functorial in $S^*$. We can compose $\D$-morphisms
$f:X \ra Y$, $g:Y \ra Z$ of orders $r$ and $s$, respectively, and get a
$\D$-morphisms $g \circ f:X \rightarrow Z$ of order $r+s$.
There is a natural map from the set of $\D$-morphisms
$X \ra Y$ of order $r$ into the set of $\D$-morphisms $X \ra Y$
of order $r+1$, induced by the maps $Y(S^r) \ra Y(S^{r+1})$
arising from the $S^r$-algebra structure of $S^{r+1}$.

Recall that given a scheme $X$ over the ring $A$, if $B$ is an
$A$-algebra we let $X(B)$ denote the set of all morphisms of
$A$-schemes ${\rm Spec}\, B \rightarrow X$, and any such morphism
is called a $B$-{\it point} of $X$. For instance, recall that if
$X={\mathbb A}^1={\mathbb G}_a={\rm Spec}\, A[y]$, then $X(B)$ is
simply the set $B$ itself because a morphism ${\rm Spec}\,  B
\rightarrow {\rm Spec}\, A[y]$ is the same as a morphism $A[y]
\rightarrow
 B$, and the latter is uniquely determined by the image of $y$ in $B$. If
on the other hand, $X={\mathbb G}_m={\rm Spec}\, A[y,y^{-1}]={\rm Spec}\,
A[x,y]/(xy-1)$, then $X(B)=B^{\times}$ because ${\rm Hom}_{A}(A[y,y^{-1}],
B) =B^{\times}$ via the map $f \mapsto f(y)$. Finally, if
$X={\rm Spec}\, A[x,y]/(f(x,y))$, then $X(B)=\{(a,b) \in B^2: \, f(a,b)=0\}$.

\begin{definition}
Let $G$ and $H$ be smooth group schemes over $A$. We say that $G
\ra H$ is a $\D$-{\it homomorphism of order} $r$ if it is a
$\D$-morphism of order $r$ such that, for any prolongation
sequence $S^{*}$, the maps $X(S^0) \ra Y(S^r)$ are group
homomorphisms. A $\D$-{\it character} of order $r$ of $G$ is a
$\D$-homomorphism $G \rightarrow {\mathbb G}_a$ of order $r$,
where ${\mathbb G}_a={\rm Spec}\, A[y]$ is the additive  group
scheme over $A$. The group of $\D$-characters of order $r$ of $G$
will be denoted by $\bX^r_{\di \ddi}(G)$.
\end{definition}

Note that the group $\bX^r_{\di \ddi}(G)$ has a natural structure
of $A[\phip,\dd]$-module. We shall view the $\D$-characters of $G$
as ``linear partial differential operators'' on it. Of course,
they are highly ``non-linear'' in the affine coordinates around
various points of $G$.

Let $\psi:G \ra {\mathbb G}_a$ be a $\D$-character of a commutative
smooth group scheme $G$ over $A$. For any $\D$-ring $B$ over $A$,
$\psi$ induces a ${\mathbb Z}$-linear map
\begin{equation}
\label{9991} \psi:G(B) \ra {\mathbb G}_a(B)=B\, .
\end{equation}

\begin{definition}
The {\it group of solutions} of $\psi$ in $B$ is the kernel of the map
{\rm (\ref{9991})}, that is to say, the group
$$\{u \in G(B)\ | \ \psi u=0\}\, .$$
Given a subgroup $\Gamma$ of $G(B)$, the {\it group of solutions} of $\psi$
in $\Gamma$ is will be the group
$$\{u \in \Gamma\ | \ \psi u=0\}\, .$$
\end{definition}

\begin{example}
\label{cucurigu}
Let us assume that $A=R$, and that $G$ is a commutative smooth group scheme
over $R$. We may consider the natural ring homomorphisms
$$\begin{array}{rcl}
R[[q]] & \ra &  R \\
U(q)=\sum_{n \geq 0} a_n q^n &  \mapsto & U(0):=\left( \sum_{n \geq 0}
a_n q^n \right)_{|q=0}= a_0
\end{array}$$
and
$$
\begin{array}{rcl}
R[[q^{-1}]] & \ra &  R \\
U(q)=\sum_{n \leq 0} a_n q^n  & \mapsto &
U(\infty):=\left( \sum_{n \leq 0} a_n q^n \right)_{|q^{-1}=0}=
a_0 \,
\end{array}$$
respectively, and the corresponding induced group homomorphisms
\begin{equation}\label{mabuze}
\begin{array}{ccc}
G(R[[q]]) &  \ra & G(R) \\
u & \mapsto &
u(0)
\end{array}
\end{equation}
and
\begin{equation}\label{mabuze2}
\begin{array}{ccc}
G(R[[q^{-1}]]) &  \ra & G(R) \\
u & \mapsto & u(\infty)
\end{array}\, .
\end{equation}
If $e$ denotes the identity, we have the natural diagram of groups
\begin{equation}\label{dia}
\begin{array}{ccccc}
G(qR[[q]]) & \ra & G(R[[q]]) & \ra & G(R((q))\wh)\\
\ua  & \  & \ua & \  & \ua\\
\{e\}  & \ra  & G(R) & \ra & G(R[q,q^{-1}]\wh)\\
\da  & \  & \da & \  & \da\\
G(q^{-1}R[[q^{-1}]]) & \ra & G(R[[q^{-1}]]) & \ra &
G(R((q^{-1}))\wh)
\end{array}\, ,
\end{equation}
where $G(qR[[q]])$ and $G(q^{-1}R[[q^{-1}]])$ are the kernels of
the maps in (\ref{mabuze}) and (\ref{mabuze2}), respectively.

We will think of the spectrum of the ring
$A=R((q))^{\widehat{\phantom{p}}}$ as the $pq$ ``plane.'' The $p$-axis
is the ``arithmetic'' direction, while the $q$-axis is ``geometric;''
$q$ will be viewed as the ``exponential of $-2 \pi it$''
where $t$ is ``time.'' The operators $\d$ and
$\dd$ are ``vector fields'' along these two directions.

Under these interpretations,  the elements of
$R((q))^{\widehat{\phantom{p}}}$ play the role of ``functions in
the variables'' $p,q$. By considering (infinitely many) negative
powers of $q$, we allow the corresponding ``function''  to have an
``essential singularity as time approaches $-i \infty$.'' The
elements $u=\sum_{n \geq 0} a_nq^n$ of the power series ring
$R[[q]]$ inside $R((q))^{\widehat{\phantom{p}}}$ are viewed as the
analogues of functions in two variables that, as time goes to
infinity, approach a well defined limit function of the spatial
variable, the function $a_0=u(0)$. Here, the elements of the
monoid $R^{\d}=\{ x\in R : \; \d x=0\}$ are to be thought of as
the ``constant functions.''

Let $\psi$ be a $\D$-character of $G$. Considering the groups of
solutions of $\psi$ in the various groups of (\ref{dia}) above, we get a
diagram
\begin{equation}
\label{dia2}
\begin{array}{ccccc}
\cU_1 & \ra & \cU_+ & \ra & \cU_{\ra}\\
\ua  & \  & \ua & \  & \ua\\
\{e\}  & \ra  & \cU_0 & \ra & \cU_{\da}\\
\da  & \  & \da & \  & \da\\
\cU_{-1} & \ra  & \cU_- & \ra & \cU_{\la}
\end{array}\, .
 \end{equation}
 We will usually denote by $\cU_*$  the groups in the above
 diagram.
The elements of $\cU_0$ will be referred to as {\it stationary
solutions} of $\psi$, and are interpreted as solutions that ``do
not depend on time.''

For $q_0 \in pR$, we may consider the ring homomorphism
$$\begin{array}{ccc}
R[[q]] & \ra & R \\
q & \mapsto & q_0
\end{array}\, .$$
This is a $\D$-ring homomorphism for $q_0=0$, but not for $q_0 \neq 0$.
The embedding ${\rm Spec}\, R \ra {\rm Spec}\, R[[q]]$ will be
viewed as a ``curve'' in the ``$pq$-plane.''

Corresponding to the homomorphism above, we may consider the specialization
homomorphism
$$\begin{array}{ccc}
G(R[[q]]) & \ra & G(R) \\
u & \mapsto & u(q_0)
\end{array}\, .$$
The ``curve'' associated to $q_0 \neq 0$ will be thought of as a
curve along which we impose ``boundary conditions.'' For if $u_0
\in G(R)$, we will consider the ``boundary value problem'' with
that initial data, that is to say, the problem of finding a
solution (possibly unique) $u \in \cU_{+}$ such that $u(q_0)=u_0$.
When $q_0=0$, $u(0) \in G(R)$ is thought of as the ``limit of $u$
as time goes to infinity,'' as mentioned earlier. This situation
can be similarly stated for $R[[q^{-1}]]$ instead of $R[[q]]$.

We claim that the following hold:
\begin{equation}
\label{moderate}
\begin{array}{rclll}
\psi(u(0)) & = & (\psi u)(0) & \text{for} & u \in G(R[[q]])\, ,\\
\psi(u(\infty)) & = & (\psi u)(\infty) & \text{for} & u \in
G(R[[q^{-1}]])\, .
\end{array}
\end{equation}
Indeed, in order
 to prove the first of these identities, it is enough to observe
that, by the functorial definition of $\psi$, we have a commutative diagram
$$
\begin{array}{ccc}
G(R[[q]]) & \stackrel{\psi}{\ra} & R[[q]]\\
\da & \  & \  \da\\
G(R) & \stackrel{\psi}{\ra} & R
\end{array}\, ,
$$
where the vertical arrows are induced by the $\D$-ring
homomorphism
$$\begin{array}{ccc}
R[[q]] & \ra  & R \\ q & \mapsto & 0
\end{array}\, .$$
The argument for the second of the identities is similar.

The identities in (\ref{moderate}) imply that
\begin{equation}
\label{sarrah} \cU_{\pm}  = \cU_{0}  \times  \cU_{\pm 1},
\end{equation}
where $\times$ denotes the ``internal direct product.''

We may unify this picture by considering the various groups
of solutions $\cU_*$ as subgroups of a single group. Indeed, let us
consider the direct sum homomorphism
$$
\psi \oplus \psi\, : \; \frac{G(R((q))\h) \oplus
G(R((q^{-1}))\h)}{G(R[q,q^{-1}]\h)} \ra \frac{R((q))\h \oplus
R((q^{-1}))\h}{R[q,q^{-1}]\h}\, ,
$$
where the denominators in the domain and range are diagonally embedded into
the numerators. We then set
$$\cU := {\rm ker}\, (\psi \oplus \psi) \, ,$$
and refer to this group as the {\it group of generalized solutions} of $\psi$.
Its elements are the analogues of the generalized solutions of Sato.

The groups $\cU_{\la}$
and $\cU_{\ra}$ naturally embed into $\cU$ via the maps $x \mapsto
(-x,0)$ and $x \mapsto (0,x)$, respectively, and the restrictions
of these two embeddings to $\cU_{\da}$ coincide. Thus, all the
groups $\cU_*$ in the diagram (\ref{dia2}) can  be identified with
subgroups of $\cU$, and we have
$$
\begin{array}{rcl}
\cU_{\ra} \cap \cU_{\la} & = & \cU_{\da}\, , \\
\cU_+ \cap \cU_- & = & \cU_0\, , \\
\cU_{1} \cap \cU_{-1} & = & \{e\}\, .
\end{array}
$$
Note that, a priori, we do not have $\cU=\cU_{\ra} \cdot
\cU_{\la}$.

Next we introduce a  concept of {\it convolution} in our setting.
 Let ${\mathbb
Z}\mu(R)$ be the group ring of the group $\mu(R)$. That is to say,
${\mathbb Z} \mu(R)$ is the set of all functions $f:\mu(R) \ra
{\mathbb Z}$ of finite support, equipped with pointwise addition,
and multiplication $\star$ given by convolution
$$(f_1 \star f_2)(\zeta):=\sum_{\zeta_1 \zeta_2=\zeta}
f_1(\zeta_1) f_2(\zeta_2)\, .$$ For any integer $\kappa \in \bZ$,
we consider the ring homomorphism
\begin{equation}
\label{cuofi}
\begin{array}{ccc}
\bZ \bmu(R) & \stackrel{[\kappa]}{\ra} & \bZ \bmu(R)\\
f & \mapsto & f^{[\kappa]}
\end{array}\, , \end{equation}
where $f^{[\kappa]}(\zeta):= f(\zeta^{\kappa})$. If ${\mathbb
Z}[\mu(R)] \subset R$ is the subring generated by $\mu(R)$, then
we have a natural surjective (but not injective) ring homomorphism
$$\begin{array}{ccc}
{\mathbb Z} \mu(R) & \rightarrow & {\mathbb Z}[\mu(R)] \\
 f & \mapsto & f^{\sharp}:=\sum_{\zeta} f(\zeta) \zeta
\end{array}\, .$$

Now let $R((q^{\pm 1}))\wh$ be either
$R((q))\wh$ or $R((q^{-1}))\wh$, respectively.
For $\zeta \in \mu(R)$, we may consider the $R$-algebra isomorphism
$$\begin{array}{ccc}
\sigma_{\zeta}:R((q^{\pm 1}))\wh & \ra &  R((q^{\pm 1}))\wh \\
u(q) & \mapsto & u(\zeta q)
\end{array}$$
that induces a group isomorphism
$$\sigma_{\zeta}:G(R((q^{\pm 1}))\wh ) \ra G(R((q^{\pm 1}))\wh)\, .$$
Hence, the group $G(R((q^{\pm 1}))\wh)$ acquires a natural structure of
${\mathbb Z}
\mu(R)$-module via convolution. Indeed, for any $f \in {\mathbb Z} \mu(R)$ and
$u \in G(R((q^{\pm 1}))\wh)$, we have
$$f \star u:=\sum_{\zeta} f(\zeta) \sigma_{\zeta}(u)\, ,$$
essentially a ``superposition" of translates of $u$.

Note now that the homomorphism $\sigma_{\zeta}:R((q^{\pm 1}))\wh \ra R
((q^{\pm 1}))\wh $
is actually a $\D$-ring homomorphism: the commutation of
$\sigma_{\zeta}$ and $\dd$ is clear, while the commutation of
$\sigma_{\zeta}$ and $\d$ follows because
$$\phip(\sigma_{\zeta}(q^{\pm1}))=\phip(\zeta^{\pm1} q^{\pm1})=
\zeta^{\pm \phi}
q^{\pm p}=\zeta^{\pm p}q^{\pm p}=\sigma_{\zeta}(q^{\pm p})=\sigma_{\zeta}(
\phip(q^{\pm}))\, .$$

If $\psi$ is a $\D$-character of $G$, by the functorial definition of
$\D$-characters we obtain a commutative diagram
$$\begin{array}{ccc}
G(R((q^{\pm 1}))\wh ) & \stackrel{\psi}{\ra} & R((q^{\pm 1}))\wh\\
\sigma_{\zeta} \da & \  & \da \sigma_{\zeta}\\
G(R((q^{\pm 1}))\wh ) & \stackrel{\psi}{\ra} & R((q^{\pm 1}))\wh
\end{array}\, .
$$
Hence, for any $f \in {\mathbb Z} \mu(R)$ and $u \in G(R((q^{\pm
1}))\wh)$, we have
$$\psi(f \star u)=f \star \psi(u)\, .$$
That is to say, $\psi$ is $\bZ \bmu(R)$- module homomorphism. In
particular, the groups $\cU_{\pm 1}$ of solutions of $\psi$ in
$G(R((q^{\pm 1}))\wh)$ are ${\mathbb Z} \mu(R)$- submodules of
$G(R(( q^{\pm 1}))\wh)$ respectively. Morally speaking, this says
that a superposition of translates of a solution is again a
solution.

Let us additionally assume that $G$ has relative dimension one
over $R$. Let $T$ be an \'{e}tale coordinate on a neighborhood of
the zero section in $G$ such the zero section is given scheme
theoretically by $T=0$. Then the ring of functions on the
completion of $G$ along the zero section is isomorphic to a power
series ring $R[[T]]$. We fix such an isomorphism. Then we have an
induced group isomorphism
\begin{equation}
\label{hacha}
\iota: q^{\pm 1}R[[q^{\pm 1}]] \ra G(q^{\pm 1}R[[q^{\pm 1}]])\,  ,
\end{equation}
where $q^{\pm 1}R[[q^{\pm 1}]]$ is a group relative to a formal
group law $\cF(T_1,T_2) \in R[[T_1,T_2]]$ attached to $G$. We
recall that the series $[p^n](T) \in R[[T]]$, defined by
multiplication by $p^n$ in the formal group $\cF$, belongs to the
ideal $(p,T)^n$. We equip $G(q^{\pm 1}R[[q^{\pm 1}]])$ with the
topology induced via $\iota$ from the $(p,q^{\pm 1})$-adic
topology on $R[[q^{\pm 1}]]$. Since $G(q^{\pm 1}R[[q^{\pm 1}]])$
is complete in this topology, the ${\mathbb Z} \mu(R)$-module
structure of $G(q^{\pm 1}R[[q^{\pm 1}]])$ extends uniquely to a
${\mathbb Z} \mu(R)\wh$-module structure on $G(q^{\pm 1}R[[q^{\pm
1}]])$ in which multiplication by scalars is continuous. Here,
${\mathbb Z} \mu(R)\wh$ is, of course, the $p$-adic completion of
${\mathbb Z} \mu(R)$. We observe that the ring homomorphism
\begin{equation}
\label{notinj} \begin{array}{rcl}
{\mathbb Z} \mu(R)\wh & \stackrel{\sharp}{\ra} &  R \\
f & \mapsto & f^{\sh}:=\sum_{\zeta} f(\zeta)\zeta
\end{array}
\end{equation}
is surjective (but not injective). Clearly, the groups $\cU_{\pm
1}$ of solutions of $\psi$ are ${\mathbb Z} \mu(R)\wh$- submodules
of $\cU_{\pm}$, respectively.

Also, we note that a Mittag-Leffler argument shows that for any
$\kappa$ coprime to $p$ the map $[\kappa]:\bZ \bmu(R)\h \ra \bZ
\bmu(R) \h$ is surjective.

We end our discussion here by introducing boundary value operators
at $q^{\pm 1}=0$. We start by considering, for any $0 \neq \kappa
\in \bZ$, operators
$$\begin{array}{rcl}
\Gamma_{\kappa}:R[[q^{\pm 1}]] & \ra & R \vspace{1mm} \\
u=\sum a_n q^n & \mapsto  & \Gamma_{\kappa} u=a_{\kappa}
\end{array}\, .
$$
We fix also a $\dd$-character $\psi_q$ of $G$. In the
applications, for $G=\bG_a$, $\psi_q$ will be the identity; and
for $G$ either $\bG_m$ or an elliptic curve $E$ over $R$, $\psi_q$
will be the ``Kolchin logarithmic derivative.'' If we fix a
collection of non-zero integers ${\mathcal K}_{\pm} \subset
\bZ_{\pm}$, and set $\rho_{\pm}:=\sharp {\mathcal K}_{\pm}$, then
we may consider the {\it boundary value operator at $q^{\pm
1}=0$},
\begin{equation}
\label{vocea}
\begin{array}{rcl}
G(R[[q^{\pm 1}]])  & \stackrel{B_{\pm}^0}{\ra} & R^{\rho_{\pm}}\\
B_{\pm}^0 u & = & (B_{\kappa}^0 u)_{\kappa \in \cK_{\pm}}
\end{array} \, ,
\end{equation}
where $B_{\kappa}^0 u:= \Gamma_{\kappa} \psi_q u$. Here we view
$R^{\rho_{\pm}}$ as a direct product
$R^{\rho_{\pm}}=\prod_{\kappa \in \cK_{\pm}} R_{\kappa}$
of copies $R_{\kappa}$ of $R$ indexed by $\cK_{\pm}$.

It is easy to check that $B_{\pm}^0$ is a $\bZ \bmu(R)$-module homomorphism
provided that $R^{\rho_{\pm}}$ be viewed as a $\bZ \bmu(R)$-module with the
following structure: for each $\kappa \in {\mathcal K}_{\pm}$, the
copy of $R$ indexed by $\kappa$ in $R^{\rho_{\pm}}$ is a $\bZ \bmu(R)$-module
via the ring homomorphism
\begin{equation}
\label{notinjj}
\begin{array}{c}
\bZ \bmu(R) \stackrel{[\kappa]}{\ra} \bZ \bmu(R)
\stackrel{\sh}{\ra} R\, .
\end{array}
\end{equation}
In other words, for $f \in \bZ \bmu(R)$ and $(r_{\kappa})_{\kappa}
\in R^{\rho_{\pm}}$, we have that $$f \cdot
(r_{\kappa})_{\kappa}:=((f^{[\kappa]})^{\sharp}
r_{\kappa})_{\kappa}.$$
 The restriction of $B_{\pm}^0$ to
$G(q^{\pm 1}R[[q^{\pm 1}]])$,
\begin{equation}
\begin{array}{rcl}
B_{\pm}^0:G(q^{\pm 1}R[[q^{\pm 1}]])  & \ra & R^{\rho_{\pm}},
\end{array}\end{equation}
is a $\bZ \bmu(R)\h$-module homomorphism where $R^{\rho_{\pm}}$
is viewed as a $\bZ \bmu(R)\h$-module with the following
structure: for each $\kappa \in {\mathcal K}_{\pm}$, the copy of
$R$  indexed by $\kappa$ in $R^{\rho_{\pm}}$ is a $\bZ \bmu(R)\h$-module via
the ring homomorphism
\begin{equation}
\label{notinjjj}
\begin{array}{c}
\bZ \bmu(R)\h \stackrel{[\kappa]}{\ra} \bZ \bmu(R)\h
\stackrel{\sh}{\ra} R\, .
\end{array}
\end{equation}
It our applications, it will sometimes be the case that the $\bZ
\bmu(R)\h$-module structure of $\cU_{\pm 1}$ induces a certain
structure of $R$-module, and the restriction of $B_{\pm}^0$ to
$\cU_{\pm 1}$ becomes an $R$-module homomorphism.
\end{example}

\begin{remark}
The discussion in the previous Example  applies to the case where
$G$ is one of the following group schemes over $A=R$:
\begin{enumerate}
\item $G={\mathbb G}_a:={\rm Spec}\, R[y]$, the additive group.
\item $G={\mathbb G}_m:={\rm Spec}\, R[y,y^{-1}]$, the multiplicative group.
\item $G=E$, an elliptic curve over $R$.
\end{enumerate}
The discussion, however, does not apply to the following
case, which will also interest us later:
\begin{enumerate}
\item[(4)] $G=E$, an elliptic curve over $A=R((q))\h$ that does not
descend to $R$, for instance, the Tate curve.
\end{enumerate}
In this latter case, we will need new definitions for some of the  groups
(\ref{dia2}) (cf. the discussion in our subsection on the Tate curve).
\end{remark}

\section{Partial differential jet spaces of schemes}
\setcounter{theorem}{0} In this section we introduce
arithmetic-geometric partial differential jet spaces of schemes. We do this
by analogy with the ordinary differential ones in geometry
\cite{hermann} and arithmetic \cite{char,difmod,book},
respectively. We also record some of their general properties
(that can be proven in a manner similar to the corresponding proofs
in the ordinary case \cite{book}, proofs that, therefore, will be omitted
here.) We recall that the $\D$-ring $A$ we have fixed is a
$p$-adically complete  Noetherian integral domain of characteristic zero.

For any scheme $X$ of
finite type over $A$, we define a projective system of $p$-adic
formal schemes
\begin{equation}
\label{fite}
\cdots \ra J^r_{\di \ddi}(X) \ra J^{r-1}_{\di \ddi}(X) \ra
\cdots \ra J^1_{\di \ddi}(X) \ra J^0_{\di \ddi}(X)=\hat{X}\, ,
\end{equation}
called the $\D$-{\it jet spaces} of $X$.

Let us assume first that $X$ is
affine, that is to say, $X={\rm Spec}\, A[y]/I$, where $y$ is a tuple of
indeterminates. We then set
$$J^n_{\di \ddi}(X):=Spf\ A[y,Dy,\ldots ,D^n y]\h/(I,DI,\ldots ,D^n I)\, .$$
The sequence of rings $\{\cO(J^n_{\di \ddi}(X))\}_{n \geq 0}$
has a natural structure of $\D$-prolongat- ion sequence
$\cO(J^*_{\di \ddi}(X))$, and the latter has the following
universality property that can be easily checked.

\begin{proposition}
\label{universality} For any $\D$-prolongation sequence $S^*$
that consists of $p$-adically complete rings, and for any homomorphism
$u:\cO(X) \ra S^0$ over $A$, there is a unique morphism of $\D$-prolongation
sequences
$$u_*=(u_n):\cO(J^*_{\di \ddi}(X)) \ra S^*\, ,$$
with $u_0=u$.
\end{proposition}

{\it Proof}. Similar to \cite{book}, Proposition 3.3. \qed

As a consequence we get that, for affine $X$, the construction
$X \mapsto J^r(X)$ is compatible with localization in the following
sense:

\begin{corollary}\label{ma duc acum}
If $X={\rm Spec}\,  B$ and $U={\rm Spec}\, B_f$, $f \in B$,
then
$$\cO(J^r_{\di \ddi}(U))=(\cO(J^r_{\di \ddi}(X))_f)\wh\, .$$
Equivalently,
$$J^r_{\di \ddi}(U)=J^r_{\di \ddi}(X) \times_{\hat{X}} \hat{U}\, .$$
\end{corollary}

{\it Proof}. Similar to \cite{book}, Corollary 3.4. \qed
\medskip

Consequently, for $X$ that is not necessarily affine, we can define
a formal scheme $J^r_{\di \ddi}(X)$ by gluing the various schemes $J^r_{\di
\ddi}(U_i)$ for $\{U_i\}$ an affine Zariski open covering of $X$.
The resulting formal scheme will have a corresponding universality
property, whose complete formulation and verification we leave to the reader.

\begin{remark}
If $\bG_a ={\rm Spec}\,  A[y]$ is the additive group scheme over $A$, then
$$J^n_{\di \ddi}(\bG_a)=Spf\ A[y,Dy,\ldots,D^n y]\wh \, .$$
If $\bG_m ={\rm Spec}\,  A[y,y^{-1}]$ is the multiplicative group scheme
over $A$, then
$$J^n_{\di \ddi}(\bG_m)=Spf\ A[y,y^{-1},Dy,\ldots,D^n y]\wh \, .$$
\end{remark}

\begin{remark}
By the universality property of jet spaces,
we have that
$$J^n_{\di \ddi}(X \times Y) \simeq J^n_{\di \ddi}(X)
\times J^n_{\di \ddi}(Y)$$
where the product on the left hand side is taken in the category
of schemes of finite type over $A$, and the product on the
right hand side is taken in the category of formal schemes over $A$.
In the same vein, if $X$ is a group scheme of finite type
over $A$, then (\ref{fite}) is a projective system of groups in the category of
$p$-adic formal schemes over $A$.
\end{remark}

\begin{remark}
By the universality property of jet spaces, the set of order $r$ $\D$-morphisms
$X \ra Y$ between two schemes of finite type over $A$
naturally identifies with the set of morphisms over $A$
of formal schemes $J^r_{\di \ddi}(X) \ra J^0_{\di \ddi}(Y)=\hat{Y}$.
 In particular, the set $\cO^r_{\di \ddi}(X)$
of all order $r$ $\D$-morphisms $X \ra \hat{\bA}^1$ identifies with the ring
of global functions $\cO(J^r_{\di \ddi}(X))$. If $G$ is a group scheme of
finite type over $A$, then
the group $\bX_{\di \ddi}^r(G)$
 of order $r$ $\D$-characters $G \ra \bG_a$ identifies
with the group of homomorphisms $J^r_{\di \ddi}(G) \ra \hat{\bG}_a$, and
thus, it identifies with an $A$-submodule of $\cO(J^r_{\di
\ddi}(G))$. Let
$$\cO^{\infty}_{\di \ddi}(X):=\lim_{\ra} \cO^r_{\di \ddi}(X)$$
be the $\D$-ring of all $\D$-morphisms $X \ra \bA^1$, and let
$$\bX^{\infty}_{\di \ddi}(X):=\lim_{\ra} \bX^r_{\di \ddi}(X)$$
be the group of $\D$-characters $G \ra \bG_a$. Then
$\cO^{\infty}_{\di \ddi}(X)$ has a natural structure of
 $A[\phip,\dd]$-module, and
$\bX^{\infty}_{\di \ddi}(X)$ is an $A[\phip,\dd]$-submodule.
\end{remark}

\begin{proposition}
\label{local}
Let $X$ be a smooth affine scheme over $A$, and let $u:A[y] \ra
\cO(X)$ be an \'{e}tale morphism, where $y$ is a $d$-tuple of
indeterminates. Let $y^{(i,j)}$ be $d$-tuples of indeterminates,
where $i,j \geq 0$. Then the natural morphism
$$\cO(\hat{X})[y^{(i,j)}\mid _{1 \leq i+j \leq n} ]\wh \ra \cO(J^n_{\di
\ddi}(X))
$$
that sends $y^{(i,j)}$ into $\d^i \dd^j (u(y))$ is an
isomorphism. In particular, we have an isomorphism of formal schemes over
$A$
$$J^n_{\di \ddi}(X) \simeq \hat{X} \times \hat{\bA}^{\frac{n(n+3)}{2}d}\, .$$
\end{proposition}

{\it Proof}. The argument is similar to that used in Proposition 3.13 of
\cite{book}. \qed

\begin{corollary}
If $Y \ra X$ is an \'{e}tale morphism of smooth schemes over $A$
then
$$J^n_{\di \ddi}(Y) \simeq J^n_{\di \ddi}(X) \times_{\hat{X}}
\hat{Y}\, .$$
\end{corollary}

\begin{remark}
The jet spaces with respect to a derivation \cite{annals} that arises in the
Ritt-Kolchin theory \cite{kolchin} \cite{ritt} are ``covered"
by of our $\D$-jet spaces here. The same holds for the $p$-jet spaces
(with respect to a $p$-derivation) constructed in \cite{char}.
More precisely, if we for affine $X/A$ we set
$$\begin{array}{rcl}
J^n_{\ddi}(X) & := & {\rm Spec}\, A[y,\dd y,\ldots,\dd^n y]/(I,\dd I,
\ldots,\dd^n I)\, , \\
J^n_{\di}(X) & := & Spf\ A[y,\d y,\ldots,\d^n y]\h/(I,\d I,\ldots,\d^n I)\, ,
\end{array}$$
then we have natural morphisms
$$\begin{array}{rcl}
J^n_{\di \ddi}(X) & \ra & J^n_{\ddi}(X)\wh \, , \\
J^n_{\di \ddi}(X) & \ra & J^n_{\di}(X)\, .
\end{array}$$
The same holds then for any scheme $X$, not necessarily affine. The elements
of $\cO^r_{\ddi}(X):= \cO(J^r_{\ddi}(X))$ identify with the
$\dd$-morphisms $X \ra \bA^1$; the elements of
 $\cO^r_{\di}(X):=\cO(J^r_{\di}(X))$ identify with the
$\d$-morphisms $X \ra \bA^1$. The homomorphisms $J^r_{\ddi}(G) \ra
\hat{\bG}_a$ identify with the $\dd$-characters $G \ra \bG_a$; the
homomorphisms $J^r_{\di}(G) \ra \hat{\bG}_a$ identify with the
$\d$-characters $G \ra \bG_a$.
\end{remark}

\begin{remark}
\label{ginge} The functors $X \mapsto J^r_{\di \ddi}(X)$ from the
category ${\mathcal C}$ of $A$-schemes of finite type to the
category $\hat{\mathcal C}$ of $p$-adic formal schemes naturally
extends to a functor from ${\mathcal B}$ to $\hat{\mathcal C}$,
where ${\mathcal B}$ is the category whose objects as the same as
those in ${\mathcal C}$, hence $Ob\, {\mathcal C}=Ob\, {\mathcal
B}$, and whose morphisms are defined by
$${\rm Hom}_{\mathcal B}(X,Y):={\rm Hom}_{\hat{\mathcal C}}(\hat{X},\hat{Y})
$$
for all $X,Y \in Ob \, {\mathcal B}$.
\end{remark}

\section{Partial differential jet spaces of formal groups}
\setcounter{theorem}{0} In this section we attach to any formal
group law $\cF$ in one variable, certain formal groups in several
variables that should be thought of as arithmetic-geometric
partial differential jets of $\cF$. We use the interplay between
these and the partial differential jet spaces of schemes
introduced in the previous section to prove that the module of
$\D$-characters is finitely generated. We also define the notion
of Picard-Fuchs symbol of a $\D$-character, which will play a key
r\^ole later.

Let $A$ be our fixed $\D$-ring and $y$ a $d$-tuple of variables.
For any $G \in A[[y,Dy,\ldots,D^n y]]$, we let $G_{|y=0} \in
A[[Dy,\ldots,D^n y]]$
denote the series obtained from $G$ by setting $y=0$, while keeping
$\d^i \dd^j y$ unchanged for $i+j \geq 1$. We recall that $L$ stands
for the field of fractions of $A$.

\begin{lemma}
\label{carmen} For $a,b \in \bZ_+$, we set $B:=A[Dy,\ldots,D^{a+b} y]\wh$.
Then:
\begin{enumerate}
\item[1)] If $G \in A[[y]][\dd y,\ldots,\dd^b y]$, then $(\phi^a G)_{|y=0}
\in B$.
\item[2)] If $F \in A[[y]]$, then $(\d^a \dd^b F)_{|y=0} \in B$.
\item[3)] In the case where $y$ is a single variable, if $F =\sum_{n \geq 1}
c_n y^n \in yL[[y]]$ and $nc_n \in A$, then $(\phip^a \dd^b F)_{|y=0} \in B$
and $(\phip^a F)_{|y=0} \in pB$.
\end{enumerate}
\end{lemma}

{\it Proof}. For the proof assertion 1), we may assume that
$G=a(y) \cdot (\dd y)^{i_1}\ldots (\dd^b y)^{i_b}$, $a(y) \in A[[y]]$.
We have that
\begin{equation}
\label{dif}
y^{\phi^a}-y^{p^a}=p \Phi_a\, , \quad
\Phi_a \in \bZ[y,\d y,\ldots,\d^a y]\, ,
\end{equation}
and so we get
$(\phip^a G)_{|y=0}  =  (\phip^a a(y))_{|y=0} \cdot [(\phip^a \dd
y)^{i_1}\ldots (\phip^a \dd^b y)^{i_b}]_{|y=0}$.

It is clear that $[(\phip^a \dd y)^{i_1}\ldots(\phip^a \dd^b y)^{i_b}]_{|y=0}
 \in \bZ[Dy,\ldots,D^{a+b} y]$. On the other hand, if $a(y)=\sum_{n \geq
0} a_n y^n$, by (\ref{dif}) we then have that
$$\begin{array}{rcl}
(\phip^a a(y))_{|y=0} & = & (\sum_{n \geq 0} a_n^{\phi^a}
 (y^{p^a}+p\Phi_a)^n)_{|y=0}\\
 & = & \sum_{n \geq 0} a_n^{\phi^a} p^n((\Phi_a)_{|y=0})^n \, ,
\end{array}
$$
which is an element of  $A[\d y,\ldots,\d^a y]\wh$.

For the proof of the second part, we first use induction to check that
$\dd^b F \in A[[y]][\dd y,\ldots,\dd^b y]$. By the part
of the Lemma already proven, we conclude that
$(\phip^a \dd^b F)_{|y=0} \in A[Dy,\ldots,D^{a+b} y]\wh $.
Assertion 2) then follows because $p^a \d^a \dd^b F$ belongs to
the ring generated by $\phip^i \dd^b F$ with $i \leq a$.

For assertion 3), we first use induction to check that
\begin{equation}
\label{ssst}
\dd^b F=\sum_{n \geq 1} (\dd^b c_n) y^n +G_b,\quad G_b \in
A[[y]][\dd y,\ldots,\dd^b y]\, .
\end{equation}
Indeed, this is true for $b=0$, with $G_0=0$, and if true for some $b$, then
$$\dd^{b+1}F=\sum_{n \geq 1} (\dd^{b+1} c_n) y^n +
[\sum_{n \geq 1} (\dd^b c_n) n y^{n-1}]\dd y +\dd G_b\, ,$$
by the induction hypothesis and the fact that $nc_n \in A$.
By (\ref{ssst}), we derive that
$$(\phip^a \dd^b F)_{|y=0}=
[\sum_{n \geq 1} (\phip^a \dd^b c_n)(y^{\phi^a})^n]_{|y=0}+
(\phip^a G_b)_{|y=0}\, .
$$
By the second part of the Lemma proven above,
$(\phip^a G_b)_{|y=0} \in A[\d^i \dd^j
y\mid_{1 \leq i+j \leq a+b} ]\wh$. On the other hand
$$\begin{array}{rcl}
[\sum_{n \geq 1} (\phip^a \dd^b c_n)(y^{\phip^a})^n]_{|y=0} & = &
[\sum_{n\geq 1} (\phip^a \dd^b c_n)(y^{p^a}+p \Phi_a)^n]_{|y=0}\\
\  & = & \sum_{n\geq 1} (\phip^a \dd^b c_n)p^n
((\Phi_a)_{|y=0})^n\, ,
\end{array}
$$
which belongs to $pA[Dy,\ldots,D^{a+b} y]\wh $ because
$c_n p^n=n c_n (p^n/n) \in pA$, and it defines a sequence that
converge to $0$ as $n \ra \infty$. \qed

Let now $T$ be one variable, $(T_1,T_2)$ a pair of ``copies'' of
$T$, and $\cF:=\cF(T_1,T_2) \in A[[T_1,T_2]]$ be a formal group
law in the  variable $T$. Then, as in \cite{book}, p. 124,  the
tuple
\begin{equation}
\label{popo} (\cF,D \cF,\ldots,D^n \cF) \in
A[[T_1,T_2,DT_1,DT_2,\ldots,D^n T_1,D^n T_2]]^{\frac{(n+1)(n+2)}{2}}
\end{equation}
is a formal group law over $A$ in the variables $T,DT,\ldots,D^nT$.
Consequently, the tuple
\begin{equation}
\label{popopo} (D \cF_{|T_1=T_2=0},\ldots,D^n \cF_{|T_1=T_2=0})
 \in A[[DT_1,DT_2,\ldots,D^nT_1,D^nT_2]]^{\frac{n(n+3)}{2}}
\end{equation}
is a formal group law over $A$ in the variables $DT,\ldots,D^nT$.
By the second part of Lemma \ref{carmen}, this series
(\ref{popopo}) belong to $A[DT_1,DT_2,\ldots,D^nT_1, D^nT_2]\wh$,
so they define a structure of group objects in the category of
formal schemes over $A$,
\begin{equation}
\label{papa}
(\hat{\bA}^{\frac{n(n+3)}{2}},[+]).
\end{equation}
Let $l=l(T) \in L[[T]]$ be the logarithm of $\cF$
(cf. \cite{sil}), and let $a+b \leq n$.
By assertion 3) of Lemma \ref{carmen}, we have that
$$L^{[a,b]}:=p^{\epsilon(b)}(\phip^a \dd^b l)_{|T=0} \in
A[DT,\ldots,D^nT]\wh \, ,$$
where $\epsilon(b)$ is either $0$ or $-1$ if $b>0$ or $b=0$, respectively.
So $L^{[a,b]}$ defines a morphism of formal schemes
$$L^{[a,b]}:\hat{\bA}^{\frac{n(n+3)}{2}} \ra \hat{\bA}^1\, .$$
As in \cite{book}, p. 125, the morphism $L^{[a,b]}$ above
is actually a homomorphism
\begin{equation}
\label{romanu}
L^{[a,b]}:(\hat{\bA}^{\frac{n(n+3)}{2}},[+]) \ra
\hat{\bG}_a=(\hat{\bA}^1,+)
\end{equation}
of group objects in the category of formal schemes.

We let $G$ be $\bG_a$, $\bG_m$, or an elliptic curve (defined by a
Weierstrass equation) over $A$, and let $\cF$ be the formal group
law naturally attached to $G$.
In the case $G=\bG_a={\rm Spec}\, A[y]$, we let $T=y$. In the case
$G=\bG_m={\rm Spec}\, A[y,y^{-1}]$,  we let $T=y-1$. And if $G=E$, we let
$T$ be an \'{e}tale coordinate in a neighborhood $U$ of the zero section $0$
such that $0$ is the zero scheme of $T$ in $U$.  The group
$${\rm ker} (J^n_{\di \ddi}(G) \ra J^0_{\di \ddi}(G)=\hat{G})$$
is isomorphic to the group (\ref{papa}).

Let $e(T) \in L[[T]]$ be the exponential of the formal group law
$\cF$ (that is to say, the compositional inverse of $l(T) \in L[[T]]$).
We have $e(pT) \in pTA[T]\h$. We may consider the map
\begin{equation}\label{prost}
\begin{array}{rcl}
A[[T]][DT,\ldots,D^r T]\wh & \stackrel{\circ
e(pT)}{\longrightarrow} & A[T,DT,\ldots,D^rT]\wh  \\
G \mapsto G \circ e(pT) & := & G(\d^i\dd^j (e(pT))\mid_{0 \leq i+j \leq r})
\end{array}\, .
\end{equation}
On the other hand, by Proposition \ref{local}
we may consider the natural map
$$\cO(J^r_{\di \ddi}(G)) \ra A[[T]][DT,\ldots,D^r T]\wh \, ,$$
which by the said Proposition is injective with torsion free
cokernel. We shall view this map as an inclusion.(Note that
this is the case, more generally, when $G$ is replaced by a smooth
scheme over $A$ and $T$ is an \'{e}tale coordinate.)

\begin{lemma}
\label{xxzz}
For any $\D$-character $\psi \in \cO(J^r_{\di
 \ddi}(G))$, the series $\psi \circ e(pT)$
is in the $A$-linear span of $\{\phip^i \dd^j T\mid_{0 \leq i+j \leq
r}\}$.
\end{lemma}

{\it Proof}. Clearly $F:=\psi \circ e(pT)$ is {\it additive}, that is to say,
$$F(T_1+T_2, \ldots,D^r(T_1+T_2))  = F(T_1,\ldots,D^rT_1)
 + F(T_2,\ldots,D^rT_2)\, .$$
But the only additive elements in
$$A[1/p][[T,DT,\ldots,D^rT]]=
A[1/p][[\phip^i \dd^jT\mid_{0 \leq i+j \leq r}]]$$
are those in the
$A[1/p]$-linear span of
 $\{\phip^i \dd^j T\mid_{0 \leq i+j \leq r}\}$. We are thus left to show that
if
$$\sum_{ij} a_{ij} \phip^i \dd^j T \in pA[T,DT,\ldots,D^rT]\wh $$
for $a_{ij} \in A$, then $a_{ij} \in pA$ for all $i,j$.
Let $\bar{a}_{ij} \in A/pA$ be the images of $a_{ij}$. Then
we have
$$\sum \bar{a}_{ij} (\dd^j T)^{p^i}=0 \in (A/pA)[T,\dd T,\ldots,\dd^r T]\, ,$$
which clearly implies $\bar{a}_{ij}=0$ for all $i,j$. \qed

\begin{corollary}
The $A$-module $\bX_{\di \ddi}^r(G)$ of $\D$-characters or order $r$
is finitely generated of rank at most
$$\frac{(r+1)(r+2)}{2}\, .$$
\end{corollary}

{\it Proof}. By Lemma \ref{xxzz}, we have that $\bX_{\di \ddi}^r(G)$ embeds
into a finitely generated $A$-module of rank at most
$\frac{(r+1)(r+2)}{2}$. The result follows because $A$ is Noetherian.
\qed
\medskip

In light of Lemma \ref{xxzz}, if $\psi$ is a $\D$-character of $G$ then,
 as an element of the ring $A[[T]]DT,\ldots,D^rT]\wh$,
$\psi$ can be identified with the series
\begin{equation}
\label{mamscu} \psi=\frac{1}{p} \sigma(\phip,\dd) l(T)\, ,
\end{equation}
where $\sigma=\sigma(\xi_p,\xi_q) \in A[\xi_p,\xi_q]$ is a
polynomial.

\begin{definition}
\label{piiccfu} The polynomial $\sigma$ is the {\it Picard-Fuchs
symbol} of $\psi$ with respect to $T$.
\end{definition}

Cf. \cite{char} for comments on the terminology. In a more general
context, we will later define what we call the {\it Fr\'{e}chet
symbol}, and will explain the relation between these two symbol
notions.

The following Lemmas will also be needed later.

\begin{lemma}
\label{caiin}
$$\left( L^{[a,b]} \circ e(pT) \right)_{|T=0}=
\left( p^{1+\epsilon(b)} \phip^a \dd^b T \right)_{|T=0}\, .$$
\end{lemma}

{\it Proof}. We have
$$\begin{array}{rcl}
\left( L^{[a,b]} \circ e(pT) \right)_{|T=0} & = &
p^{\epsilon(b)}(\phip^a \dd^b
l)(0,\d(e(pT)),\dd(e(pT)),\ldots)_{|T=0} \vspace{1mm} \\
 & = & p^{\epsilon(b)}(\phip^a \dd^b
l)(e(pT),\d(e(pT)),\dd(e(pT)),\ldots)_{|T=0} \vspace{1mm} \\
 & = & (p^{\epsilon(b)} \phip^a \dd^b (l(e(pT))))_{|T=0} \vspace{1mm} \\
 & = & \left( p^{1+\epsilon(b)} \phip^a \dd^b T \right)_{|T=0}\, .
\end{array}
$$
\qed

\begin{lemma}
\label{fu}
The family $\{L^{[a,b]}\mid_{1 \leq a+b\leq r}\}$ is $A$-linearly independent.
\end{lemma}

{\it Proof}. By Lemma \ref{caiin}, it is enough to check that the family
$$\{(\phip^a \dd^b T)\mid_{T=0, 1 \leq a+b \leq r}\}$$
is $A$-linearly independent. Let us assume that
$$\sum_{a+b \geq 1} \lambda_{ab} (\phip^a \dd^b T)_{|T=0}=0\, .$$
This implies that
$$\sum_{a+b \geq 1} \lambda_{ab} \phip^a \dd^b T \in TA[T,DT,D^2T,\ldots
] \subset TL[\phip^i \dd^j T\mid_{i\geq 0,j\geq 0}] \, ,$$
which clearly implies $\lambda_{ab}=0$ for all $a,b$.
\qed

\section{Fr\'{e}chet derivatives and symbols}
\setcounter{theorem}{0}

We now develop arithmetic-geometric analogues of some
classical constructions in the calculus of variations, including
Fr\'{e}chet derivatives and Euler-Lagrange equations. We use
the Fr\'{e}chet derivatives to define Fr\'{e}chet symbols of $\D$-characters,
and we relate Fr\'{e}chet symbols to the previously defined Picard-Fuchs
symbols.

\subsection{Fr\'{e}chet derivative}
We recall that for a smooth scheme $X$ over $A$, we denote by
$$T(X): ={\rm Spec}\, {\rm Symm}(\Omega_{X/A})$$
the tangent scheme of $X$. Also, we set
$$\cO_{\di
\ddi}^{\infty}(X)=\lim_{\ra} \cO_{\di \ddi}^r(X)\, .$$
Let $\pi:T(X) \ra X$ be the canonical projection. Using ideas in
\cite{book}, we construct a natural compatible sequence of
$A$-derivations 
\begin{equation} \label{classsoon} \Theta:\cO_{\di
\ddi}^r(X) \ra \cO_{\di \ddi}^r(T(X)) 
\end{equation}
inducing an $A$-derivation
\begin{equation}
\label{air1} \Theta:\cO_{\di \ddi}^{\infty}(X) \ra
\cO^{\infty}_{\di \ddi}(T(X)),
\end{equation}
respecting the filtration by orders. We call $\Theta f$ the
$\D$-{\it tangent map} or {\it Fr\'{e}chet derivative} of $f$. (In
the ordinary case treated in \cite{book} $\Theta$ was denoted by
$T$.) The construction is local
and natural, so it provides a global concept. Thus, we may assume that $X$
is affine.

For any ring $S$, we denote by $S[\epsilon]$ (where
$\epsilon^2=0$) the ring $S \oplus \epsilon S$
 of {\it dual numbers} over $S$.
Note that any prolongation sequence $S^*=\{S^r\}$ can
be uniquely extended to a
 prolongation sequence $S^*[\epsilon]=\{S^r[\epsilon]\}$ where
$\d \epsilon=\epsilon$, and so $\phip(\epsilon)=p \epsilon$, and
$\dd \epsilon=0$, respectively. In
particular, we have a $\D$-prolongation sequence
$\cO_{\di \ddi}^*(T(X))[\epsilon]=
\{\cO_{\di \ddi}^r(T(X))[\epsilon]\}$.

On the other hand, we have a natural inclusion $\cO(X) \subset \cO(T(X))$,
 and a natural derivation
\begin{equation}
\label{irinushe} d:\cO(X) \ra
\cO(T(X))={\rm Symm}(\Omega_{\cO(X)/A})
\end{equation}
induced by the universal K\"{a}hler derivation
$d:\cO(X) \ra \Omega_{\cO(X)/A}$. Hence, we have an $A$-algebra map
$$\begin{array}{rcl}
\cO(X) & \ra & \cO(T(X))[\epsilon] \\ f & \mapsto & f+\epsilon \cdot df
\end{array}\, .
$$
By the universality property of $\cO_{\di \ddi}^*(X)$, there are
naturally induced ring homomorphisms
\begin{equation}
\label{air3} \cO_{\di \ddi}^r(X) \ra \cO_{\di
\ddi}^r(T(X))[\epsilon]
\end{equation}
whose composition with the first projection
$$\begin{array}{ccc}
\cO_{\di \ddi}^r(T(X))[\epsilon] & \stackrel{pr_1}{\ra } &
\cO_{\di \ddi}^r(T(X)) \\ a+\epsilon b & \mapsto & a
\end{array}
$$
is the identity. Composing the morphism (\ref{air3}) with the second projection
$$\begin{array}{ccc}
\cO_{\di \ddi}^r(T(X))[\epsilon] & \stackrel{pr_2}{\ra } &
\cO_{\di \ddi}^r(T(X)) \\ a+\epsilon b & \mapsto & b
\end{array}\, ,
$$
we get $A$-derivations as in (\ref{classsoon}), which agree with each other
as $r$ varies, hence induce an $A$-derivation as in (\ref{air1}).
Note that $\Theta$ restricted to $\cO(X)$ equals $d$. Also the map
$f \mapsto f+ \epsilon \Theta f$ in (\ref{air3})
commutes with $\phip$ and $\dd$ (by universality). In particular
 we get that
 \begin{equation}
 \label{iriplea}
 \begin{array}{rcl}
 \Theta \circ \phip & = & p \cdot \phip \circ \Theta\, ,\\
 \Theta \circ \dd & = & \dd \circ \Theta \, .
 \end{array}
 \end{equation}
Clearly we have the following:

 \begin{proposition}
The mapping $\Theta$ is the unique $A$-derivation
$$\cO^{\infty}_{\di \ddi}(X) \ra \cO^{\infty}_{\di \ddi}(T(X))$$
extending $d$ in {\rm (\ref{irinushe})} and satisfying the
commutation relations in {\rm (\ref{iriplea})}.
\end{proposition}

\begin{corollary}
\label{bibby} For any $A$-derivation $\partial:\cO_X \ra \cO_X$,
there exists a unique derivation
\begin{equation}
\label{parst}
\partial_{\infty}:\cO^{\infty}_{\di \ddi}(X) \ra
\cO^{\infty}_{\di \ddi}(X)\end{equation}
 extending $\partial$ and
satisfying the commutation relations
\begin{equation}
\begin{array}{rcl}
\partial_{\infty} \circ \phip & = & p \cdot \phip \circ
\partial_{\infty}\, ,\\
\partial_{\infty} \circ \dd & = & \dd \circ
\partial_{\infty}\, .
\end{array}
\end{equation}
 \end{corollary}

The first of the two conditions above says that
$\partial_{\infty}$ is a $\d$-derivation. The derivation
$\partial_{\infty}$ extends the derivation $\partial_*$ on
$\cO^{\infty}_p(X)$ in \cite{book}, Definition 3.40.

{\it Proof}. The uniqueness is clear. For the proof of existence, we may
assume that $X$ is affine. By the universality property of the K\"{a}hler
differentials, $\partial$ induces an $\cO(X)$-module map
$$\<\partial,\ \>:\Omega_{\cO(X)/A} \ra \cO(X)$$
such that $\<\partial,df\>=\partial f$ for all $f \in \cO(X)$. By
the universality property of the symmetric algebra, we get an
induced $\cO(X)$-algebra map
$$\cO(T(X)) \ra \cO(X)\, .$$
Composing this map with the inclusion $\cO(X) \subset
\cO^{\infty}_{\di \ddi}(X)$, we get a homomorphism
$$\cO(T(X)) \ra \cO^{\infty}_{\di \ddi}(X)\, .$$
By the universality property of $\D$-jet spaces, we get a $\D$-ring
homomorphism
$$\<\partial,\ \>_{\infty}:\cO^{\infty}_{\di \ddi}(T(X)) \ra
\cO^{\infty}_{\di \ddi}(X)\, ,$$
which is an $\cO^{\infty}_{\di
\ddi}(X)$- algebra homomorphism. Composing the latter with the
Fr\'{e}chet derivative
$$\Theta:\cO_{\di \ddi}^{\infty}(X) \ra \cO_{\di
\ddi}^{\infty}(T(X))\, ,$$
we get an $A$-derivation $\partial_{\infty}$ as in (\ref{parst}),
$$\partial_{\infty} f=\<\partial,\Theta f\>_{\infty}\, , \quad
f \in \cO^{\infty}_{\di \ddi}(X)\, ,$$
which clearly satisfies all the required properties.
\qed

\begin{definition}
The derivation $\partial_{\infty}$ in Corollary \ref{bibby} is the
{\it prolongation} of $\partial$. An {\it infinitesimal $\infty$-symmetry} of
an element $f \in \cO^{\infty}_{\di \ddi}(X)$ is an
$A$-derivation $\partial:\cO_X \ra \cO_X$ such that
$\partial_{\infty} f=0$.
\end{definition}

The concept introduced above is an analogue of that of an
infinitesimal symmetry in differential geometry \cite{olver}. Note
that the concept of infinitesimal $\infty$-symmetry does not
reduce in the ordinary arithmetic case to the concept in
\cite{book}, Definition 3.48. Again, the difference lies in the
powers of $p$ occurring in these two definitions.

Next we examine the image of the Fr\'{e}chet derivative. For any
affine $X$, let us define
\begin{equation}
\label{air4} \cO_{\di \ddi}^r(T(X))^+:=\sum_{i+j=0}^r \sum_{f \in
\cO(X)} \cO_{\di \ddi}^r(X) \cdot \phip^i \dd^j (df) \subset
\cO^r(T(X))\, .
\end{equation}
Note that if $\omega_1,\ldots,\omega_d$ is a basis of $\Omega_{X/A}$
(for instance, if  $T_1,\ldots,T_d \in \cO(X)$ is a system of \'{e}tale
coordinates and $\omega_i=dT_i$), then $\cO^r(T(X))^+$ is a free
$\cO_{\di \ddi}^r(X)$-module with basis
$\{\phip^i \dd^j \omega_m\mid_{1 \leq m \leq d, 0 \leq i+j \leq r}\}$.

\begin{proposition}
\label{air5} The image of the map $\Theta:\cO_{\di \ddi}^r(X) \ra
\cO_{\di \ddi}^r(T(X))$ is contained in $\cO_{\di
\ddi}^r(T(X))^+$.
\end{proposition}

{\it Proof}.
Since $\Theta$ is an $A$-derivation and $\cO_{\di \ddi}^r(X)$ is
topologically generated by $A$ and elements of the form $\d^i
\dd^j f$  with $i+j  \leq r$, $f \in \cO(X)$, it suffices to check
that for any such $i,j, f$, we have that
$\Theta (\d^i \dd^j  f) \in \cO_{\di \ddi}^r(T(X))^+$.

But
$$\begin{array}{rcl}
\Theta (\d^i \dd^j  f) & = & pr_2(\d^i \dd^j (f+\epsilon df))\\
 & = & pr_2(\d^i (\dd^j f+\epsilon \dd^j(df)))\, .
\end{array}
$$
By \cite{book}, Lemma 3.34, the latter has the form
$$\sum_{l=0}^i \Lambda_{il}(\dd^j f,\d \dd^j f,\ldots,\d^{i-1}\dd^j
f) \phip^l\dd^j(df)$$
where $\Lambda_{il}$ are polynomials with $A$-coefficients, and we are done.
\qed

\begin{remark}
\label{dablr} The Fr\'{e}chet derivative is functorial with
respect to pull back in the following sense. Let $f:\hat{Y} \ra
\hat{X}$ be a morphism of $p$-adic formal schemes between the
$p$-adic completions of two schemes, $X$ and $Y$, of finite type
over $A$. Cf. Remark \ref{ginge}. Then we have a natural
commutative diagram
\begin{equation}
\label{shaa}
\begin{array}{ccc}
\cO^r_{\di \ddi} (X) & \stackrel{f^*}{\ra} & \cO^r_{\di \ddi} (Y)
\vspace{1mm} \\
\Theta \downarrow & \  & \downarrow \Theta\\
\cO^r_{\di \ddi} (T(X)) & \stackrel{f^*}{\ra} & \cO^r_{\di \ddi}
(T(Y))
\end{array} \, .
\end{equation}
\end{remark}

\subsection{Fr\'{e}chet symbol}
We assume in what follows that $X$ has relative dimension $1$ over
$A$. Let $f \in \cO^r_{\di \ddi}(X)$, and let $\omega$ be a basis
of $\Omega_{X/A}$. Then we can write
\begin{equation}
\label{numst} \Theta f=\theta_{f,\omega}(\phip,\dd) \omega=\sum
a_{ij} \phip^i \dd^j \omega \, ,
\end{equation}
where $\theta_{f,\omega}$ is a polynomial
$$\theta_{f,\omega} =\theta_{f,\omega} (\xi_p,\xi_q)=\sum a_{ij}
\xi_p^i \xi_q^j \in \cO^r_{\di \ddi}(X)[\xi_p,\xi_q]\, .$$
If, in addition, we have given an $A$-point $P \in X(A)$, by the universality
property of $\D$-jet spaces, we obtain a naturally induced lift $P^r \in
J^r_{\di \ddi}(X)(A)$ of $P$.
For any $g \in \cO^r_{\di
\ddi}(X)=\cO(J^r_{\di \ddi}(X))$, we denote by $g(P) \in A$ the
image of $g$ under the ``evaluation'' homomorphism $\cO^r_{\di
\ddi}(X) \ra A$ induced by $P^r$. Then we may consider the
polynomial
$$\theta_{f,\omega,P}(\xi_p,\xi_q)=\sum a_{ij}(P)
 \xi_p^i \xi_q^j\in A[\xi_p,\xi_q]\, .
$$

\begin{definition}
\label{centeicit} The polynomial $\theta_{f,\omega}$ is the {\it
Fr\'{e}chet symbol} of $f$ with respect to $\omega$. The
polynomial $\theta_{f,\omega,P}$ is the {\it Fr\'{e}chet symbol}
of $f$ at $P$ with respect to $\omega$.
\end{definition}
\medskip

The polynomials $\theta_{f,\omega}$ and $\theta_{f,\omega,P}$
have a certain  covariance property with respect to $\omega$,
which we explain next.

\begin{definition}
Let $B$ be any  $\D$-ring in which $p$ is a non-zero divisor. We
define the (right) action of the group $B^{\times}$  on the ring
$B[\xi_p,\xi_q]$ as follows.
 If $b \in B^{\times}$ and $\theta \in B[\xi_p,\xi_q]$, then
$\theta \cdot b \in B[\xi_p,\xi_q]$ is the unique polynomial with
the property that for any $\D$-ring extension $C$ of $B$, and any
$x \in C$, we have
$$((\theta \cdot b)(\phip,\dd))(x)=\theta(\phip,\dd) \cdot (bx)\, .$$
\end{definition}

Now taking $B$ to be either $\cO^{\infty}_{\di \ddi}(X)$ or $A$, it
is easy to see that \begin{equation} \label{iksplod}
\begin{array}{rcl}
\theta_{f,g \omega} & = & \theta_{f,\omega} \cdot g^{-1} \\
\theta_{f,g \omega,P} & = & \theta_{f,\omega,P} \cdot (g(P))^{-1}
\end{array}.\end{equation}
In particular, $\theta_{f,\omega,P}$ only depends on the image
$\omega(P)$ of $\omega$ in the {\it cotangent space of $X$ at
$P$}, $\Omega_{X/A} \otimes_{\cO(X),P} A$, where $A$ here is viewed as an
$\cO(X)$-module via the evaluation map $P:\cO(X) \ra A$.

\begin{remark}
The proof of Corollary \ref{bibby} and (\ref{numst}) show
that if $\omega$ is a basis of $\Omega_{X/A}$ and $\partial:\cO_X
\ra \cO_X$ is an $A$-derivation, then for any $f \in
\cO^{\infty}_{\di \ddi}(X)$ we have that
\begin{equation}
\label{illle}
\partial_{\infty} f=\<
\partial,\theta_{f,\omega}(\phip,\dd)\omega\>_{\infty}=
\theta_{f,\omega}(\phip,\dd)(\<\partial,\omega\>),
\end{equation}
where $\theta_{f,\omega}$ is the Fr\'{e}chet symbol of $f$ with
respect to $\omega$.
\end{remark}
\medskip

Our next proposition relates Fr\'{e}chet symbols of
$\D$-characters to the Picard-Fuchs symbol; cf. Definition
\ref{piiccfu}. This will be useful later, as in the applications
to come, the Picard-Fuchs symbols will be easy to compute.

We let $G$ be either $\bG_a$, or $\bG_m$, or an elliptic curve over $A$.

\begin{proposition}
\label{inquests}
 Assume that $\omega$ is an $A$-basis for  the
invariant $1$-forms on $G$, and let $\psi$ be a $\D$-character of
$G$. Let $\theta=\theta_{\psi,\omega} \in \cO_{\di
\ddi}^{\infty}[\xi_p,\xi_q]$ be the Fr\'{e}chet symbol of $\psi$
with respect to $\omega$,  and let $\sigma \in A[\xi_p,\xi_q]$ be
the Picard-Fuchs symbol of $\psi$ with respect to an \'{e}tale
coordinate $T$ at the origin $0$ such that $\omega(0)=dT$. Then
$$\sigma^{(p)}=p \theta\, .$$
\end{proposition}
\medskip

\begin{remark}
\begin{enumerate}
\item[1)] The relation $\sigma^{(p)}=p \theta$ reads
$\sigma(p\xi_p,\xi_q)=p\theta(\xi_p,\xi_q)$.
\item[2)] The expression above implies that, in particular, $\theta \in
A[\xi_p,\xi_q]$, and $\Theta \psi=\frac{1}{p} \sigma(p \phip,\dd) \omega$.
\item[3)] If ${\rm Lie}(G)$ denotes the $A$-module of invariant $A$-derivations
of $\cO_G$, and if $\partial \in {\rm Lie}(G)$, by (\ref{illle}) we then have
that
\begin{equation}
\label{durre}
\partial_{\infty} \psi=\frac{1}{p} \sigma(p
\phip,\dd)(\<\partial,\omega\>) \in A\, .
\end{equation}
In particular, $\partial$ is an infinitesimal $\infty$-symmetry of
$\psi$ if, and only if, the element $u:=\<\partial,\omega\> \in
A=\bG_a(A)$ is a solution to the $\D$-character of $\bG_a={\rm Spec}\,
A[y]$ defined by $\frac{1}{p} \sigma(p \phip,\dd) y \in A\{y\}$.
\end{enumerate}
\end{remark}

{\it Proof of Proposition \ref{inquests}}.
Consider the diagram in (\ref{shaa}) with $X=G$, $Y=\bG_a$,
and $f=e(pT)$. By the definition of the Picard-Fuchs symbol, if
$\sigma=\sum a_{ij} \xi_p^i \xi_q^j$ then
$$e(pT)^* \psi=\sum a_{ij} \phip^i \dd^j T\, .$$
Applying the Fr\'{e}chet derivative to this identity,
and using the commutativity of the diagram (\ref{shaa}), we get that
\begin{equation}
\label{eric} e(pT)^* \Theta \psi  =  \Theta(e(pT)^*
\psi)=\Theta(\sum a_{ij} \phip^i \dd^j T)=\sum a_{ij} p^i \phip^i
\dd^j (dT)\, .
\end{equation}
On the other hand, we may write
$\theta_{\psi,\omega}=\sum b_{ij} \xi_p^i \xi_q^j$,
with $b_{ij} \in \cO^r_{\di \ddi}(G)$. By the definition of the
Fr\'{e}chet symbol, we have
$\Theta \psi=\sum b_{ij} \phip^i \dd^j \omega$, and since
$\omega(0)=dT$, we have that $\omega=d l(T)$, where $l(T)$ is the
compositional inverse of $e(T)$. Therefore, $e(pT)^* \omega=
d(l(e(pT)))=pdT$. We obtain that
\begin{equation}
\label{erric} e(pT)^* \Theta \psi=e(pT)^* \left( \sum b_{ij}
\phip^i \dd^j \omega \right)=p \sum b_{ij} \phip^i \dd^j (dT)\, .
\end{equation}
The identities (\ref{eric}) and (\ref{erric}) imply that
$\left( \sum a_{ij}p^i \phip^i \dd^j \right) (dT)=
p \left( \sum b_{ij} \phip^i \dd^j \right) (dT)$, and this finishes the
proof. \qed

\subsection{Euler-Lagrange equations}
The Fr\'{e}chet symbol can be also used to introduce an
Euler-Lagrange formalism in our setting.
In order to explain this, we begin by making the following definition.

\begin{definition}
Let $B$ be any $\D$-ring in which $p$ is a non-zero divisor. For
any $r \geq 1$, we denote by $B[\xi_p,\xi_q]_r$ the submodule of the
polynomial ring $B[\xi_p,\xi_q]$ consisting of all polynomials of
degree $\leq r$. The {\it adjunction map}
$$Ad^r:B[\xi_p,\xi_q]_r \ra B[\xi_p,\xi_q] [1/p]$$
is defined by
\begin{equation}  Ad^r(\sum_{ij}
b_{ij} \xi_p^i \xi_q^j):= \sum_{ij} (-1)^j p^{-ij} \xi_p^{r-i}
\xi_q^j \cdot b_{ij}\, .
\end{equation}
This map induces an {\it
adjunction map}
$$ad^r:B[\xi_p,\xi_q]_r \ra B [1/p]$$
by
$$ad^r(Q):=(Ad^r(Q)(\phip,\dd))(1)\, ,$$
which is explicitly given by
\begin{equation} \label{defadsus} ad^r(\sum_{ij}
b_{ij} \xi_p^i \xi_q^j):= \sum_{ij} (-1)^j p^{-ij} \phip^{r-i}
\dd^j b_{ij}\, .
\end{equation}
For any $Q \in B[\xi_p,\xi_q]_r$, we have
$$ad^{r+1}(Q)=(ad^r(Q))^{\phi}\, .$$
\end{definition}
\medskip

This adjunction map is a hybrid between the ``familiar''
adjunction map for usual linear differential operators with
respect to $\dd$ (as encountered after the usual integration by parts
argument in the calculus of variations), and the adjunction
map for $\d$-characters (as defined in \cite{book}). The
definition of the adjunction map above might seem ad hoc, and
somewhat complicated, but it is justified by the following
covariance property.

\begin{lemma}
\label{laiz} For any $Q \in B[\xi_p,\xi_q]_r$ and $b \in B$, we have
that
$$
ad^r(Q \cdot b)=b^{\phi^r} \cdot ad^r(Q)\, .
$$
\end{lemma}

{\it Proof}. This follows from a direct computation. \qed
\medskip

The mapping $ad^r$ is related (but does not coincide) with the mapping
$ad_r$ defined in \cite{book}, p. 92. Indeed, the powers of $p$ in
the definitions of $ad^r$ and $ad_r$ are different.

Returning to our geometric setting, let $X/A$ be a smooth
scheme of relative dimension $1$, let $f \in \cO^r_{\di \ddi}(X)$,
and let $\partial:\cO_X \ra \cO_X$ be an $A$-derivation. We define an element
$\epsilon^r_{f,\partial} \in \cO^{\infty}_{\di \ddi}(X)[1/p]$ as follows.
Let us assume first that $X$ is affine and $\Omega_{X/A}$ is free with basis
$\omega$. We let $B$ in the discussion above be the $\D$-ring
$\cO^{\infty}_{\di \ddi}(X)$. We consider the Fr\'{e}chet symbol
$\theta_{f,\omega} \in \cO^{\infty}_{\di \ddi}(X)[\xi_p,\xi_q]_r$, and its
image under $ad^r$, $ad^r(\theta_{f,\omega}) \in \cO^{\infty}_{\di \ddi}(X)
[1/p]$. Then we set
\begin{equation}
\epsilon^r_{f,\partial}:=\<\partial,\omega\>^{\phi^r} \cdot
ad^r(\theta_{f,\omega})\in \cO^{\infty}_{\di \ddi}(X) [1/p]\, .
\end{equation}
By Lemma \ref{laiz} and (\ref{iksplod}), $\epsilon^r_{f,\partial}$ does not
depend on the choice of $\omega$. Therefore, this definition globalizes to one
in the case when $X$ is not necessarily affine and $\Omega_{X/A}$ is not
necessarily free.

\begin{definition}
We say that the element $\epsilon^r_{f,\partial}$ is the {\it Euler-Lagrange
equation} attached to the {\it Lagrangian} $f$ and the {\it vector field}
$\partial$.
\end{definition}

\begin{definition}
An {\it energy function} of order $r$ on $G$ is a $\D$-morphism
$H:G \ra \bA^1$ that can be written as
\begin{equation}
\label{biir}H=\sum_{ij} h_{ij} \psi_i \psi_j \, ,
\end{equation}
where $h_{ij} \in A$ and $\psi_i$ are $\D$-characters of $G$ of order
$r$.
\end{definition}

\begin{proposition}
If $H$ is an energy function of order $r$ on $G$ and
$\partial:\cO_G \ra \cO_G$ is an $A$-derivation that constitutes a basis
of ${\rm Lie}(G)$, then the Euler-Lagrange equation
$\epsilon^r_{H,\partial}$ is a $K$-multiple of a $\D$-character.
\end{proposition}

{\it Proof}. Let us assume that $H$ is an in (\ref{biir}), and let $\omega$
be an $A$-basis of the invariant $1$-forms on $G$.  Let
\begin{equation}
\label{altazi} \theta_{\psi_i,\omega}=\sum_{mn} a_{imn} \xi_p^m
\xi_q^n\, .
\end{equation}
Then
$$\begin{array}{rcl}
\Theta H & = & \sum_{ij} h_{ij} ( \psi_i \Theta \psi_j+\psi_j
\Theta \psi_i)\\
 & = & \sum_{mnij}  [h_{ij}(a_{jmn}\psi_i+a_{imn} \psi_j)]\phip^m \dd^n
\omega \, .
\end{array}
$$
Hence
\begin{equation}
\label{ele} \epsilon^r_{H,\partial}=\<\partial,\omega\>^{\phi^r}
\cdot \sum_{mnij} (-1)^n p^{-mn} \phip^{r-m} \dd^n[h_{ij}(a_{jmn}
\psi_i+a_{imn} \psi_j)]\, ,
\end{equation} and we are done.
\qed

\begin{definition}
A {\it boundary element} in $\cO_{\di \ddi}^{\infty}(X) \otimes K$
is an element of the form
$\dd a +\phip b -b$, for some $a,b \in \cO_{\di \ddi}^{\infty}(X) \otimes K$.
\end{definition}
\medskip

The following can be interpreted as an analogue of Noether's
Theorem in mechanics \cite{olver}. It is a hybrid between the
usual Noether Theorem and the ``arithmetic Noether Theorem'' in
\cite{book}. We state it in the affine case.

\begin{proposition}
Let $X$ be an affine smooth scheme over $A$ of relative dimension
$1$, and let $\partial:\cO_X \ra \cO_X$ be an $A$-derivation that is
a basis for the $\cO(X)$-module of all $A$-derivations. Then,
for any $f \in \cO^r_{\di \ddi}(X)$ we that
$\epsilon^r_{f,\partial}-\partial_{\infty}f \in
\cO_{\di \ddi}^{\infty}(X) \otimes K$ is a boundary element. In
particular, if $\partial$ is an infinitesimal $\infty$-symmetry
of $f$, then $\epsilon^r_{f,\partial}$ is a boundary element.
\end{proposition}

{\it Proof}. Let $\omega$ be a basis of $\Omega_{X/A}$,  and let
$\theta_{f,\omega}=\sum b_{ij} \xi_p^i \xi_q^j$ for $b_{ij} \in
\cO^{\infty}_{\di \ddi}(X)$. We set
$v:=\<\partial,\omega\> \in \cO(X)^{\times}$ (in what follows we
might assume that $v=1$, but this would not simplify the computation), and
denote by $\equiv$ the congruence relation in $\cO^{\infty}_{\di
\ddi}(X) \otimes K$ modulo the group of boundary elements.
Then we have
$$
\begin{array}{rcl}
\partial_{\infty}f & = & \sum b_{ij} \phip^i \dd^j v = \sum
b_{ij}p^{-ij} \dd^j \phip^i v\\
 & \equiv & \sum (-1)^j p^{-ij} (\dd^j b_{ij})(\phip^i v)\, .
\end{array}
$$
We also we have that
$$
\begin{array}{rcl}
\epsilon^r_{f,\partial} & = & \sum (-1)^j p^{-ij} (\phip^r
v)(\phip^{r-i}\dd^j b_{ij})\\
 & = & \sum (-1)^j p^{-ij} \phip^{r-i}[(\dd^j
b_{ij})(\phip^i v)]\, .
\end{array}
$$
Consequently,
$$
\epsilon^r_{f,\partial}-\partial_{\infty} f \equiv  \sum (-1)^j
p^{-ij} (\phip^{r-i}-1)[(\dd^j b_{ij})(\phip^i v)] \equiv
0
$$
because $\phip^k-1$ is divisible by $\phip-1$ in $A[\phip]$.
\qed

\section{Additive group}
\setcounter{theorem}{0}

In this section we prove our main results about $\D$-characters and
their space of solutions in the case where $G$ is the additive group.

 Let $\bG_a={\rm Spec}\, A[y]$ be the additive group over our fixed 
$\D$-ring $A$. We equip $\bG_a$ with the invariant $1$-form
$$\omega:=dy\, .$$

\begin{proposition}
\label{09876}
 The $A$-module $\bX_{\di \ddi}^r(\bG_a)$ of
$\D$-characters of order $r$ on $\bG_a$ is free with basis
$$\{\phip^i \dd^j y \mid_{0 \leq i+j
\leq r}\}\, .$$
Hence the $A[\phip,\dd]$-module $\bX_{\di
\ddi}^{\infty}(\bG_a)$ of $\D$-characters of $\bG_a$ is free of
rank one with basis $y$.
\end{proposition}

{\it Proof}. Same argument as in the proof of Lemma \ref{xxzz}. \qed

Throughout the rest of this section, we let $A=R$. In particular,
the notation and discussion in Example \ref{cucurigu} applies to
$\bG_a$ over $R$.

By Proposition \ref{09876}, any
 $\D$-character of $\bG_a$ can be written uniquely as
\begin{equation}
\label{psia} \psi_a:=\psi_a^{\mu} :=\mu(\phip,\dd)y \in
R[y,Dy,\ldots,D^n y]\h,
\end{equation}
where \[\mu(\xi_p,\xi_q)=\sum \mu_{ij} \xi_p^i\xi_q^j
 \in R[\xi_p,\xi_q]\] is a polynomial. Note that the
Picard-Fuchs symbol $\sigma(\xi_p,\xi_q)$ of $\psi_a$ with respect
to the \'{e}tale coordinate $T=y$ is given by
\[\sigma(\xi_p,\xi_q)=p\mu(\xi_p,\xi_q).\]
The Fr\'{e}chet symbol of $\psi_a$ with respect to $\omega=dy$ is
\[\theta(\xi_p,\xi_q)=\mu(p\xi_p,\xi_q).\]

\begin{definition}
We say that $\mu(\xi_p,\xi_q)$ is the {\it characteristic
polynomial} of the character $\psi_a$. We say that the
$\D$-character $\psi_a$ is {\it non-degenerate} if $\mu(0,0) \in
R^{\times}$. Given a non-degenerate character $\psi_a$, we say
that ${\kappa} \in \bZ$ is a {\it characteristic integer} of
$\psi_a$  if $\mu(0,\kappa)=0$. (Note that any characteristic
integer of a non-degenerate $\D-$ character $\psi_a$ must be coprime to
$p$.) We say that $\kappa \in \bZ$ is a {\it totally
non-characteristic} integer if $\kappa \not\equiv 0$ mod $p$ and
$\mu(0,\kappa) \not\equiv 0$ mod $p$. We denote by ${\mathcal K}$
the set of all characteristic integers of $\psi_a$, and set
$\cK_{\pm}:=\cK \cap \bZ_{\pm}$. We denote by $\cK'$ the set of
all totally non-characteristic integers. For all $0 \neq \kappa
\in \bZ$ and $\alpha \in R$, we define the  {\it basic series} of
$\psi_a$ by
\begin{equation}
\begin{array}{rcl}
\label{ua}  u_{a,\kappa,\alpha}:=  u_{a,\kappa,
\alpha}^{\mu} & := & \sum_{n \geq 0} b_{n,\kappa}\phip^n(\alpha
q^{\kappa})\\
\  & = & \sum_{n \geq 0} b_{n,\kappa} \alpha^{\phi^n} q^{\kappa
p^n}\\
\  &  = &  \alpha q^{\kappa}+\cdots \in q^{\pm 1}R[[q^{\pm 1}]],\end{array}
\end{equation}
where $\{b_{n,\kappa}\}_{n \geq 0}$ is the sequence of elements in
$R$ defined inductively by $b_{0,\kappa}=1$,
\begin{equation}
\label{aniinvatza}
b_{n,\kappa}:=-\frac{\sum_{s=1}^n \left(\sum_{j \geq 0}\mu_{sj} \kappa^j
p^{j(n-s)}\right) b_{n-s,\kappa}^{\phi^s}}{\sum_{j \geq 0}
\mu_{0j} \kappa^j p^{jn}}\, ,\;  n \geq 1\, .
\end{equation}
In this last expression, the denominator is congruent to
$\mu(0,0)$ mod $p$ (and is therefore an element of $R^{\times}$).
\end{definition}

The next Lemma intuitively says that the two collections of series
$\{u_{a,\kappa,1}\ |\ \kappa \neq 0\}$ and $\{q^{\kappa}\ |\
\kappa \neq 0\}$ ``diagonalize'' $\psi_a$.

\begin{lemma}
\label{unidul}
For all $0 \neq \kappa \in \bZ$ and $\alpha \in R$
we have $\psi_a u_{a,\kappa,\alpha}=\mu(0,\kappa) \cdot \alpha
q^{\kappa}$.
\end{lemma}

{\it Proof}. The desired result follows from the following computation:
\begin{equation}
\label{inima} \begin{array}{rcl} \psi_a u_{a,\kappa,\alpha} & = &
\sum_{i,j,m \geq 0} \mu_{ij} \phip^i \dd^j(b_{m,\kappa}
\phip^m(\alpha q^{\kappa}))
\vspace{1mm} \\
 & =  & \sum_{i,j,m \geq 0} \mu_{ij} b_{m,\kappa}^{\phi^i}
\kappa^j
p^{jm} \phip^{i+m}(\alpha q^{\kappa}) \vspace{1mm} \\
  & = & \sum_{n \geq 0} \left( \sum_{s \geq 0} \left( \sum_{j
\geq 0} \mu_{sj} \kappa^j p^{j(n-s)} \right)
b_{n-s,\kappa}^{\phi^s} \right)
\phip^n(\alpha q^{\kappa}) \vspace{1mm} \\
\  & = & \left(\sum_{j \geq 0} \mu_{0j} \kappa^j\right) \alpha q^{\kappa}\\
\  & = & \mu(0,\kappa) \alpha q^{\kappa}\, .
\end{array}
\end{equation}
\qed

For any series $u \in R[[q^{\pm 1}]]$  let $\bar{u} \in k[[q^{\pm
1}]]$ denote the reduction of $u$ mod $p$; we recall that
$k=R/pR$. Also we denote by \[\bar{\mu}(\xi_p,\xi_q)=\sum
\bar{\mu}_{ij} \xi_p^i \xi_q^j \in k[\xi_p,\xi_q]\] the reduction
mod $p$ of the characteristic polynomial $\mu(\xi_p,\xi_q)$. It is
convenient to introduce the following terminology.

\begin{definition}
We say that a polynomial $\mu \in R[\xi_p,\xi_q]$ is {\it unmixed}
if $\bar{\mu}_{ij}=0$ for $ij\neq 0$ and there exists $i \neq 0$
such that $\bar{\mu}_{i0}\neq 0$. Equivalently, $\mu$ is unmixed if
$$
\bar{\mu}(\xi_p,\xi_q)=\bar{\mu}(\xi_p,0)+\bar{\mu}(0,\xi_q)-
\bar{\mu}(0,0)\, ,
$$
and
$$
\bar{\mu}(\xi_p,0) \neq \bar{\mu}(0,0)\, .
$$
\end{definition}

\begin{definition}
We say that $S \subset \bZ_+ \backslash \{0\}$ is {\it short} if
$$\frac{\max S}{\min S} < \frac{p}{2}\, .$$
We say that $S \subset \bZ_-\backslash\{0\}$ is {\it short} if the set 
$-S$ is short. In particular, a set $S \subset \bZ \backslash \{0\}$ 
consisting of a single element is short.
\end{definition}

\begin{lemma}
\label{irinusescoala} Let $S \subset \bZ_{\pm} \backslash p\bZ_{\pm}$
be a non-empty finite set of either positive or negative integers, and
let
$$
u:=\sum_{\kappa \in S} u_{a,\kappa,\alpha_{\kappa}}\, ,
$$
where $\alpha_{\kappa} \in R$, not all of them in $pR$.
Then the following hold:
\begin{enumerate}
\item $\bar{u}$ is integral over $k[q^{\pm 1}]$, and the field
extension $k(q) \subset k(q,\bar{u})$ is Abelian with Galois group
killed by $p$. \item If $\mu$ is unmixed, then $\bar{u} \not\in
k(q)$. \item If $\mu$ is unmixed and $S$ is short, then $u$ is
transcendental over $K(q)$.
\end{enumerate}
\end{lemma}

{\it Proof}. Let us assume that $S \subset \bZ_+$. The case $S \subset \bZ_-$
is treated in a similar manner.

We prove assertion 1.
By (\ref{aniinvatza}), for $n \geq 1$, we have that
$$\begin{array}{rcl}
\bar{b}_{n,\kappa} & = & -(\bar{\mu}_{00})^{-1} \left(
\sum_{s=1}^{n-1} \bar{\mu}_{s0} \bar{b}_{n-s,\kappa}^{p^s}
+\sum_{j \geq 0} \bar{\mu}_{nj} \bar{\kappa}^j \right)\\
  & = & -(\bar{\mu}_{00})^{-1} \left( \sum_{s\geq 1}
\bar{\mu}_{s0} \bar{b}_{n-s,\kappa}^{p^s} +\sum_{j \geq 1}
\bar{\mu}_{nj} \bar{\kappa}^j \right)\, .
\end{array}
$$
Also
$$
\bar{u}=\sum_{\kappa \in S} \sum_{n \geq 0} \bar{b}_{n,\kappa}
\bar{\alpha}_{\kappa}^{p^n} q^{\kappa p^n}\, .
$$
Let us denote by $F_p:k \ra
k$ the $p$-th power Frobenius map. Consider the polynomial $g(q) \in
k[q]$ given by
$$\begin{array}{rcl}
g(q) & = & \sum_{\kappa \in S}[\bar{\mu}(F_p,\kappa)-
\bar{\mu}(F_p,0)-\bar{\mu}(0,\kappa)]
(\alpha q^{\kappa})\\
 & = & \sum_{\kappa \in S} \sum_{n \geq 0} \sum
_{j \geq 1} \bar{\mu}_{nj}
 \bar{\kappa}^j\bar{\alpha}_{\kappa}^{p^n} q^{\kappa p^n}-
 \sum_{\kappa \in S} \sum_{j \geq 0} \bar{\mu}_{0j}
 \bar{\kappa}^j \bar{\alpha}_{\kappa}q^{\kappa}\, ,
\end{array}
$$
and define the polynomial $G(t) \in k(q)[t]$ by
$$
G(t):=\sum_{s \geq 0} \bar{\mu}_{s0}t^{p^s}+g(q)\, .
$$
We have that $G(t)$ has an invertible leading coefficient, is
separable ($dG/dt=\bar{\mu}_{00}$), and has $\bar{u}$ as a root:
$$
\begin{array}{rcl}
G(\bar{u}) & = & \sum_{\kappa \in S} \sum_{s \geq 0} \sum_{n \geq
0} \bar{\mu}_{s0}
 \bar{b}_{n,\kappa}^{p^s} \bar{\alpha}_{\kappa}^{p^{n+s}} q^{\kappa p^{n+s}}\\
 &  & + \sum_{\kappa \in S} \sum_{n \geq 0} \sum_{j \geq 1}
\bar{\mu}_{nj} \bar{\kappa}^j \bar{\alpha}_{\kappa}^{p^n}
q^{\kappa p^n}-\sum_{\kappa \in S} \sum_{j\geq 0}
 \bar{\mu}_{0j}\bar{\kappa}^j \bar{\alpha}_{\kappa} q^{\kappa}\\
  & = & \sum_{\kappa \in S} \sum_{n \geq 0} \left( \sum_{s \geq
0} \bar{\mu}_{s0} \bar{b}_{n-s,\kappa}^{p^s}+\sum_{j \geq 1}
 \bar{\mu}_{nj}\bar{\kappa}^j \right)
\bar{\alpha}_{\kappa}^{p^n} q^{\kappa p^n}\\
 &   & - \sum_{\kappa \in S} \sum_{j\geq 0}
 \bar{\mu}_{0j}\bar{\kappa}^j \bar{\alpha}_{\kappa} q^{\kappa}
\\  & = & 0\, .
\end{array}
$$
The difference of any two roots of $G(t)$ is in $k$. Hence,
$k(q,\bar{u})$ is Galois over $k(q)$ and its Galois group $\Sigma$
embeds into $k$ via the map
$$
\begin{array}{rcl}
\Sigma & \ra & k\\
\sigma  & \mapsto & \sigma \bar{u}-\bar{u}\, .
\end{array}
$$

We prove assertion 2.
In this case we have
$$
g(q)=-\bar{\mu}(0,0) \cdot \sum_{\kappa \in S}
\bar{\alpha}_{\kappa} q^{\kappa}\, ,
$$
and $G(t)$ has degree $p^e$ in
$t$ for some $e \geq 1$. Assume that $\bar{u} \in k(q)$. Since
$\bar{u}$ is integral over $k[q]$ it follows that $\bar{u} \in
k[q]$. Let $d$ be the degree of $\bar{u}$. Since
$G(\bar{u})=0$, we see that $d \geq 1$. Since the integers in $S$ are
not divisible by $p$, the coefficient of $q^{dp^e}$ in the
polynomial $G(\bar{u})\in k[q]$ is non-zero, a contradiction.

We prove assertion 3. Let $\kappa_1<\kappa_2< \cdots <\kappa_s$
be the integers in $S$. By our assumption, there is a real
 $\epsilon > 0$ such that $(2+\epsilon)\kappa_s<\kappa_1 p$.
We will show that
\[u(q^{-1})=\sum_{j=1}^s \sum_{n \geq 0} b_{n,\kappa_j}
\alpha_{\kappa_j}^{\phi^n}
q^{-\kappa_j p^n} \in K((q^{-1}))\]
is transcendental over $K(q)$.

For any $\varphi=\sum_{n=-\infty}^f c_n q^n \in K((q^{-1}))$
with $c_f \geq 0$, set $|\varphi|:=e^f$. Also, set $|0|=0$.
Roth's theorem for characteristic zero function fields
\cite{uchi} states that if $\varphi$ is algebraic over $K(q)$
then, for any $\epsilon >0$, the inequality
$$0<|\varphi-P/Q|<|Q|^{-2-\epsilon}$$
has only finitely many solutions $P/Q$ with $P,Q \in K[q]$.
On the other hand, for any $n$, we have
$$\begin{array}{rcl}
0 & < & |u(q^{-1})-\sum_{j=1}^s \sum_{n=0}^N b_{n,\kappa_j}
\alpha_{\kappa_j}^{\phi^n}
q^{-\kappa_j p^n}|\\
  & \leq & e^{-\kappa_1 p^{n+1}}\\
  & < & (e^{\kappa_s p^n})^{-2-\epsilon}\\
  & = & |q^{\kappa_s p^n}|^{-2-\epsilon}\, .
\end{array}
$$
It follows that $u(q^{-1})$ is transcendental over $K(q)$, and this
completes the proof.
 \qed

\begin{remark}
\begin{enumerate}
\item[1)] If $\mu_{ij} \in pR$ for $i \geq 1$, then
$u_{a,\kappa,\alpha}
\in R[q,q^{-1}]\h$.
\item[2)] The mapping
$$\begin{array}{rcc}
R & \ra & R[[q^{\pm 1}]] \\
\alpha & \mapsto & u_{a,\kappa,\alpha}^{\rs}
\end{array}$$
is an injective group homomorphism.  \item[3)] Recall that,
attached to a $\dd$-character $\psi_q$  we defined in
(\ref{vocea}) operators $B_{\kappa}^0$ and $B_{\pm}^0$. Let
$\psi_q$ be, in our case, the identity; hence, for $0 \neq \kappa
\in \bZ$,
$$\begin{array}{ccl}
R[[q^{\pm 1}]] & \stackrel{B_{\kappa}^0}{\ra} & R \\
B_{\kappa}^0\left(\sum a_n q^n \right) & = & a_{\kappa},
\end{array}$$
and
$$\begin{array}{ccl}
R[[q^{\pm 1}]] & \stackrel{B_{\pm}^0}{\ra} & R^{\rho_{\pm}} \\
B_{\pm}^0\left(\sum a_n q^n \right) & = & (a_{\kappa})_{\kappa \in
\cK_{\pm}}.
\end{array}$$
 Note  that if
$\kappa_1, \kappa_2 \in \bZ \backslash p\bZ$, we have
\begin{equation}
\label{meah}
B_{\kappa_1}^0 u_{a,\kappa_2,\alpha}=\alpha \cdot
 \delta_{\kappa_1 \kappa_2},\end{equation}
 where $\delta_{\kappa_1 \kappa_2}$ is the Kronnecker delta.
 \item[4)] We have the
identity
\begin{equation}
\label{sasesi} \dd u_{a,\kappa,\alpha}^{\mu}=\kappa \cdot
u_{a,\kappa,\alpha}^{\mu^{(p)}}\, ,
\end{equation}
where we recall that $\mu^{(p)}(\xi_p,\xi_q):=\mu(p\xi_p,\xi_q)$.
\item[5)] We have the identity
\begin{equation}
\label{picc} u_{a,\kappa,\zeta^{\kappa}
\alpha}(q)=u_{a,\kappa,\alpha}(\zeta q)
\end{equation}
for all $\zeta \in \bmu(R)$; this holds because
$\zeta^{\phi}=\zeta^p$. In particular, if
 $\alpha=\sum_{i=0}^{\infty} m_i \zeta_i^{\kappa}$, $\zeta_i \in \bmu(R)$,
  $m_i \in \bZ$, $v_p(m_i) \ra \infty$, then
$$u_{a,\kappa,\alpha}(q)
= \sum_{i=0}^{\infty} m_i u_{a,\kappa,1}(\zeta_i q)\, .$$ Thus, if
$f \in \bZ \bmu(R) \h$ is such that $(f^{[\kappa]})^{\sh}=\alpha
\in R$, then $u_{a,\kappa,\alpha}$ can be expressed using
convolution:
$$
u_{a,\kappa,\alpha}=f \star
u_{a,\kappa,1}\, .
$$
Note that $\{u_{a,\kappa,\alpha}\ |\ \alpha \in
R\}$ is a $\bZ \bmu(R)\h$-module (under convolution). If $\kappa
\in \bZ \backslash p\bZ$, this module structure comes from  an
$R$-module structure, still denoted by $\star$, by a base change
via the surjective homomorphism
$$\bZ \bmu(R)\h \stackrel{[\kappa]}{\ra} \bZ \bmu(R)\wh
\stackrel{\sh}{\ra} R $$ (see (\ref{notinjjj})),  and the
$R$-module $\{u_{a,\kappa,\alpha}\ |\ \alpha \in R\}$ is free with
basis $u_{a,\kappa,1}$. Hence, for $g \in \bZ \bmu(R)\wh$,
$\beta=(g^{[\kappa]})^{\sh}$, we have that
$$\beta \star u_{a,\kappa,\alpha}=g \star u_{a,\kappa,\alpha}\, .$$
In particular,
$$
u_{a,\kappa,\alpha}  = \alpha \star u_{a,\kappa,1}\, .
$$
\item[6)] We have the following ``rationality'' property: if $\alpha,
\mu_{ij} \in \bZ_{(p)}$, then $u_{a,\kappa,\alpha} \in
\bZ_{(p)}[[q^{\pm 1}]]$.
\end{enumerate}
\end{remark}

\begin{definition}
We say that the mapping
$$\begin{array}{rcc}
R & \ra &  R[[q]] \\ \alpha & \mapsto & v_{\alpha}
\end{array}$$
is a {\it pseudo $\d$-polynomial map} if for any integer $n \geq 0$ there
exists an integer $r_n \geq 0$  and a polynomial $P_n \in R[x_0,x_1,\ldots,
x_{r_n}]$ such that, for all $\alpha \in R$ we have that
\begin{equation}
v_{\alpha}=\sum P_n(\alpha,\d \alpha,\ldots,\d^{r_n}
\alpha) q^n\, . \label{tempo}
\end{equation}
If $X$ is a scheme over $R[[q]]$, then a map $R \ra X(R[[q]])$ is
said to be a {\it pseudo $\d$-polynomial map} if there exists an
open subscheme $U \subset X$, and a closed embedding $U \subset
\bA^n$ such that the image of $R \ra X(R[[q]])$ is contained in
$U(R[[q]])$, and the  maps $$R \ra U(R[[q]]) \subset
\bA^N(R[[q]])=R[[q]]^N \stackrel{pr_i}{\ra} R[[q]],$$ are pseudo
$\d$-polynomial. Here $pr_i$ are the various projections.

Similarly, a mapping
$$\begin{array}{rcc}
R & \ra &  R[[q^{-1}]] \\ \alpha & \mapsto & v_{\alpha}
\end{array}$$
is said to be a {\it pseudo $\d$-polynomial map} if for
any integer $n \leq 0$ there exists an
integer $r_n \geq 0$ and a polynomial $P_n \in
R[x_0,x_1,\ldots,x_{r_n}]$ such that (\ref{tempo}) holds for all
$\alpha \in R$, and given a scheme $X$ over $R[[q^{-1}]]$, a
{\it pseudo $\d$-polynomial} map $R \ra
X(R[[q^{-1}]])$ is defined as above, with the r\^ole of
$R[[q]]$ now being played by $R[[q^{-1}]]$.
\end{definition}

The prefix {\it pseudo} was included in order to suggest an
analogy with ``differential operators of infinite order.'' This is not
to be confused with the pseudo-differential operators in micro-local
analysis.

\begin{example}
The basic series mappings
$$\begin{array}{rcc}
R & \ra &  R[[q^{\pm 1}]] \\ \alpha & \mapsto & u_{a,\kappa,\alpha}
\end{array}$$
are pseudo $\d$-polynomial maps. In particular, when
interpreted as mapping $R \ra \bG_a(R[[q^{\pm
1}]])$, they are pseudo $\d$-polynomial maps.
\end{example}

The following example is as elementary as they come. More interesting ones 
will be given later on, while studying $\bG_m$ and elliptic curves.

\begin{theorem}
\label{addeq} Let $\psi_a$ be a non-degenerate $\D$-character of
$\bG_a$, and let $\cU_*$ be the corresponding groups of solutions.
Let ${\mathcal K}$ be the set of characteristic integers, and
$u_{a,\kappa,\alpha}$ be the basic series. Then the following
hold:
\begin{enumerate}
\item[1)] If ${\mathcal K} =\emptyset$, then $\cU_{\la}=\cU_{\ra}=\cU_0$.
\item[2)] We have
$$\begin{array}{rcl}
\cU_{\pm 1} & = & \bigoplus_{\kappa \in {\mathcal K}_{\pm}}
\{u_{a,\kappa,\alpha}\, | \; \alpha \in R\}\, ,
\end{array}$$
where $\oplus$ denotes internal direct sum. In particular,
$\cU_{\pm 1}$  are free $R$-modules under convolution, with bases
$\{u_{a,\kappa,1}\, |\; \kappa \in \cK_{\pm}\}$, respectively.
\item[3)] $\cU_{\ra}+\cU_{\la}=\cU_+ +\cU_-$.
\end{enumerate}
\end{theorem}

{\it Proof}. Let $u=\sum_{n=-\infty}^{\infty}a_n q^n$ be either an element of
$R((q))\wh$ or of $R((q^{-1}))\wh$. We express the $\D$-character $\psi_a$ as
$\psi_a=\mu(\phip,\dd)y$, for some polynomial
$\mu(\xi_p,\xi_q)=\sum_{i,j \geq 0}
 \mu_{ij} \xi_p^i \xi_q^j \in R[\xi_p,\xi_q]$.
Then $\psi^{\rs}_a u=0$ if, and only if,
\begin{equation}
\label{tutu} \mu(0,n) a_n+\sum_{j \geq 0} \sum_{i \geq 1} \mu_{ij}
(n/p^i)^j a_{n/p^i}^{\phi^i}=0.
\end{equation}
for all $n \in \bZ$. In this last expression, $a_{n/p^i}=0$ if
$n/p^i \not \in \bZ$. Thus, if $\mu(0,n) \neq 0$ for all $n \in \bZ$, we
derive by induction that $a_n=0$ for all $n \neq 0$, which proves the first
part.

In order to prove 2), we first note that, by Lemma \ref{unidul},
$u_{a,\kappa,\alpha} \in \cU_{\pm 1}$ according as   $\kappa \in
{\mathcal K}_{\pm}$ respectively.
 Now if $u^* \in \cU_{1}$, that is to say, if $u^*=\sum_{n \geq 1} a_n q^n$,
we set
$$u^{**}:=u^*-\sum_{\kappa \in {\mathcal K}_+}u_{a,\kappa,a_{\kappa}}
 \in \cU_{1}\, .$$
 Set  $\rho_+:=\sharp \cK_+$. Using (\ref{tutu}),
 one easily checks that
 the map $B_+^0:\cU_+ \ra R^{\rho_+}$ defined by
 \[B_+^0(\sum a_n q^n)=(a_{\kappa})_{\kappa \in \cK_+}\]
 is injective. On the other hand,
by (\ref{meah}), we have
  $B_+^0 u^{**}=0$.
   Thus, $u^{**}=0$, and $u^{*}=\sum u_{a,\kappa.a_{\kappa}}$. A
similar argument holds for $\cU_{-1}$. This completes the proof of
the second part.

In order to prove 3), let us note that if
$$u=\sum_{n=-\infty}^{\infty} a_n q^n \in \cU_{\ra}\, ,$$
it is then clear that
$$\psi_a\left( \sum_{n<0} a_n q^n \right)=0\, , \; \text{and}\;
 \psi_a\left( \sum_{n\geq 0} a_n q^n
\right)=0\, .$$
Thus, $\sum_{n<0} a_n q^n \in \cU_-$ and $\sum_{n\geq 0}
a_n q^n \in \cU_+$. Therefore, $u \in \cU_-+\cU_+$, and so $\cU_{\ra}
\subset \cU_-+\cU_+$. A similar argument shows that $\cU_{\la} \subset
\cU_-+\cU_+$. \qed

\begin{example}
\label{duclacasino}
\label{rupdi} Let us examine a special case of Theorem
\ref{addeq}. For integers $r,s \geq 1$ and $\lambda \in R^{\times}$, we
consider the $\D$-character
\begin{equation}
\label{gilmore} \psi_a:=(\dd^r+\lambda \phip^s-\lambda)y\, ,
\end{equation}
If $(r,s)$ is any one of the pairs $(1,1), (1,2), (2,1), (2,2)$,
then $\psi_a$ can be viewed as an analogue of the convection
equation, heat equation, sideways heat equation, or wave equation,
respectively. The characteristic polynomial of $\psi_a$ is
$$\mu(\xi_p,\xi_q)=\xi_q^r+\lambda \xi_p^s -\lambda\, .$$
Clearly $\mu$ is unmixed. The characteristic integers are
the integer roots of the equation
$$\xi_q^r-\lambda=0\, .$$
Thus, if $\lambda \not\in \{n^r\ |\ n \in \bZ\}$, there are no characteristic
integers, and $\cU_{\ra}=\cU_{\la}=\cU_0=R^{\phi^s}$.

Assume in what follows that $\lambda=n^r$ for some $n \in \bZ$. For even
$r$, we may assume further that $n>0$. Then $\cK=\{n\}$ for $r$ odd,
and $\cK=\{-n,n\}$ for $r$ even. The basic series for
$\kappa \in \cK$ are
\begin{equation}
\label{maomor} u_{a,\kappa,\alpha}= \sum_{j \geq 0} (-1)^j
\frac{1}{F_j(p^{sr})} \phip^{sj}(\alpha q^{\kappa})\, ,
\end{equation}
where $F_j(x) \in \bZ[x]$ are the polynomials $F_0(x)=1$,
$$ F_j(x):=\prod_{i=1}^j (x^i-1)\, , \; j \geq 1\, .$$
Notice that the integers $F_j(p^{sr})$ have a nice simple
interpretation in terms of flags:
$$F_j(p^{sr})=(p^{sr}-1)^j \cdot \sharp(GL_j(\bF_{p^{sr}})/
B_j(\bF_{p^{sr}}))\, ,$$
where $B_j(\bF_{p^{sr}})$ is the subgroup of $GL_j(\bF_{p^{sr}})$
consisting of all upper triangular matrices.
Also notice that, for $\kappa \in \cK$, we have that
 $\bar{u}_{a,\kappa,\alpha} \in k[[q^{\pm 1}]]$,
the reduction mod $p$ of $u_{a,\kappa,\alpha}$, is given by
$$
\bar{u}_{a,\kappa,\alpha}=\sum_{n \geq 0} \bar{\alpha}^{p^n} q^{\kappa p^n}\, ,
$$
so $\bar{u}_{a,\kappa,\alpha}$ is a root of the Artin-Schreier polynomial
$$
t^p-t+\bar{\alpha} q^{\kappa} \in k(q)[t]\, .
$$

We have the following:
\begin{enumerate}
\item[a)]  For $n>0$ and odd $r$,
$$\begin{array}{ccl}
\cU_{-1} & = & 0\, ,\\
\cU_{1} & = & \{u_{a,n,\alpha}\ |\ \alpha \in R\}\, .
\end{array}
$$
\item[b)] For $n <0$ and odd $r$,
$$\begin{array}{ccl}
\cU_{-1} & = & \{u_{a,n,\alpha}\, |\;\alpha \in R\}\, ,\\
\cU_{1} & = & 0\, .
\end{array}
$$
\item[c)] For $r$ even,
$$\begin{array}{ccl}
\cU_{-1} & = & \{u_{a,-n,\alpha}\, | \; \alpha \in R\}\, , \\
\cU_{1} & = & \{u_{a,n,\alpha}\, | \;\alpha \in R\}\, .
\end{array}
$$
\item[d)] $\cU_{\ra}=\cU_+$, \hspace{2mm} $\cU_{\la}=\cU_-$.
\end{enumerate}

Indeed, (a), (b) and (c) follow directly by Theorem \ref{addeq}. The two
families of solutions in (c) should be viewed as analogues of the two waves
traveling in opposite directions in the case of the classical wave equation.
In contrast to this, we have only one ``wave'' in (a), which is the case
of the ``convection'' equation.

We prove the first equality in (d). The second follows by a similar
argument. Let
$$u=\sum_{n=-\infty}^{\infty} a_n q^n \in \cU_{\ra}\, .$$
It is clear that
$$\psi_a\left( \sum_{n<0} a_n q^n \right)=0\, .$$
By (a), (b) and (c), we must have that
$$\sum_{n<0} a_n q^n=u_{a,-|\kappa|,\alpha}$$
for some $\alpha \in R$. But $a_n \ra 0$ as $n \ra -\infty$, and this
is the case for $u_{a,-|\kappa|,\alpha}$ only when $\alpha=0$ (see
(\ref{maomor})). Thus, $\alpha=0$, and $u \in \cU_+$. \qed
\end{example}

We derive here some consequences of Theorem \ref{addeq}.

\begin{corollary}
Under the hypotheses of Theorem {\rm \ref{addeq}} let $u \in
\cU_{\pm 1}$. Then the following hold:
\begin{enumerate}
\item
 The series
$\bar{u} \in k[[q^{\pm 1}]]$ is integral over $k[q^{\pm 1}]$ and
the field extension $k(q) \subset k(q,\bar{u})$ is Abelian with
Galois group killed by $p$. \item If the characteristic polynomial
of $\psi_a$ is unmixed and $\cK_{\pm}$ is short then $u$ is
transcendental over $K(q)$.
\end{enumerate}
\end{corollary}

{\it Proof}.
This follows directly from Theorem \ref{addeq}
and Lemma \ref{irinusescoala}.
\qed

\begin{corollary}
\label{narrenume}
Under the hypotheses of Theorem {\rm \ref{addeq}} the maps
$B_{\pm}^0:\cU_{\pm 1} \ra R^{\rho_{\pm}}$ are $R$-module
isomorphisms. Furthermore, for any $u \in \cU_{\pm 1}$  we have
$$u=\sum_{\kappa \in {\mathcal K}_{\pm}} (B_{\kappa}^0 u) \star
u_{a,\kappa,1}\, .$$
\end{corollary}

In particular  the ``boundary value problem at $q^{\pm 1}=0$'' is
well posed.

The next Corollary says that the ``boundary value problem at $q
\neq 0$'' is well posed.

\begin{corollary}
\label{cosi} Under the hypotheses of Theorem {\rm \ref{addeq}},
assume further that $\cK_+=\{\kappa \}$. Then for any $q_0 \in
p^{\nu}R^{\times}$ with $\nu \geq 1$, and any $g \in p^{\kappa
\nu}R$, there exists a unique $u \in \cU_{1}^{\rs}$ such that
$u(q_0)=g$.
\end{corollary}

{\it Proof}. We need to show that the map
$$\begin{array}{ccl}
R  & \ra  & p^{\kappa \nu}R \\ \alpha & \mapsto & \sum_{n \geq 0}
b_{n,\kappa} \alpha^{\phi^n}  q_0^{\kappa p^n}
\end{array}
$$
is bijective. This follows by Lemma \ref{ajutator} below. \qed

\begin{lemma}
\label{ajutator} Let $c_0 \in R^{\times}$, $c_1,c_2,c_3,\ldots  \in pR$
and $c_n \ra 0$ $p$-adically as $n \ra \infty$. Then the map
$$\begin{array}{ccl}
R & \ra  & R \\
\alpha & \mapsto & \sum_{n \geq 0} c_n
\alpha^{\phi^n}
\end{array}
$$ is bijective.
\end{lemma}

{\it Proof}. The injectivity is clear. And surjectivity follows by a
Hensel-type argument. \qed

The following Corollary is concerned with the inhomogeneous
equation $\psi_a u=\varphi$.  Recall that for $\varphi
 \in q^{\pm 1}R[[q^{\pm 1}]]$ we define the {\it support}
of $\varphi$ as the set
$\{n ; c_n \neq 0\}$. This notion of support is standard
for series but note that it is not a direct analogue of
the notion of support in real analysis if one pursues the analogy
according to which
 $q$ is an analogue of the exponential of complex time.
Also, for any series $v \in R((q^{\pm 1}))\h$ we denote by
$\bar{v} \in k((q^{\pm 1}))$ the reduction of $v$ mod $p$.

\begin{corollary}
Let $\psi_a$ be a non-degenerate $\D$-character of $\bG_a$ and let
$\varphi \in q^{\pm 1}R[[q^{\pm 1}]]$ be a series whose support
is contained
 in the set $\cK'$ of totally non-characteristic
integers of $\psi_a$. Then the following hold:
\begin{enumerate}
\item The
equation
$\psi_a u=\varphi$ has a unique solution $u \in \bG_a(q^{\pm
 1}R[[q^{\pm 1}]])$
such that $u$ has support disjoint from the set $\cK$ of
characteristic integers. \item If  $\bar{\varphi} \in k[q^{\pm
1}]$ then $\bar{u} \in k[[q^{\pm 1}]]$ is integral over $k[q^{\pm
1}]$ and the field extension $k(q) \subset k(q,\bar{u})$ is
Abelian with Galois group killed by $p$. \item If the
characteristic polynomial of $\psi_a$ is unmixed and the support
of $\varphi$ is short then $u$ is transcendental over $K(q)$.
\end{enumerate}
\end{corollary}

\begin{proof}
The existence in assertion 1 follows from Lemma
\ref{unidul}. Uniqueness follows from Lemma \ref{narrenume}.
Assertions 2 and 3 follows from
\ref{irinusescoala}.
\end{proof}

\begin{remark}
Corollary \ref{cosi} implies that if $\psi_a$ is non-degenerate
and $\cK_+=\{1\}$, for any $q_0 \in pR^{\times}$ the group
homomorphism
$$\begin{array}{ccl}
R & \stackrel{S_{q_0}}{\ra} &  \bG_a(R)=R \\
\alpha & \mapsto &
\frac{1}{p}u_{a,1,\alpha}(q_0)
\end{array}
$$
is an isomorphism. Thus, for any $q_1,q_2 \in pR^{\times}$, we have an
isomorphism
$$S_{q_1,q_2}:=S_{q_2} \circ S_{q_1}^{-1}:\bG_a(R) \ra \bG_a(R)\, .$$
mapping that can be viewed as the ``propagator'' attached to
$\psi_a$. Note that if $\zeta \in \bmu(R)$ and $q_0 \in p
R^{\times}$,  then by (\ref{picc}) we have that
$$S_{\zeta q_0}(\alpha)=\frac{1}{p} u_{a,1,\alpha}(\zeta
q_0)=\frac{1}{p} u_{a,1,\zeta \alpha}(q_0)=S_{q_0}(\zeta \alpha)\, ,$$
so
$$S_{\zeta q_0}=S_{q_0} \circ M_{\zeta}\, ,$$
where $M_{\zeta}:R \ra R$ is the mapping defined by
$M_{\zeta}(\alpha):=\zeta \alpha$.
Thus, for $\zeta_1, \zeta_2 \in \bmu(R)$, we get that
$$S_{\zeta_1 q_0,\zeta_2 q_0}=S_{q_0} \circ M_{\zeta_2/\zeta_1}
\circ S_{q_0}^{-1}\, .$$
In particular,
$$S_{q_0,\zeta_1 \zeta_2 q_0}=S_{q_0,\zeta_2 q_0} \circ
S_{q_0,\zeta_1 q_0}\, .$$
This latter equality can be interpreted as a
(weak)  incarnation of ``Huygens principle'' (\cite{rauch}, p. 104).
\end{remark}

\section{Multiplicative group}
\setcounter{theorem}{0}

In this section we prove our main results about $\D$-characters
and their space of solutions in the case where $G$ is the multiplicative group.

 Let $\bG_m:={\rm Spec}\, A[y,y^{-1}]$ be the multiplicative group over
our fixed $\D$-ring $A$. We equip $\bG_m$ with the invariant
$1$-form
$$\omega:=\frac{dy}{y}\, .$$

Let us consider the $\D$-characters
$$\psi_{\di}, \psi_{\ddi} \in \bX^1_{\di \ddi}(\bG_m) \subset
 A[y,y^{-1}, \d y, \dd y]\wh $$
defined by
$$
\begin{array}{rcl}
\psi_{\di} = \psi_{m,\di} & := & {\displaystyle \frac{1}{p}\log{ \left(
\frac{\phip(y)}{y^p} \right)}= \frac{1}{p}\log{ \left( 1+p \frac{\d
y}{y^p} \right)}=\frac{\d y}{y^p}-\frac{p}{2}\left(\frac{\d
y}{y^p}\right)^2+\cdots }\, , \vspace{2mm} \\
 \psi_{\ddi} = \psi_{m,\ddi} & := & {\displaystyle \dd \log{y}:= \frac{\dd
y}{y}}\, .
\end{array}
$$
Here, if $y=1+T$, then
$$\log{y}:=l(T)=T-\frac{T^2}{2}+\frac{T^3}{3}-\cdots $$
is the logarithm of the formal group of $\bG_m$. We clearly have
$\psi_{\di} \in \bX_{\di}^1(\bG_m)$, and $\psi_{\ddi} \in
\bX_{\ddi}^1(\bG_m)$. The images of $\psi_{\di}$ and $\psi_{\ddi}$
in $A[[T]][\d T, \dd T]\h$ are
$$\begin{array}{rcl}
\psi_{\di} & = & {\displaystyle \frac{1}{p}(\phip-p)l(T)}\, , \vspace{1mm} \\
\psi_{\ddi} & = & \dd l(T)\, .
\end{array}
$$

\begin{lemma}
\label{999} We have that $\dd \psi_{\di}=(\phip-1)\psi_{\ddi}$ in $\bX_{\di
\ddi}^2(\bG_m)$.
\end{lemma}

{\it Proof}. By a direct calculation,
$$ \dd \psi_{\di}=\dd \left( \frac{1}{p} \left( \phip-p\right)
l(T) \right)=(\phip-1)\dd l(T)=(\phip-1) \psi_{\ddi}\, .$$ \qed

\begin{proposition}
\label{roset}
 For each $r \geq 1$, the $L$-vector space $\bX_{\di \ddi}^r(\bG_m)
\otimes_A L$ has basis
$$\{\phip^i \psi_{\di}\mid_{0 \leq i \leq r-1}\} \cup
\{\phip^i \dd^j \psi_{\ddi}\mid_{0 \leq i+j \leq r-1}\}\, .$$
In particular,
$\psi_{\di}$ and $\psi_{\ddi}$ generate the $L[\phip,\dd]$-module
$\bX^{\infty}_{\di \ddi}(\bG_m) \otimes L$.
\end{proposition}

{\it Proof}. There is an exact sequence
of homomorphisms of groups in the category of $p$-adic formal schemes
$$0={\rm Hom}(\hat{\bG}_m,\hat{\bG}_a) \stackrel{\pi_r^*}{\ra}
{\rm Hom}(J^r_{\di \ddi}(\bG_m),\hat{\bG}_a) \stackrel{\rho}{\ra}
{\rm Hom}(N^r,\hat{\bG}_a)\, ,$$
where $\pi_r:J^r_{\di \ddi}(\bG_m) \ra
\hat{\bG}_a$ is the natural projection, $N^r={\rm ker}\, \pi_r$, and
$\rho$ is defined by restriction. We recall that ${\rm Hom}(J^r_{\di
\ddi}(\bG_m),\hat{\bG}_a)$ identifies with the module of $\D$-characters
$\bX^r_{\di \ddi}(\bG_m)$. Looking at the level of Lie algebras, we
see that the rank of ${\rm Hom}(N^r,\hat{\bG}_a)$ over $A$ is at most equal
to its dimension, $r(r+3)/2$. Since $\rho$ is injective, it is enough to
show that the family in the statement of the Proposition is
$A$-linearly independent. Thus, it is enough to show that the image
of this family via the map (\ref{prost}) is $A$-linearly
independent. But
$$
\begin{array}{rcl}
(\phip^i \psi_{\di}) \circ e(pT) & = & {\displaystyle
\phip^i \left( \frac{1}{p} (\phip-p) l(e(pT)) \right)=\phip^{i+1}T-p
\phip^i T}\, ,\vspace{1mm} \\
(\phip^i \dd^j \psi_{\ddi}) \circ e(pT) & = & \phip^i \dd^{j+1}
l(e(pT))=p\phip^i
 \dd^{j+1} T\, ,
 \end{array}
$$
and it is rather clear that these elements are $A$-linearly independent. \qed

From now on, we let $A=R$, hence $L=K$, and we use the notation
and discussion in Example \ref{cucurigu} applied to $\bG_m$ over
$R$.
 We have a natural embedding
$$\begin{array}{ccc}
\iota:q^{\pm 1}R[[q^{\pm 1}]] & \ra & \bG_m(q^{\pm 1}R[[q^{\pm 1}]]) \\
u & \mapsto & \iota(u)=1+u
\end{array}$$

By Proposition \ref{roset}, any $\D$-character of
$\bG_m$ is a $K$-multiple of a $\D$-character of the form
\begin{equation}
\label{rosett} \psi_m:= \nu(\phip,\dd)\psi_{\ddi}+\lambda(\phip)
\psi_{\di}\, ,
\end{equation}
where $\nu(\xi_p,\xi_q) \in R[\xi_p,\xi_q]$, $\lambda(\xi_p) \in
R[\xi_p]$. The Picard-Fuchs symbol of $\psi_m$ with respect to the
\'{e}tale coordinate $T=y-1$ is trivially seen to be
\[\sigma(\xi_p,\xi_q)=p\nu(\xi_p,\xi_q)\xi_q+\lambda(\xi_p)(\xi_p-p).\]
Hence, the Fr\'{e}chet symbol with respect to the invariant form
$\omega=\frac{dy}{y}=\frac{dT}{1+T}$ is given by
$$
\theta(\xi_p,\xi_q)=
\frac{\sigma(p\xi_p,\xi_q)}{p}
=\nu(p\xi_p,\xi_q)\xi_q+\lambda(p\xi_p)(\xi_p-1)\, .$$

\begin{definition}
Let $\psi_m$ be a $\D$-character of $\bG_m$ of the form {\rm
(\ref{rosett})}. We define the {\it characteristic polynomial}
$\mu(\xi_p,\xi_q)$ of $\psi_m$ to be the Fr\'{e}chet symbol
$\theta(\xi_p,\xi_q)$ of $\psi_m$ with respect to $\omega$. We say
that the $\D$-character $\psi_m$ is {\it non-degenerate} if
$\mu(0,0) \in R^{\times}$, or equivalently, if $\lambda(0) \in
R^{\times}$. For a non-degenerate character $\psi_m$ with symbol
$\mu(\xi_p,\xi_q)$, we define the {\it characteristic integers} to
be the  integers $\kappa$ that are solutions of the equation
$\mu(0,\kappa)=0$. Any such characteristic integer $\kappa$ must
be coprime to $p$. We say that $\kappa \in \bZ$ is {\it totally
non-characteristic} if $\kappa \not\equiv 0$ mod $p$ and
$\mu(0,\kappa) \not\equiv 0$ mod $p$. We denote by $\cK$ the set
of all characteristic integers and set $\cK_{\pm}:=\cK \cap
\bZ_{\pm}$. Also we denote by $\cK'$ the set of totally
non-characteristic integers.  For any  $0 \neq \kappa \in \bZ$ and
 $\alpha \in R$, we define the {\it basic series}
\begin{equation}
\label{val}
 u_{m,\kappa,\alpha}:=\exp{ \left(
 \int u_{a,\kappa,\alpha}^{\mu} \frac{dq}{q} \right)}
\in 1+ \alpha \frac{q^{\kappa}}{\kappa} +\cdots \in K[[q^{\pm 1}]]\, ,
\end{equation}
where $u_{a,\kappa,\alpha}^{\mu}$ is as in (\ref{ua}),
$$\exp{(T)}=1+T+\frac{T^2}{2!}+\frac{T^3}{3!}+\cdots\, ,$$
and
$$v \mapsto \int vdq $$
is the usual indefinite integration $K[[q]] \ra qK[[q]]$ or
the integration $q^{-2}K[[q^{-1}]] \ra q^{-1}K[[q^{-1}]]$,
according to the cases $\kappa$ positive or
negative, respectively.
\end{definition}

\begin{example} Consider the $\D$-characters
$$\psi_m:=\dd^{r-1} \psi_{\ddi}+(\lambda_{s-1}\phip^{s-1}+\cdots
+\lambda_1 \phip + \lambda_0) \psi_{\di}\, ,$$ where $\lambda_0
\in R^{\times}$, $\lambda_1,\ldots ,\lambda_{s-1} \in R$. If $(r,s)$
is any one of the pairs $(1,1)$, $(1,2)$, $(2,1)$, $(2,2)$, then
$\psi_m$ can be viewed as an analogue of the convection equation,
heat equation, sideways heat equation, or wave equation
respectively. The characteristic polynomial of $\psi_m$ is
$$\mu(\xi_p,\xi_q)=\xi_q^r+(p^{s-1}\lambda_{s-1}\xi_p^{s-1}+\cdots +p \lambda_1
\xi_p+ \lambda_0)(\xi_p-1)\, ,$$
hence $\mu$ is unmixed, and
the characteristic  integers are the integer roots of the equation
$$\xi_q^r-\lambda_0=0\, .$$
In particular, for
\begin{equation} \label{speetro}
\psi_m:=\dd^{r-1} \psi_q+\lambda \psi_p\, ,
\end{equation}
the characteristic polynomial equals
$$
\mu(\xi_p,\xi_q)=\xi_q^r+\lambda \xi_p-\lambda\, ,
$$
which is the case $s=1$ of the Example \ref{duclacasino}. So in
this case $u_{a,\kappa,\alpha}^{\mu}$ in (\ref{val}) is given by
\begin{equation} \label{spetro}
u_{a,\kappa,\alpha}^{\mu}=\sum_{n \geq 0} (-1)^n
\frac{1}{F_n(p^r)} \phip^n(\alpha q^{\kappa})\, .
\end{equation}
\end{example}

\begin{example}
\label{disom} Let us consider the ``simplest'' possible energy function
$$H=a(\psi_q)^2+2b\psi_p \psi_q +c(\psi_p)^2\, ,$$
for $a,b,c \in R$. Using the obvious equalities
$$\begin{array}{rcl}
\theta_{\psi_q,\omega} & = & \xi_q\, ,\\
\theta_{\psi_p,\omega} & = & \xi_p-1\, ,
\end{array}$$
and (\ref{ele}), we obtain the following formula for the
Euler-Lagrange equation associated to $H$ and the vector field
$\partial=y \partial_y$:
$$\epsilon^1_{H,\partial}=(-2a^{\phi} \phip \dd -2 b^{\phi}
\phip^2+2b) \psi_q+(-2c^{\phi}\phip +2c)\psi_p\, .$$ This
$\D$-character $\epsilon^1_{H,\partial}$ is non-degenerate if, and
only if, $c \in R^{\times}$, and its characteristic  polynomial is
given by
$$\mu(\xi_p,\xi_q)=(-2b^{\phi}p^2\xi_p^2-2a^{\phi}p\xi_p\xi_q+2b)\xi_q+
(-2c^{\phi}p \xi_p+2c)(\xi_p-1)\, .$$
The set of characteristic integers of $\epsilon^1_{H,\partial}$ is
$$\cK=\left\{\frac{c}{b}\right\} \cap \bZ\, .$$
\end{example}

\begin{lemma}
\label{ol} Let $\kappa \in \bZ \backslash p\bZ$. Then, for any
$\alpha \in R$, we have that $u_{m,\kappa,\alpha} \in 1+q^{\pm
1}R[[q^{\pm 1}]]$. Furthermore, the mapping
$$\begin{array}{rcl}
R & \ra &  R[[q^{\pm 1}]]^{\times}=\bG_m(R[[q^{\pm 1}]]) \\
\alpha & \mapsto & u_{m,\kappa,\alpha}
\end{array}
$$
is a pseudo $\d$-polynomial map.
\end{lemma}

{\it Proof}. We let $u:=u_{m,\kappa,\alpha}^{\rs}=\exp{h}$, where
$h:=\int u_{a,\kappa,\alpha}^{\mu} \frac{dq}{q}$. By Dwork's Lemma \ref{dwor},
in order to show that
$u_{m,\kappa,\alpha}^{\rs} \in 1+q^{\pm 1}R[[q^{\pm 1}]]$ it is
enough to prove that
$u^{\phi}/u^p \in 1+pq^{\pm 1}R[[q^{\pm 1}]]$. Since
$u^{\phi}/u^p=\exp((\phi-p)h)$, we just need to show that
\begin{equation}
\label{cevreu} (\phip-p)h \in pq^{\pm 1}R[[q^{\pm 1}]]\, .
\end{equation}

Let us start by observing that
\begin{equation}
\label{tatty}
\dd
h=u_{a,\kappa,\alpha}^{\mu}\, ,
\end{equation}
and consider the $\D$-character of $\bG_a$ defined by
$$\psi_a:=\psi_a^{\mu}:=\mu(\phip,\dd)y\, .$$
By  Lemma \ref{unidul}, we obtain that
$$\begin{array}{rcl}
\dd(p\mu(0,\kappa)\alpha \kappa^{-1} q^{\kappa}) & = & p
\mu(0,\kappa) \alpha q^{\kappa}\\
  & = & p
\psi_a^{\mu} u_{a,\kappa,\alpha}^{\mu}\\
\  & = & p[\nu(p
\phip,\dd)\dd+\lambda(p \phip)(\phip -1)]\dd h \vspace{1mm} \\
 & = & \dd[p \nu(\phip,\dd) \dd+\lambda(\phip)(\phip -p)]h\, .
\end{array}
$$
We deduce that
$$[p \nu(\phip,\dd)\dd+\lambda(\phip)(\phip-p)]h -
p \mu(0,\kappa)\alpha \kappa^{-1}q^{\kappa}
\in q^{\pm 1} K[[q^{\pm
1}]] \cap K=0\, ,$$ and consequently,
$$
\begin{array}{rcl}
\lambda(\phip)(\phip -p)h & = & -p \nu(\phip,\dd)\dd h
+p \mu(0,\kappa)\alpha \kappa^{-1}q^{\kappa}\\
 & = & -p \nu(\phip,\dd) u_{a,\kappa,\alpha}^{\mu}
 +p \mu(0,\kappa)\alpha \kappa^{-1}q^{\kappa}\\
 & \in  & pq^{\pm 1} R[[q^{\pm 1}]]\, .
\end{array}
$$
By Lemma \ref{floricica}, (\ref{cevreu}) then follows.

In order to prove the remaining assertion, let us assume that $\kappa >0$.
The opposite case can be argued similarly. Note that we can find polynomials
$Q_n \in R[z_0,z_1,\ldots,z_{r_n}]$, and integers
$m_n \geq 0$, such that
$$u_{m,\kappa,\alpha}=\sum_{n \geq 1}
p^{-m_n} Q_n(\alpha,\alpha^{\phi},\ldots,\alpha^{\phi^{r_n}})q^n$$
for all $\alpha \in R$. By the part of the Lemma already proven, we
have
$$p^{-m_n}Q_n(\alpha,\alpha^{\phi},\ldots,\alpha^{\phi^{r_n}})\in R$$
for all $\alpha \in R$. We then apply Corollary 3.21 in \cite{book}, and
conclude that there exists a polynomial $P_n \in R[x_0,x_1,\ldots,x_{r_n}]$
such that
$$p^{-m_n}Q_n(\alpha,\alpha^{\phi},\ldots,\alpha^{\phi^{r_n}})=
P_n(\alpha,\alpha^{\phi},\ldots,\alpha^{\phi^{r_n}})$$
for all $\alpha \in R$. Thus, the mapping
$$\begin{array}{rcl}
R & \ra & R[[q]] \\ \alpha & \mapsto & u_{m,\kappa,\alpha}
\end{array}$$
is pseudo $\d$-polynomial. Then
$$\begin{array}{rcl}
R & \ra & R[[q]] \\ \alpha & \mapsto &
\frac{1}{u_{m,\kappa,\alpha}}
\end{array}$$
is pseudo $\d$-polynomial also. By considering the embedding
$$\begin{array}{rcc}
\bG_m & \ra & \bA^2 \\
x & \mapsto & (x,x^{-1})
\end{array}\, ,
$$
the last assertion of the Lemma follows. \qed

As in the additive case, we have the following ``diagonalization''
result.

\begin{lemma}
\label{gainuse}
\label{saduslacaf} For all $\kappa \in \bZ \backslash p\bZ$ and $\alpha \in
R$, we have that
\[\begin{array}{rcl}
\psi_q u_{m,\kappa,\alpha} & = & u_{a,\kappa,\alpha}^{\mu}\\
\psi_m u_{m,\kappa,\alpha} & = & \mu(0,\kappa) \cdot \alpha
\kappa^{-1} q^{\kappa}.\end{array}\]
\end{lemma}

{\it Proof}. The first equality is clear.
Let us assume now that $\kappa >0$. A similar argument can be used
to handle the case $\kappa<0$. By Lemma \ref{999}, we obtain that
$$
\begin{array}{rcl}
\dd \psi^{\rs}_m & = & \dd \left[\nu(\phip,\dd)\psi_{\ddi} +
\lambda(\phip) \psi_{\di} \right]\\  &  = &
[\nu(p\phip,\dd)\dd+\lambda(p\phip)(\phip-1)]\psi_{\ddi}
  \\  & = & \psi_a^{\mu} \psi_{\ddi}\, ,
 \end{array}
$$
where $\psi_a^{\mu}:=\mu(\phip,\dd)y$. In particular, by Lemma
\ref{unidul} we have that
$$\begin{array}{rcl} \dd \psi_m u_{m,\kappa,\alpha}
& = & \psi_a^{\mu} \psi_{\ddi} u_{m,\kappa,\alpha}\\
  & = & \psi_a^{\mu} u_{a,\kappa,\alpha}^{\mu}\\
  & = & \mu(0,\kappa) \alpha q^{\kappa}\\
  & = & \dd(\mu(0,\kappa) \alpha
\kappa^{-1}q^{\kappa})\, ,
\end{array}
$$
and so  $v:=\psi_m
u_{m,\kappa,\alpha} -\mu(0,\kappa)\alpha \kappa^{-1} q^{\kappa}
\in R$.
 On the other hand,
$$v(0)=(\psi_m u_{m,\kappa,\alpha})(0)=\psi_m(u_{m,\kappa,\alpha}(0))=
\psi_m(1)=0\, ,
$$
hence $v=0$, and we are done. \qed

\begin{remark}
\begin{enumerate}
\item[a)] For all $\kappa \in \bZ \backslash p\bZ$ the map
$$\begin{array}{rcc}
R & \ra & R[[q^{\pm 1}]]^{\times} \\
\alpha & \mapsto & u_{m,\kappa,\alpha}^{\rs}
\end{array}$$
is an injective homomorphism.
\item[b)] We recall (see
(\ref{vocea})) the natural group homomorphisms attached to
$\psi_q$,
$$\begin{array}{ccc}
B_{\kappa}^0:R[[q^{\pm 1}]]^{\times} & \ra & R \\
B_{\kappa}^0 u & = & \Gamma_{\kappa} \psi_q u,
\end{array}$$
where $\kappa \in \bZ \backslash p\bZ$ and the homomorphism
$$\begin{array}{ccc}
B_{\pm}^0:R[[q^{\pm 1}]]^{\times} & \ra & R^{\rho_{\pm}} \\
B_{\pm}^0 u & = & (\Gamma_{\kappa} \psi_q u)_{\kappa \in
\cK_{\pm}}.
\end{array}$$
For integers $\kappa_1,\kappa_2 \in \bZ \backslash p\bZ$
 we get that
\begin{equation} \label{pueah} B_{\kappa_1}^0
u_{m,\kappa_2,\alpha} =\Gamma_{\kappa_1}
 u^{\mu}_{a,\kappa_2,\alpha}=
\alpha \cdot \delta_{\kappa_1 \kappa_2}\, .
\end{equation}
 \item[c)] For $\kappa \in \bZ \backslash p\bZ$ we see that
\begin{equation}
\label{pim} u_{m,\kappa,\zeta^{\kappa}
\alpha}(q)=u_{m,\kappa,\alpha}(\zeta q)
\end{equation}
for all $\zeta \in \bmu(R)$. In particular, if
$\alpha=\sum_{i=0}^{\infty} m_i \zeta^{\kappa}_i$, $\zeta_i \in
\bmu(R)$, $m_i \in \bZ$, $v_p(m_i) \ra \infty$, we have
$$u_{m,\kappa,\alpha}(q)
=\prod_{i=0}^{\infty}(u_{m,\kappa,1}(\zeta_i q))^{m_i}\, .$$ (Note
that the right hand side of the equality above converges in the
$(p,q^{\pm 1})$-adic topology of $R[[q^{\pm 1}]]$.) Thus, if $f
\in \bZ \bmu(R)\h$ is such that $(f^{[\kappa]})^{\sh}=\alpha \in
R$, then $u_{m,\kappa,\alpha}$ can be expressed via convolution by
$$u_{m,\kappa,\alpha}=f \star u_{m,\kappa,1}\, .$$
Under convolution, the set $\{u_{m,\kappa,\alpha}\ |\ \alpha \in
R\}$ is a $\bZ \bmu(R)\h$-module, whose module structure
 arises from an $R$-module structure
(still denoted by $\star$) induced by base change via the morphism
$$\bZ \bmu(R)\h \stackrel{[\kappa]}{\ra}
\bZ \bmu(R)\h \stackrel{\sh}{\ra} R\, .$$ (Cf. to
(\ref{notinjjj})). The $R$-module $\{u_{m,\kappa,\alpha}\, | \;
\alpha \in R\}$ is free with basis $u_{m,\kappa,1}$. So, for $g
\in \bZ \bmu(R)\h$, $\beta=(g^{[\kappa]})^{\sh}$, we have that
$$\beta \star u_{m,\kappa,\alpha}=g \star u_{m,\kappa,\alpha}\, ,$$
and in particular
$$
u_{m,\kappa,\alpha} = \alpha \star u_{m,\kappa,1}\, .$$
\item[d)] We have the following ``rationality'' property: if $\alpha \in
\bZ_{(p)}$, then $u_{m,\kappa,\alpha} \in \bZ_{(p)}[[q^{\pm 1}]]$
for $\kappa \in \bZ \backslash p\bZ$.
\end{enumerate}
\end{remark}

\begin{theorem} \label{666}
Let $\psi_m$ be a non-degenerate $\D$-character
of $\bG_m$, $\cU_*$ the corresponding groups of solutions, $\cK$
the set of characteristic integers, and $u_{m,\kappa,\alpha}$ the
basic series. We set
$$
\begin{array}{rclll}
\cU_{tors} & := & \bmu(R) & \subset & \cU_0\\
\cU_{\sim} & := & (R^{\times} \cdot q^{\bZ}) \cap \cU  & \subset  &
\cU_{\da}\, .
\end{array}
$$
Then the following hold:
\begin{enumerate}
\item[1)] $\cU_{\ra}=\cU_{\sim} \cdot \cU_{1}$, and
$\cU_{\la}=\cU_{\sim} \cdot \cU_{-1}$. \item[2)] We have that
$$\begin{array}{ccc}
\cU_{\pm 1} & = & \prod_{\kappa \in \cK_{\pm}}
\{u_{m,\kappa,\alpha}\, |\; \alpha \in R\}\, ,
\end{array}$$
where $\prod$ stands for internal direct product. In particular,
 $\cU_{\pm 1}$  are free $R$-modules under
convolution, with bases $\{u_{m,\kappa,1}\, |\; \kappa \in
\cK_{\pm}\}$ respectively. \item[3)] We have that
$\cU_0=\cU_{tors}$.
\end{enumerate}
\end{theorem}

\begin{remark}
\label{funnyterm}
 The elements $u=bq^n$ of $R^{\times} \cdot
q^{\bZ}$ deserve to be called {\it plane waves} with
$$\begin{array}{ll}
\textit{frequency} & n
\in \bZ\, ,\\
 \textit{wave number} & \gamma:=\psi_{m,\di}(b)={\displaystyle \frac{1}{p}\log
\frac{\phip(b)}{b^p}}\in R\, ,\\
 \textit{wave length}  & 1/\gamma \in K \cup \{\infty\}\, ,\\
  \textit{propagation speed}
 & n /\gamma \in K \cup \{\infty\}\, .
 \end{array}
$$
This can be justified by the analogy behind the chosen terminology. Indeed,
if $u=bq^n$, then
$$n=\frac{\dd u}{u}\, ,$$
and so $n$ is to be interpreted as the ``logarithmic derivative'' of $u$
with respect to $\dd$. Also, we have that
$$\gamma=\frac{1}{p}\log \frac{\phip(u)}{u^p}\, ,$$
and so $\gamma$ is the analogue of a ``logarithmic derivative of $u$ with
respect to $\d$.'' These observations suffice to justify the assertion.

For in the classical theory, the complex valued
function $u(t,x)$ of real variables $t$ and $x$ defined by
$$u(t,x):=a(x)e^{-2 \pi i \nu t}$$
is viewed as a plane wave with  frequency $\nu$, wave number (at
$x$) $\gamma=\gamma(x)=\frac{1}{2 \pi i}
\frac{a'(x)}{a(x)}$, wave length $1/\gamma$, and propagation speed
$\nu /\gamma$. (The standard situation is that in which $a(x)=e^{2
\pi i \gamma x}$, where $\gamma \neq 0$ is a constant.) Notice that
these  $\nu$ and $\gamma$ are equal to
$$-(2 \pi i)^{-1} \frac{\partial_t u}{u}\quad \text{and} \quad
(2 \pi i)^{-1} \frac{\partial_x u}{u}\, ,$$
respectively. Thus, our terminology correspond to these classical
definitions of frequency and wave number provided that $\dd$ and $\d$ are
viewed as analogues of $-(2 \pi i)^{-1}
\partial_t$ and $(2 \pi i)^{-1}\partial_x$, respectively.

The solutions in $\cU_{\sim}$ correspond to those solutions that in the
classical case are obtained by ``separation of variables.''
\end{remark}

{\it Proof of Theorem \ref{666}}. By Lemma \ref{saduslacaf} we
have
\[\psi_m u_{m,\kappa,\alpha}=0\]
for $\kappa \in \cK$.
 We now show that if $u_m \in
\cU_{\ra}$ or $u_m \in \cU_{\la}$, then $u_m \in \cU_{\sim} \cdot
\cU^{\rs}_{\pm 1}$, respectively, which will end the proof of
assertion 1). Indeed, let us just treat the case where $u_m \in
\cU_{\ra}$. The case $u_m \in \cU_{\la}$ follows by a similar
argument.

Note that
$$0=\dd \psi_m^{\rs} u_m=\psi^{\mu}_a \psi_{\ddi} u_m\, ,$$
and so by Theorem \ref{addeq}, there exist $\alpha_1,\ldots,\alpha_{s}
\in R$ and $\lambda \in R$ such that
$$
q  \frac{du_m/dq}{u_m}=\psi_{\ddi} u_m=\lambda+\sum_{i=1}^s u_{a,\kappa_i,
\alpha_i}^{\mu}
=\lambda+  \sum_{i=1}^s q \frac{du_{m,\kappa_i,
\alpha_i}/dq}{u^{\rs}_{m,\kappa_i,
\alpha_i}}\, ,
$$
where $\cK_+=\{\kappa_1,\ldots,\kappa_s\}$. Hence, if
$$v:=\frac{u_m}{u_{m,\kappa_1,\alpha_1}\ldots u_{m,\kappa_s,\alpha_s}}
=\sum_{n=-\infty}^{\infty}
b_nq^n\, ,
$$
then
$$ q\frac{dv/dq}{v}=\lambda \, ,$$
and so $nb_n=\lambda b_n$ for all $n$. Thus, $v=bq^n$ with $n \in \bZ$
and $b \in R^{\times}$. This shows that $u_m \in \cU_{\sim} \cdot
 \cU_{1}$, which completes the  argument.

We now prove assertion 2) for the case where $u_m \in \cU_{1}$.
The case $u_m \in \cU_{-1}$ is treated by a similar argument. We
write
$$u_m=bq^n \prod u_{m,\kappa_i,\alpha_i}$$
as above. Since $u_m(0)=1=\prod u_{m,\kappa_i,\alpha_i}(0)$, it
follows that $n=0$ and $b=1$ so $u_m=\prod
u_{m,\kappa_i,\alpha_i}$. This representation is unique due to
formula (\ref{pueah}).

The last assertion of the theorem is clear. \qed

\begin{corollary}
\label{numagidila}
 Under the hypotheses of Theorem {\rm
\ref{666}}, let $u \in \cU_{\pm}$. Then the following hold:
\begin{enumerate}
\item The series $\overline{\psi_q u} \in k[[q^{\pm 1}]]$  is
integral over $k[q^{\pm 1}]$  and the field extension $k(q)\subset
k(q,\overline{\psi_q u})$ is Abelian with Galois group killed by
$p$. \item If the characteristic polynomial of $\psi_m$ is unmixed
and $\cK_{\pm}$ is short then $u$ is transcendental over $K(q)$.
\end{enumerate}
\end{corollary}

{\it Proof}. Assume $u \in \cU_+$; the case $u \in \cU_-$ is similar.
By  Theorem \ref{666}
we may write
$$
u=a \prod_{\kappa \in \cK_+} u_{m,\kappa,\alpha_{\kappa}}\, ,
$$
with $\alpha_{\kappa} \in R$, $a \in R^{\times}$.
So by Lemma \ref{gainuse},
$$
\psi_q u=\sum_{\kappa \in \cK_+} u^{\mu}_{a,\kappa,\alpha_{\kappa}}\, .
$$
By Lemma  \ref{irinusescoala}, assertion 1 follows. To check assertion
 2 note that if $u$ were algebraic over $K(q)$ then the same would hold
for $\psi_q u=\dd u/u$ and we would get a contradiction  by Lemma
\ref{irinusescoala}.
\qed

\begin{corollary}
\label{niciel}
Under the hypotheses of Theorem {\rm \ref{666}}, the maps
$B_{\pm}^0:\cU_{\pm 1} \ra R^{\rho_{\pm}}$ are $R$-module
isomorphisms. Furthermore, for any $u \in \cU_{\pm 1}$ we have
$$u=\sum_{\kappa \in {\mathcal K}_{\pm}} (B_{\kappa}^0 u) \star
u_{m,\kappa,1}\, .$$
\end{corollary}

So, in particular, the ``boundary value problem at $q^{\pm 1}=0$''
is well posed. We now address the ``boundary value  problem  at $q
\neq 0$'', and the ``limit at $q=0$'' issue:

\begin{corollary}
\label{maisus} Under the hypotheses of Theorem {\rm \ref{666}}, the
following hold:
\begin{enumerate}
\item If $\cK_+=\{\kappa\}$, then for any $q_0 \in
p^{\nu}R^{\times}$ with $\nu \geq 1$, and any $g \in 1+p^{\kappa
\nu}R$, there exists a unique
 $u \in \cU_{1}$ such that $u(q_0)=g$.
\item If $\cK_+=\{1\}$, then for any $q_0 \in pR^{\times}$, and
any $g \in R^{\times}$, there exists a unique
 $u \in \cU_{tors} \cdot \cU_{+}$ such that $u(q_0)=g$.
Furthermore, $u(0)$ is the unique root of unity in $R$ that is
congruent to $u(q_0)$ mod $p$.
\end{enumerate}
\end{corollary}

{\it Proof}. We recall that $\log: 1+p^NR \ra p^N R$ is a bijection for all
$N \geq 1$. So, in order to prove assertion 1), we need to check that there
is a unique $\alpha \in R$ such that
$$\log u^{\rs}_{m,\kappa,\alpha}(q_0)=\log\ g\, .$$
By the proof of Proposition \ref{666}, we have that
$$\log\ u_{m,\kappa,\alpha}^{\rs}(q_0)=\sum_{n \geq 0} c_n
\alpha^{\phi^n}\, ,
$$
where  the $p$-adic valuation of $c_n$ is $\kappa \nu p^n-n$. Hence, by
Lemma \ref{ajutator}, the mapping
$$\begin{array}{rcl}
R & \ra & p^{\kappa \nu}R \\ \alpha & \mapsto &
 \sum c_n \alpha^{\phi^n}
\end{array}
$$
is bijective, and we are done.

In order to prove assertion 2)  let $g=\gamma_0 \cdot v_0$ with
$\gamma_0$ be a root of unity, and $v_0 \in 1+pR$. By assertion
1), there exists $v \in \cU_{1}$ such that $v(q_0)=v_0$. Hence, if
$u:=\gamma_0 \cdot v$, then $u(q_0)=g$, which proves the existence
part. Now, if $\gamma_1 \in \cU_{tors}$, $v_1 \in \cU_{1}$, and
$\gamma_1 v_1 (q_0)=\gamma_0 v_0$, we get that $\gamma_1 \equiv
\gamma_0$ mod $p$, so $\gamma_1=\gamma_0$ and $v_1(q_0)=v_0$. By
the uniqueness in the first part above, we see that $v_1=v$, which
proves the uniqueness part of the assertion.

The claim about $u(0)$ is clear. \qed

The following Corollary is concerned with the inhomogeneous
equation $\psi_m u=\varphi$.

\begin{corollary}
\label{magidila}
 Let $\psi_m$ be a non-degenerate $\D$-character
of $\bG_m$, and let $\varphi \in q^{\pm 1}R[[q^{\pm 1}]]$ be a
series whose support  is contained in the set $\cK'$ of totally
non-characteristic integers of $\psi_m$. Then the following hold:
\begin{enumerate}
\item The equation $\psi_m u=\varphi$ has a unique solution $u \in
\bG_m(q^{\pm 1}R[[q^{\pm 1}]])$ such that the support of $\psi_q
u$  is disjoint from the set $\cK$ of characteristic integers.
\item If $\bar{\varphi} \in k[q^{\pm 1}]$, the series
$\overline{\psi_q u} \in k[[q^{\pm 1}]]$ is integral over
$k[q^{\pm 1}]$, and the field extension $k(q) \subset
k(q,\overline{\psi_q u})$ is Abelian with Galois group killed by
$p$. \item If the characteristic polynomial of $\psi_m$ is unmixed
and the support of $\varphi$ is short then $u$ is transcendental
over $K(q)$.
\end{enumerate}
\end{corollary}

{\it Proof}. Let us assume that $\varphi \in qR[[q]]$. The case
$\varphi \in q^{\pm 1}R[[q^{\pm 1}]]$ is similar. We express $\varphi$
as
$$\varphi=\sum_{\kappa \in S} a_{\kappa}q^{\kappa}\, ,$$
where $S$ is the support of $\varphi$  and set
$\alpha_{\kappa}:=\kappa (\mu(0,\kappa))^{-1} a_{\kappa}$. We
define
$$u:=\prod_{\kappa \in S} u_{m,\kappa,\alpha_{\kappa}}
\in \bG_m(qR[[q]])\, ,
$$
which converges $q$-adically. By Lemma \ref{saduslacaf}, $\psi_m
u=\varphi$. And observe that $u$ is the unique solution in
$\bG_m(qR[[q]])$ subject to the condition that $\psi_qu$ has
support disjoint from $\cK$; cf. Lemma  \ref{niciel}. By Lemma
\ref{saduslacaf},
$$\psi_q u=\sum_{\kappa \in S} u_{a,\kappa,\alpha_{\kappa}}\, ,$$
so we may reach the desired conclusion by Lemma \ref{irinusescoala}.
\qed

\begin{remark}
The second part of Corollary \ref{maisus} implies that, if $\psi_m$
is a non-degenerate $\D$-character of $\bG_m$, and $\cK_+=\{1\}$,
then for any $q_0 \in pR^{\times}$ the group homomorphism
$$\begin{array}{ccl}
S_{q_0}:\bmu(R) \times R & \ra & \bG_m(R)=R^{\times} \\
(\xi,\alpha) & \mapsto & \xi \cdot u_{m,1,\alpha}(q_0)
\end{array}$$
is an isomorphism. Thus, for any $q_1,q_2 \in p R^{\times}$, we have an
isomorphism
$$S_{q_1,q_2}:=S_{q_2} \circ S_{q_1}^{-1}:\bG_m(R) \ra \bG_m(R)\, ,$$
which should be viewed as the ``propagator'' associated to $\psi_m$. Note
that if $\zeta \in \bmu(R)$ and $q_0 \in p R^{\times}$, by (\ref{pim}),
$$S_{\zeta q_0}(\xi,\alpha)=\xi \cdot u_{m,1,\alpha}(\zeta
q_0)=\xi \cdot u_{m,1,\zeta \alpha}(q_0)=S_{q_0}(\xi,\zeta
\alpha)\, ,$$
and so
$$S_{\zeta q_0}=S_{q_0} \circ M_{\zeta}\, ,$$
where
$$\begin{array}{ccl}
M_{\zeta}:\bmu(R) \times R & \ra & \bmu(R) \times R \\
M_{\zeta}(\xi,\alpha) & = & (\xi,\zeta \alpha)
\end{array}\, .
$$
Hence, as for the case of $\bG_a$, if $\zeta_1, \zeta_2 \in \bmu(R)$, we get
that
$$S_{\zeta_1 q_0,\zeta_2 q_0}=S_{q_0} \circ M_{\zeta_2/\zeta_1}
\circ S_{q_0}^{-1}\, .$$
In particular,
$$S_{q_0,\zeta_1 \zeta_2 q_0}=S_{q_0,\zeta_2 q_0} \circ
S_{q_0,\zeta_1 q_0}\, ,$$
which can be interpreted  as a (weak)
``Huygens principle.''
\end{remark}

\section{Elliptic curves}
\setcounter{theorem}{0}

We begin this section by proving some general results about the
space of $\D$-characters of an elliptic curve over a general
$\D$-ring $A$. We determine bases for these spaces; the answer
depends on whether certain Kodaira-Spencer classes (an arithmetic
one, and a geometric one) vanish or not. We then prove our main
results about $\D$-characters and their solution spaces for Tate
curves over $R((q))$, and for elliptic curves over $R$.

\subsection{General results}
 We have as before a fixed $\D$-ring ring $A$ that is $p$-adically complete
Noetherian integral domain of characteristic zero. For technical
reasons, we assume hereafter that the prime $p$ is greater than
$3$. We consider an elliptic curve $E$ over $A$, equipped with an
invertible $1$-form $\omega$. All the definitions and
constructions below apply to the pair $(E,\omega)$.

Since $p \geq 5$, the integer $6$ is invertible in $A$, and therefore we
have unique elements  $a_4,a_6 \in A$ such that $4a_4^3+27a_6^2 \in
A^{\times}$, and such that, up to isomorphism, the
pair $(E, \omega)$ consists of the closure in the projective plane over $A$ of
the affine plane curve
\begin{equation}
\label{cubic} y^2=x^3+a_4x+a_6\, ,
\end{equation}
endowed with the $1$-form
\begin{equation}
\label{cubomega} \omega=\frac{dx}{y}\, .
\end{equation}

We set
\begin{equation}
\label{Nsusr} N^r:={\rm ker}(J^r_{\di \ddi}(E)
\stackrel{\pi_r}{\ra} \hat{E}).
\end{equation}
Exactly as in the theory of algebraic group extensions in the last
Chapter of \cite{serre}, we have an exact sequence
\begin{equation}
\label{papica} 0={\rm Hom}(\hat{E},\hat{\bG}_a) \stackrel{\pi_r^*}{\ra}
{\rm Hom}(J^r_{\di \ddi}(E),\hat{\bG}_a) \stackrel{\rho}{\ra}
{\rm Hom}(N^r,\hat{\bG}_a) \stackrel{\partial}{\ra} H^1(E,\cO)\, ,
\end{equation}
where $\rho$ is the restriction. Let us recall that
${\rm Hom}(J^r_{\di \ddi}(E),\hat{\bG}_a)$ identifies with the
module $\bX_{\di \ddi}^r(E)$ of $\D$-characters. We review and use this
sequence by closely following \cite{book}, pp. 191-196, where the case of
$J^r_{\di}(E)$ was considered.

First of all, note that
the projection $\pi_r:J^r_{\di \ddi}(E) \ra \hat{E}$ has a natural structure of
a principal homogeneous space under $\hat{E} \times N^r \ra \hat{E}$.
By Proposition \ref{local}, $\pi_r$
has local sections in the Zariski topology. We consider a
Zariski open covering and a corresponding family of
sections of $\pi_r$;
\begin{equation}
\label{sec}
E=\bigcup U_{\mu},\ \ s_{\mu}:\hU_{\mu} \ra \pi_r^{-1}(\hU_{\mu})\, .
\end{equation}
We may and assume that there exists an index $\mu_0$ such that the zero section
$0 \in E(A)$ belongs to $U_{\mu_0}(A)$, and $s_{\mu_0}(0)=0$.
Then the morphism
\begin{equation}
\label{noapte adinca}
\tau_{\mu}:\hU_{\mu} \times N^r \ra \pi_r^{-1}(\hU_{\mu})\, ,
\end{equation}
which at the level of $S$-points ($S$ any $A$-algebra) is given by
$(A,B) \mapsto s_{\mu}(A)+B$, is
an isomorphism of principal homogeneous spaces under
$\hU_{\mu} \times N^r \ra \hU_{\mu}$.
The isomorphism $\tau_{\mu_0}$
induces the identity
$0 \times N^r \ra \pi_r^{-1}(0)=N^r$, and if
$T \in \cO(U_{\mu_0})$ is an \'{e}tale coordinate on $U_{\mu_0}$,
then we have an induced $\cO(\hat{U}_{\mu_0})$-isomorphism
\begin{equation}
\label{extradoo}
\begin{array}{c}
\tau_{\mu_0}^*: \cO(\hat{U}_{\mu_0})[T^{(i,j)}\mid_{1 \leq
i+j \leq r} ]\wh
\stackrel{can}{\simeq} \cO^r(\pi_r^{-1}(\hat{U}_{\mu_0}))
\stackrel{\circ \tau_{\mu_0}}{\ra} \vspace{2mm} \\
\stackrel{\circ \tau_{\mu_0}}{\ra}  \cO^r(\hat{U}_{\mu_0}
 \times N^r)=
\cO^r(\hat{U}_{\mu_0})[T^{(i,j)}\mid_{1 \leq i+j \leq r}]\wh \, ,
\end{array}
\end{equation}
where $can$ is the unique isomorphism sending $T^{(i,j)}$ into
$\d^i \dd^j T$. Furthermore, $\tau_{\mu_0}^*$ is the identity modulo $T$.

Let $\hU_{\mu \nu}=\hU_{\mu} \cap \hU_{\nu}$.
 The sections (\ref{sec}) induce maps
$s_{\mu}-s_{\nu} \, : \, \hU_{\mu \nu} \ra J^r_{\di \ddi}(E)$
that clearly factor through maps
\begin{equation}
\label{siminussj}
s_{\mu \nu}:\hU_{\mu \nu}\ra N^r.
\end{equation}
In particular, at the level of $S$-points we have
$(\tau_{\nu}^{-1} \circ \tau_{\mu})(A,B)=(A,s_{\mu \nu}(A)+B)$.

We define the map $\partial$ in (\ref{papica}) by attaching to
any homomorphism $\Theta:N^r \ra \hat{\bG}_a$ the cohomology class
$\varphi:=\partial(\Theta) \in H^1(\hat{E}, \cO)=H^1(E,\cO)$ that is
represented by the cocycle
$\varphi_{\mu \nu}:=\Theta \circ s_{\mu \nu} \in \cO(\hU_{\mu \nu})$.
Then, as in \cite{book}, p. 192, we check that
the sequence (\ref{papica}) is exact.

\begin{proposition} 
\label{indi} The rank of the $A$-module $\bX_{\di \ddi}^r(E)$ is
$r(r+3)/2$ if $\partial=0$, and  $r(r+3)/2-1$ if $\partial \neq 0$.
\end{proposition}

{\it Proof}. Looking at the Lie algebras, it is clear that
${\rm Hom}(N^r,\hat{\bG}_a)$ is an $A$-module of rank at most
${\rm dim}\, N^r=r(r+3)/2$. By Lemma \ref{fu}, the homomorphisms
$L^{[a,b]} \in {\rm Hom}(N^r,\hat{\bG}_a)$, $1 \leq a+b \leq r$,
are $A$-linearly independent. Thus, ${\rm Hom}(N^r,\hat{\bG}_a)$
has rank $r(r+3)/2$ over $A$, and the conclusion follows. \qed

Since $H^1(E,\cO)$ has rank $1$ over $A$, by the exact sequence (\ref{papica})
we conclude that

\begin{corollary}
The $A$-module $\bX_{\di \ddi}^1(E)$ is non-zero.
\end{corollary}

This result contrasts deeply with the ``ode'' story for both, the geometric
\cite{ajm} and arithmetic \cite{char} case. Indeed, for ``general'' $E$ we then
have that
$$\bX_{\di}^1(E)=\bX_{\ddi}^1(E)=0\, .$$

In the sequel, we construct and analyze $\D$-characters of $E$ in more detail.
We use the invertible $1$-form $\omega$ on $E$, and assume that the closed
subscheme of $U_{\mu_0}$ defined by $T$ is the zero section $0$. We
identify the $T$-adic completion of $\cO(U_{\mu_0})$ with the ring of power
series $A[[T]]$, and furthermore, choose $T$ such that $\omega \equiv dT$
mod $(T)$. Using the cubic (\ref{cubic}), we may take, for
instance, $T=-\frac{x}{y}$. We set $W:=-\frac{1}{y}$.

This affine coordinate $T$, and $W$ around the zero section, are mapped into
$T \in A[[T]]$ and
\begin{equation}
\label{dub} W(T)=T^3+a_4T^7+\cdots \in A[[T]]\, ,
\end{equation}
respectively; cf. \cite{sil2}, p. 111.

For all integers $a,b \in \bZ_+$ with $a+b \geq 1$, we define
{\it $\D$-Kodaira-Spencer classes}
$f^{[a,b]} \in A$ as follows. For $r \geq a+b$, we have a natural isomorphism
$N^r \simeq (\hat{\bA}^{\frac{r(r+3)}{2}},[+])$, where
$\hat{\bA}^{\frac{r(r+3)}{2}}=Spf\ A[T^{(i,j)}\mid_{1 \leq i+j
\leq r}]\wh$, which we view hereafter as an identification.
Thus, we obtain the identification
$$\begin{array}{c}
s_{\mu \nu}=(\alpha_{\mu \nu}^{i,j}\mid_{1 \leq i+j \leq r}) \in
\cO(\hU_{\mu \nu})^{\frac{r(r+3)}{2}}\, , \\
\alpha_{\mu \nu}^{i,j}:=T^{(i,j)} \circ s_{\mu \nu}\, .
\end{array}
$$
By (\ref{romanu}), we have the
series $L^{[a,b]} \in A[T^{(i,j)}\mid_{1 \leq i+j \leq r}]\wh$
defining homomorphisms
$$L^{[a,b]}:N^r \simeq (\hat{\bA}^{\frac{r(r+3)}{2}},[+]) \ra \hat{{\mathbb G}
}_a=(\hat{\bA}^1,+)\, .$$
We define elements
$$\varphi^{[a,b]}_{\mu \nu}:=L^{[a,b]}(\alpha_{\mu \nu}^{i,j}\mid_{1 \leq i+j
 \leq r}) \in \cO(\hU_{\mu \nu})\, .$$
As we vary $\mu \nu$, the collection of such sections defines
a cohomology class
$$\varphi^{[a,b]} \in H^1(\hat{E},\cO_{\hat{E}})=H^1(E,\cO_E) \, ,$$
which is, of course, the class $\gamma(L^{[a,b]})$
defined by the exact sequence (\ref{papica}).
Let
$$\<\ ,\ \>:H^1(E,\cO) \times H^0(E,\Omega^1) \ra A $$
denote the Serre duality pairing, and define
\begin{equation}
\label{def of fr}
f^{[a,b]}:=\<\varphi^{[a,b]},\omega\> \in A\, .
\end{equation}
It is clear that $f^{[a,b]}$ depends only on $a$ and $b$ but not
on $r$. It is also clear that $\partial$ in the exact sequence
(\ref{papica}) is $0$ if and only if $f^{[a,b]}=0$ for all $1 \leq
a+b \leq r$. Proceeding verbatim as in \cite{book}, we check that
$f^{[a,b]}$ depends on the pair $(E,\omega)$, and not on the
choice of $T$, as long as $T$ satisfies the condition that $\omega
\equiv dT$ mod $(T)$.

Now let $a,b,c,d \in \bZ_+$ with $1 \leq a+b,c+d \leq r$, and
and consider the homomorphism
$$\Theta:=\Theta^{[a,b]}_{[c,d]}
:N^r \ra \hat{\bG}_a $$
 given by
$$\Theta:=f^{[a,b]} L^{[c,d]}-f^{[c,d]}L^{[a,b]}\in A[T^{(i,j)}\mid_{1
\leq i,j \leq r}]\wh\, .$$
Then we have that $\partial(\Theta)=0$, for $\partial(\Theta)$ is the class
$[\gamma_{\mu \nu}]$ in $H^1(\hat{E},\cO)$ of the cocycle
$$\gamma_{\mu \nu}:=
f^{[a,b]}\varphi_{\mu \nu}^{[c,d]}-f^{[c,d]}\varphi_{\mu \nu}^{[a,b]}
 \in \cO(\hU_{\mu \nu})\, ,$$
and $[\gamma_{\mu \nu}]=0$ because
$$\<[\gamma_{\mu \nu}],
\omega\>=f^{[a,b]} \cdot \<\varphi^{[c,d]},\omega\> -
 f^{[c,d]}
 \cdot
\<\varphi^{[a,b]},\omega\>=0\, .$$
By the exactness of the sequence (\ref{papica}), $\Theta$ lifts to a
unique homomorphism
$$\psi=\psi^{[a,b]}_{[c,d]}:J^r_{\di \ddi}(E) \ra \hat{\bG}_a $$
which we interpret as a $\D$-character
$$\psi \in \bX^r_{\di \ddi}(E)\, .$$
This character $\psi$ depends only on $a,b,c,d$, but not on $r$
in the sense that if we change $r$ to $r+s$, then the new $\psi$
is obtained from the old one by composition with the
projection $J^{r+s}_{\di \ddi}(E) \ra J^r_{\di \ddi}(E)$. Incidentally,
$\psi$ is obtained by gluing functions
\begin{equation}
\psi_{\mu} \circ \tau_{\mu}^{-1} \in \cO(\pi_r^{-1}(\hat{U}_{\mu}))\, ,
\end{equation}
\begin{equation}
\label{prea mult de munca}
\psi_{\mu}:=f^{[a,b]} \cdot L^{[c,d]} - f^{[c,d]}
 \cdot L^{[a,b]} +\gamma_{\mu} \in
\cO(\hat{U}_{\mu})[T^{(i,j)}\mid_{1 \leq i+j \leq r}]\wh\, ,
\end{equation}
with $\gamma_{\mu} \in \cO(\hU_{\mu})$, and we have
the identities
$$\begin{array}{rcl}
\psi^{[a,b]}_{[a,b]} & = & 0 \, , \vspace{1mm}\\
\psi^{[a,b]}_{[c,d]}+\psi^{[c,d]}_{[a,b]} & = & 0\, , \vspace{1mm} \\
f^{[a_1,b_1]} \psi^{[a_2,b_2]}_{[a_3,b_3]}+
f^{[a_2,b_2]} \psi^{[a_3,b_3]}_{[a_1,b_3]}+
f^{[a_3,b_3]} \psi^{[a_1,b_1]}_{[a_2,b_2]} & = & 0\, ,
\end{array}$$
which follow from the very same identities that are obtained when the $\psi$s
are replaced by  the $\Theta$s.

Since $\tau_{\mu_0}$ is the identity modulo $T$, we have that
$\psi \circ e(pT)$ is congruent to $\psi_{\mu_0} \circ e(pT)$ modulo $T$.
By Lemmas \ref{caiin} and \ref{xxzz}, we obtain
\begin{equation}
\label{supitza} \left( \psi^{[a,b]}_{[c,d]} \right) \circ e(pT)=
p^{1+\epsilon(d)} f^{[a,b]} \phip^c \dd^d T - p^{1+\epsilon(b)}
f^{[c,d]} \phip^a \dd^b T +\tf^{[a,b]}_{[c,d]} T\, ,
\end{equation}
where $\tf^{[a,b]}_{[c,d]} \in A$. So $\psi^{[a,b]}_{[c,d]}$,
viewed as an element of $A[[T]][\d^i \dd^j T\mid_{1 \leq i+j \leq r}]\wh$,
has the form
\begin{equation}
\label{incaosupitza} \psi^{[a,b]}_{[c,d]}=\frac{1}{p}
\left[p^{1+\epsilon(d)} f^{[a,b]} \phip^c \dd^d  -
p^{1+\epsilon(b)} f^{[c,d]} \phip^a \dd^b  +\tf^{[a,b]}_{[c,d]}
\right] l(T)\, .
\end{equation}

\begin{remark}
We will use the following notation:
$$\begin{array}{rclrclrcl}
f^a_{\di} & := & f^{[a,0]}\, , & f^a_{\ddi} & := & f^{[0,a]}\, , &  & \\
\tf^2_{\di} & := & \tf^{[1,0]}_{[2,0]}\, , & \tf^2_{\ddi} & := &
\tf^{[0,1]}_{[0,2]}\, , & \tf^1_{\di \ddi} & := & \tf^{[1,0]}_{[0,1]}\, .
\end{array}
$$
The elements $f^a_{\di} \in A$ are (the images of) the elements $f_{jet}^a$
in \cite{book}.  The elements $\tf_{\di}^2$ are (the images of) the
elements $pf^{1,2}_{jet}$ in \cite{book}. By \cite{book},
Proposition 7.20 (and Remark 7.21) and Corollary 8.84 (and Remark 8.85), we
have
\begin{equation}
\label{4res} \tf_{\di}^2=p(f^1_{\di})^{\phi}\, .
\end{equation}

The element $f^1_{\di}$ was interpreted in \cite{book} as an {\it
arithmetic Kodaira-Spencer class} of $E$. The element $f^1_{\ddi}$
is easily seen to be  (an incarnation of) the usual
Kodaira-Spencer class of $E$; cf. \cite{difmod}. The element
$f^1_{\ddi} \in A$, and the reduction mod $p$ of $f^1_{\di}$, were
explicitly computed in \cite{hurl}. The $\dd$-character
$$\psi^2_{\ddi}:=\psi^{[0,1]}_{[0,2]}\in \bX_{\ddi}^2(E)$$
is (an ``incarnation'' of) the {\it Manin map} of $E$ \cite{man},
constructed as in \cite{annals}. If $f^1_{\ddi} \neq 0$, then
$\psi^2_{\ddi} \neq 0$. (Unlike the construction in
\cite{annals} that was done over a field, our construction here
is carried over the ring $A$.) The $\d$-character
$$\psi^2_{\di}:=\psi^{[1,0]}_{[2,0]}\in \bX_{\di}^2(E)$$
is the {\it arithmetic Manin map} in \cite{char}. And if $f^1_{\di}
\neq 0$, then $\psi^2_{\di} \neq 0$. Both of these Manin maps are
``ode'' maps with respect to the geometric and the arithmetic
direction separately. On the other hand, we can consider  the
$\D$-character
$$ \psi^1_{\di \ddi}:=\psi^{[1,0]}_{[0,1]}\in \bX_{\di \ddi}^1(E)\, .$$
If either  $f^1_{\di} \neq 0$ or $f^1_{\ddi} \neq 0$, then
$\psi_{\di \ddi}^1 \neq 0$. If both $f^1_{\di} \neq 0$ and
$f^1_{\ddi} \neq 0$, then $\psi^1_{\di \ddi}$
 is a ``purely pde'' operator (in the sense that it is not a sum of a
$\d$-character and a $\dd$-character).
\end{remark}

Indeed, we have the
following consequence of Proposition \ref{indi}.

\begin{corollary}
If $f^1_{\di} \neq 0$ and  $f^1_{\ddi} \neq 0$, then:
\begin{enumerate}
\item $\psi^1_{\di \ddi}$ is an $L$-basis of  $\bX^1_{\di \ddi}(E)\otimes
L$.
\item $\bX_{\di}^1(E)=\bX_{\ddi}^1(E)=0$.
\end{enumerate}
\end{corollary}

We may view $\psi^1_{\di \ddi}$ as a canonical ``convection
equation'' on $E$.

Using the notation above and the commutation relations in $A[\phip,\dd]$,
 we have the following equalities:
\begin{equation}
\label{nh}
\begin{array}{rcl}
\psi^1_{\di \ddi} & = & \frac{1}{p} \left[ pf^1_{\di} \dd-f^1_{\ddi}
\phip +\tf^1_{\di \ddi} \right] l(T)\, , \vspace{1mm} \\
\dd \psi^1_{\di \ddi} & = & \frac{1}{p} \left[ p(\dd f^1_{\di})
\dd+pf^1_{\di}\dd^2 - (\dd f^1_{\ddi}) \phip -f^1_{\ddi} p \phip
\dd +
\dd \tf_{\di \ddi}^1 + \tf_{\di \ddi}^1 \dd \right] l(T)\, ,\vspace{1mm} \\
\phip \psi^1_{\di \ddi} & = & \frac{1}{p} \left[ p
(f^1_{\di})^{\phi} \phip \dd-(f^1_{\ddi})^{\phi}\phip^2+ (\tf_{\di
\ddi}^1)^{\phi} \phip \right] l(T)\, , \vspace{1mm} \\
\psi^2_{\di} & = & \frac{1}{p} \left[f^1_{\di} \phip^2-f^2_{\di}
\phip
+\tf^2_{\di} \right] l(T)\, ,\vspace{1mm} \\
\psi_{\ddi}^2 & = & \frac{1}{p} \left[ pf^1_{\ddi} \dd^2-pf^2_{\ddi}
\dd+\tf_{\ddi}^2 \right] l(T)\, .
\end{array}
\end{equation}
Thus, if we represent the ordered elements
$\dd \psi^1_{\di \ddi},\ \phip \psi^1_{\di \ddi}, \ \psi^2_{\di},
\ \psi^1_{\di \ddi},\  \psi^2_{\ddi}$ as $L$-linear
combinations of the series
$p \dd^2 l(T),\ p \phip \dd l(T),\ \phip^2 l(T),\ p \dd l(T), \
\phip l(T),\ l(T)$, then the matrix of $L$-coefficients is equal to
$\frac{1}{p}M$, where
$$M:=\left(
\begin{array}{cccccc}
f^1_{\di} & -f^1_{\ddi} & 0 & \dd f^1_{\di}+\frac{\tf^1_{\di
\ddi}}{p} & - \dd f^1_{\ddi} & \dd \tf^1_{\di \ddi} \vspace{2mm} \\
0 & (f^1_{\di})^{\phi} &
-(f^1_{\ddi})^{\phi} & 0 & (\tf^1_{\di
\ddi})^{\phi} & 0 \vspace{2mm} \\
0 & 0 & f^1_{\di} & 0 & -f_{\di}^2 & \tf_{\di}^2 \vspace{2mm} \\
0 & 0 & 0 & f^1_{\di} & -f^1_{\ddi} & \tf^1_{\di \ddi}\vspace{2mm}  \\
f^1_{\ddi} & 0 & 0 & -f^2_{\ddi} & 0 & \tf_{\ddi}^2
\end{array}
 \right)\, .
$$

\begin{proposition}
\label{narnia} Let us assume that $f^1_{\di} \neq 0$ and $f^1_{\ddi}\neq 0$.
Then the following hold:
\begin{enumerate}
\item[1)] The elements
$$\psi_{\di \ddi}^1,\  \dd \psi^1_{\di \ddi},\
\phip \psi^1_{\di \ddi},\  \psi^2_{\di}
$$
form an $L$-basis of $\bX_{\di \ddi}^2(E) \otimes L$.
\item[2)] The elements
$$\psi_{\di \ddi}^1,\  \phip \psi^1_{\di \ddi},\
\psi^2_{\ddi},\ \psi^2_{\di}
$$
form an  $L$-basis of $\bX_{\di \ddi}^2(E) \otimes L$.
\item[3)] There exists a $5$-tuple $(\alpha,\beta,\gamma,\nu,\lambda) \in
A^5$, which is unique up to scaling by an element of $A$, satisfying
$$(\alpha \dd+\beta \phip +\gamma)\psi^1_{\di \ddi}=\nu
\psi^2_{\ddi}+\lambda \psi_{\di}^2,\quad \alpha \neq 0, \; \nu
\neq 0\, .
$$
\item[4)] All $5 \times 5$ minors of the matrix $M$ vanish.
\end{enumerate}
\end{proposition}

We may view the character $\nu \psi^2_{\ddi}+\lambda \psi_{\di}^2$ as a
canonical ``wave equation'' on $E$.

{\it Proof}.
By the form of the matrix $M$, each of the $4$-tuples in assertions
1) and 2) are $L$-linearly independent. By Proposition
\ref{indi}, $\bX_{\di \ddi}^2(E)$ has rank $4$ over $A$. All the
statements in the Proposition then follow. \qed

Similar arguments (cf. also to the case of $\bG_m$) yield:

\begin{proposition}
Let us assume that $f^1_{\di} \neq 0$ and $f^1_{\ddi}\neq 0$. Then, for any $r
\geq 2$, the $L$-vector space $\bX_{pq}^r(E) \otimes L$ has
basis
$$\{\phip^i \dd^j \psi^1_{pq}\mid_{0 \leq i+j \leq r-1}\}
\cup \{\phip^i \psi^2_p\mid_{0 \leq i \leq r-2}\}\, .
$$
\end{proposition}

When either $f^1_{\di}=0$ or $f^1_{\ddi}=0$, the picture above changes, and
in fact, it simplifies. For if with more generality we assume that
$f^{[a,b]}=0$ for some $a,b$ with $a+b \leq r$, then
$0=\varphi^{[a,b]}=\gamma(L^{[a,b]})$, and by the exact sequence
(\ref{papica}), $L^{[a,b]}$ lifts uniquely to a homomorphism
$\psi^{[a,b]}:J^r_{\di \ddi}(E) \ra \hat{\bG}_a$, which we interpret as a
$\D$-character $\psi^{[a,b]} \in \bX_{\di \ddi}^r(E)$.
Let us observe in passing that $\psi^{[a,b]}$ is obtained by gluing functions
$$\begin{array}{c}
\psi_{\mu} \circ \tau_{\mu}^{-1} \in
\cO(\pi_r^{-1}(\hat{U}_{\mu}))\, , \\
 \psi_{\mu}:=L^{[a,b]}  +\gamma_{\mu} \in
\cO(\hat{U}_{\mu})[T^{(i,j)}\mid_{1 \leq i+j \leq r}]\wh \, ,
\end{array}
$$
with $\gamma_{\mu} \in \cO(\hU_{\mu})$. As before, we obtain
$$\psi^{[a,b]} \circ e(pT)=p^{1+\epsilon(b)} \phip^a \dd^b
T+\tf^{[a,b]}\, $$
where $\tf^{[a,b]} \in A$. So $\psi^{[a,b]}$, viewed as a series, has the form
$$\psi^{[a,b]}=\frac{1}{p}\left[ p^{1+\epsilon(b)} \phip^a \dd^b
+\tf^{[a,b]} \right] l(T)\, .$$

Incidentally, if $f^1_{\ddi}=0$, then
\begin{equation}
\label{psi1q} \psi^1_{\ddi}:=\psi^{[0,1]} \in \bX^1_{\ddi}(E)
\end{equation}
is (an ``incarnation" of) the {\it Kolchin logarithmic derivative}
of $E$ \cite{kolchin}; cf. Proposition \ref{logder} below.

On the other hand, if $f^1_{\di}=0$ we set
\begin{equation}
\label{psi1p} \psi^1_{\di}:=\psi^{[1,0]} \in \bX_{\di}(E)\, .
\end{equation}

As before, we find the following bases for the set of
$\D$-characters:

\begin{proposition}
Let us assume that $f^1_{\di} = 0$ and $f^1_{\ddi}\neq 0$. Then, for each $r
\geq 2$, the $L$-vector space $\bX_{pq}^r(E) \otimes L$ has basis
$$\{\phip^i \dd^j \psi^2_q\mid_{\ 0 \leq i+j \leq r-2}\} \cup
 \{\phip^i \psi^1_p\mid_{0 \leq i \leq r-1}\}
 \cup \{\phip^i \dd \psi^1_p\mid_{0 \leq i \leq r-2}\}\, .$$
\end{proposition}

\begin{proposition}
\label{ileana1} Let us assume that $f^1_{\di} \neq 0$ and $f^1_{\ddi}=0$.
Then, for each $r \geq 2$, the $L$-vector space $\bX^r_{\di \ddi}(E)
\otimes L$ has basis
$$\{\phip^i \psi_p^2\mid_{0 \leq i \leq r-2}\} \cup \{\phip^i \dd^j
\psi^1_{\ddi}\mid_{0 \leq i+j \leq r-1}\}\, .
$$
\end{proposition}

\begin{proposition}
\label{ileana2}
Let us assume $f^1_{\di} = 0$ and $f^1_{\ddi}= 0$.  Then,
for each $r \geq 1$, the $L$-vector space $\bX^r_{\di \ddi}(E)
\otimes L$ has basis
$$\{\phip^i \psi_p^1\mid_{0 \leq i \leq r-1}\} \cup \{\phip^i \dd^j
\psi^1_{\ddi}\mid_{0 \leq i+j \leq r-1}\}\, .
$$
\end{proposition}

\subsection{Tate curves}
Let $E_q$ be  the {\it Tate curve with parameter} $q$ over
$A:=R((q))\wh$, equipped with its canonical $1$-form $\omega_q$.
This curve $E_q$ is defined as the elliptic curve in the
projective plane over $A$ whose affine plane equation is
$$y^2=x^3-\frac{1}{48}E_4(q)x+\frac{1}{864}E_6(q)\, ,$$
where $E_4$ and $E_6$ are the Eisenstein series
$$\begin{array}{rcl}
E_4(q) & = & 1+240 \cdot s_3(q)\, ,\\
E_6(q) & = & 1-504 \cdot s_5(q)\, .
\end{array}
$$
In here, for $m \geq 1$, we follow the usual convention and write
$$s_m(q):=\sum_{n \geq 1} \frac{n^mq^n}{1-q^n} \in R[[q]]\, .$$
Also,  the canonical form is defined by
$$\omega_q=\frac{dx}{y}\, .$$

More generally, let $\beta \in R^{\times}$ be an invertible
element that we shall view as a varying parameter, and let
$(E_{\beta q},\omega_{\beta q})$ be the pair obtained by base
change from $(E_q,\omega_q)$ via the isomorphism
$$\begin{array}{rcl}
R((q))\h & \stackrel{\sigma_{\beta}}{\ra} & R((q))\wh \\
\sigma_{\beta}(\sum a_nq^n)
& =  & \sum a_n \beta^n q^n
\end{array}\, .
$$
Thus, the pair $(E_{\beta q},\omega_{\beta
q})$ is the elliptic curve
$$y^2=x^3-\frac{1}{48}E_4(\beta q)x+\frac{1}{864}E_6(\beta q)$$
equipped with the $1$-form
$$\omega_{\beta q}=\frac{dx}{y}\, .$$
(Let us observe that $\sigma_{\beta}$ is a $\dd$-ring
homomorphism, but not a $\d$-ring homomorphism.)
We shall refer to $(E_{\beta q},\omega_{\beta q})$ as the {\it
Tate curve with parameter $\beta q$}.

Let us observe that the discussion and notation in Example \ref{cucurigu} that
concerns groups of solutions does not apply to $E_{\beta q}$ over $A$, as
$A=R((q))\wh \neq R$.

From now on, we assume that all our quantities in the general
theory (like the $f$'s, the $\tilde{f}$'s,  and the $\psi$'s) are
associated to the pair $(E_{\beta q},\omega_{\beta q})$.

\begin{lemma}
\label{tate1}  There exists $c \in \bZ_p^{\times}$ such that, for
any $\beta \in R^{\times}$, we have
$$f^1_{\di}=\eta,\quad f^2_{\di}=\eta^{\phi}+p\eta,\quad \tf_{\di}^2=p
\eta^{\phi}\, ,$$
where
$$\eta:=\frac{c}{p}\log \frac{\phi(\beta)}{\beta^p}\, .$$
\end{lemma}

{\it Proof}. We consider the ring $R((q))\h[q',q'',\ldots]\wh$ with its unique
$\d$-ring structure such that $\d q=q'$, $\d q'=q''$, etc. Let us note
that $(E_{\beta q},\omega_{\beta q})$ is obtained by base change
from $(E_q,\omega_q)$ via the composition
\begin{equation}
\label{tzutz} R((q))\h \subset R((q))\h[q',q'',\ldots]\h
\stackrel{\tilde{\sigma}_{\beta}}{\ra} R((q))\h[q',q'',\ldots]\h
\stackrel{\pi}{\ra}R((q))\wh \, ,
\end{equation}
 where
$\tilde{\sigma}_{\beta}(q)=\beta q$,
$\tilde{\sigma}_{\beta}(q')=\d(\beta q)$,
 $\tilde{\sigma}_{\beta}(q'')=\d^2(\beta q),\ldots $, and
$\pi(q)=q$, $\pi(q')=0$, $\pi(q'')=0, \ldots$  And let us note also that
$\tilde{\sigma}_{\beta}$ and $\pi$ are $\d$-ring
homomorphisms. By the functoriality of $f^r_{\di}$ plus
\cite{difmod}, Corollary 7.26, \cite{barcau}, Corollary 6.1, and
\cite{fermat}, Lemma 6.14, there exists $c \in \bZ_p^{\times}$
such that, for any $\beta$, we have
$$\begin{array}{rcl}
f^1_{\di} & = & {\displaystyle \frac{c}{p} \cdot\log\frac{\phip(\beta q)}
{(\beta q)^p}=\frac{c}{p}\log \frac{\phip(\beta)}{\beta^p}}\, , \vspace{1mm} \\
f^2_{\di} & = & (f^1_{\di})^{\phi}+pf^1_{\di}\, .
\end{array}$$
Thus, $f^1_{\di}$ and $f^2_{\di}$ have the desired values. The value of
$\tf_{\di}^2$ follows by (\ref{4res}). \qed

\begin{lemma}
\label{tate2} The following hold:
\begin{enumerate}
\item[1)] $f^1_{\ddi}$ does not depend on $\beta$.
\item[2)]  $f^1_{\ddi} \in \bZ_p^{\times}$, and $f^2_{\ddi}=\tf^2_{\ddi}=0$.
\end{enumerate}
\end{lemma}

{\it Proof}. We start by noting that $(E_{\beta q},\omega_{\beta q})$ is
obtained by a base change from $(E_q,\omega_q)$ via the isomorphism
$\sigma_{\beta}:R((q))\h \ra R((q))\wh$, and recall that $\sigma_{\beta}$ is
a $\dd$-ring homomorphism. By functoriality, the elements $f^1_{\ddi}$,
$f^2_{\ddi}$, $\tilde{f}^2_{\ddi}$ corresponding to a given $\beta$ are the
images via $\sigma_{\beta}$ of the corresponding quantities for
$\beta=1$. Thus, proving assertion 2) for $\beta=1$ suffices to
conclude both assertions 1) and 2) for arbitrary $\beta$.

We prove next assertion 2) for $\beta=1$.

The fact
that $f^1_{\ddi} \in \bZ_p^{\times}$ follows from \cite{difmod},
Corollary 7.26. In order to prove that $f^2_{\ddi}=\tf^2_{\ddi}=0$,
it is enough to show that there exists $\psi \in \bX_{\ddi}^2(E)$
such that $\psi=\dd^2 l(T)$. This is because $\bX^2_{\ddi}(E)$ has
rank $1$ over $A$ \cite{ajm}.

The existence of $\psi$ can be argued analytically, and we briefly sketch the
argument next. (A purely algebraic argument is also available but the
analytic one is simpler and classical, going back to Fuchs and
Manin \cite{man}.) For in our problem we can replace
$R((q))$ by $\bQ((q))$, and then $\bQ((q))$ by $\bC((q))$. We view the
Tate curve over $\bC((q))$ as arising from an analytic family
$E_q^{an} \ra \Delta^*$ of elliptic curves over the punctured disk
$\Delta^*$. The fiber $E_{q_{\tau}}$ of this family over a point
$q_{\tau}=e^{2 \pi i \tau} \in \Delta^*$ identifies with
$\bC/\<1,\tau\>$. Let $z$ be the coordinate on $\bC$. Then, for
any local analytic section $P$, $q_{\tau} \mapsto P(q_{\tau})$, of
$E_q^{an} \ra \Delta^*$ that is close to the zero section
$q_{\tau} \mapsto P_0(q_{\tau})=0$, we have that
\begin{equation} \label{su} l(P(q_{\tau}))=
\int_{P_0(q_{\tau})}^{P(q_{\tau})} \omega_{q_{\tau}}
\end{equation}
where $\omega_{q_{\tau}}$ is the $1$-form on $E_{q_{\tau}}$ whose
pull back to $\bC$ is $dz$. The periods of $\omega_{q_{\tau}}$
on $E_{q_{\tau}}$ are $1$ and $\tau$, so they are annihilated by
$$\left( \frac{d}{d \tau} \right)^2=-4 \pi^2 \left( q_{\tau}
 \frac{d}{dq_{\tau}}\right)^2\, .$$
Hence the map
$$P  \mapsto \left( q_{\tau}
 \frac{d}{dq_{\tau}}\right)^2 \left(
 \int_{P_0(q_{\tau})}^{P(q_{\tau})} \omega_{q_{\tau}} \right)$$
is well defined for all local analytic sections $P$ of $E_q^{an}
\ra \Delta^*$ (not necessarily close to the zero section $P_0$).
This map arises from a $\dd$-character $\psi$ of $E_q$ over
$\bC((q))$ \cite{man}. By (\ref{su}), $\psi$ coincides with
$P \mapsto \dd^2 l(P)$ for $P$ close to the zero section, and this completes
the argument. \qed

\begin{remark}
\label{tate1.5} From now on, we choose  $c$ as in Lemma \ref{tate1}, and we set
$$\begin{array}{rcl}
\eta & := & f^1_{\di}=\frac{c}{p}\log \frac{\phip(\beta)}{\beta^p} \in R\,
, \vspace{1mm} \\
 \gamma & := & f^1_{\ddi} \in
\bZ^{\times}_p\, .
\end{array}
$$
Notice that we have $\eta=0$ if, and only if,
$\beta$ is a root of unity. Later on, the
quotient of these Kodaira-Spencer classes,
$$\frac{\eta}{\gamma}=\frac{f^1_{\ddi}}{f^1_{\di}}\, ,$$
will play a key r\^ole.
\end{remark}

\begin{lemma}
\label{tate3} $\tf^1_{\di \ddi}=p \gamma$.
\end{lemma}

{\it Proof}. Let us assume first that $\beta$ is not a root of unity, and so,
$\eta \neq 0$. The $5 \times 5$ minor of the matrix $M$ above obtained
by removing the $6$-th column has determinant zero (cf. Proposition
\ref{narnia}), and using Lemma \ref{tate1}, Lemma \ref{tate2}, and
Remark \ref{tate1.5}, we get the condition
$$pf^1_{\di}(\tf^1_{\di \ddi} -p\gamma)^{\phi}+(f^1_{\di})^{\phi}
(\tf^1_{\di \ddi}-p\gamma)=0\, .
$$
We take the valuation $v_p:R((q))\h \ra \bZ_+ \cup \{\infty\}$ in this
identity, and use the fact that $v_p \circ \phi=v_p$ to get
$\tf^1_{\di \ddi}-p\gamma=0$, completing the
proof of this case.

If instead now $\beta \in R^{\times}$ is arbitrary, proceeding as in the
proof of Lemma \ref{tate1}, we have that $\tilde{f}^1_{\di \ddi}$
is the image of an element
$F \in R((q))\h[q',q'',\ldots ]\wh$ via the map
$\pi \circ \tilde{\sigma}_{\beta}$ in (\ref{tzutz}).
We know that
\begin{equation}
\label{tzinee} \pi \tilde{\sigma}_{\beta} F=p \gamma
\end{equation}
 for all $\beta \in R^{\times} \backslash \bmu(R)$.
Since the map
$$\begin{array}{rcl}
R^{\times} & \ra  & R((q))\wh \\ \beta & \mapsto & \pi
\tilde{\sigma}_{\beta} F
\end{array}$$
is $p$-adically continuous, we have that that (\ref{tzinee}) holds for all
$\beta \in R^{\times}$. This completes the proof. \qed

Combining Lemma \ref{tate1}, Lemma \ref{tate2}, and (\ref{nh}), we derive
the following.

\begin{proposition}
\label{pomu}  Let $E=E_{\beta q}$ be the Tate curve over
$A=R((q))\h$ with parameter $\beta q$, $\beta \in R^{\times}$.
When viewed as elements of the ring $$A[[T]][\d T, \dd T, \d^2 T, \d \dd T,
\dd^2 T]\wh\, ,$$
the characters  $\psi^1_{\di \ddi}$, $\psi^2_{\ddi}$,
$\psi_{\di}^2$ are equal to
$$\begin{array}{rcl}
\psi^1_{\di \ddi} & = & {\displaystyle \frac{1}{p} [p \eta \dd -\gamma
\phip+p\gamma ]l(T)}\, ,
\vspace{2mm} \\ \psi^2_{\ddi} & = & \gamma \dd^2 l(T)\, , \vspace{2mm} \\
\psi^2_{\di} & = & {\displaystyle \frac{1}{p} [\eta
\phip^2-(\eta^{\phi}+p\eta)\phip+p\eta^{\phi}]l(T)}\, .
\end{array}
$$
In particular, we have the relations
$$\begin{array}{rcl}
(\gamma \eta^{\phi+1}\dd+ \gamma^2 \eta \phip-\gamma^2
\eta^{\phi}) \psi^1_{\di \ddi} & = & \eta^{\phi+2}
\psi^2_{\ddi}- \gamma^3 \psi_{\di}^2 \, , \vspace{2mm} \\
\gamma \dd^2 \psi^1_{\di \ddi} & = & (\eta \dd -p \gamma
\phip+\gamma) \psi^2_{\ddi}\, .
\end{array}$$
\end{proposition}

\begin{remark}
Of course, if $\beta$ is a root of unity, the first of
the last two relations in the Proposition above reduces to the identity $0=0$.
In this case, we also have that
\begin{equation}
\label{fome} \psi^1_{\di \ddi}=-\gamma \psi^1_{\di}\, .
\end{equation}
\end{remark}

\begin{remark}
\label{spallate} By the Proposition above, the Fr\'{e}chet symbol
of $\psi^1_{\di \ddi}$ with respect to $\omega$ is
$\theta_{\psi^1_{\di \ddi}, \omega}=\eta \xi_q-\gamma
\xi_p+\gamma$. Thus, if we consider the ``simplest'' of the energy functions
on $E$ given by
$$H:=(\psi^1_{\di \ddi})^2\, ,$$
a direct computation shows the following expression for the
Euler-Lagrange equation attached to $H$ and the vector field
$\partial$ dual to $\omega$:
$$\epsilon^1_{H,\partial}=(-2 \eta^{\phi} \phip \dd+2 \gamma
\phip-2 \gamma)\psi^1_{\di \ddi}\, .
$$
In particular, any solution of $\psi^1_{\di \ddi}$ is a solution of the
Euler Lagrange equation $\epsilon^1_{H,\partial}$.
\end{remark}

In the sequel, we will need to use some other facts about Tate's
curves that we now recall. Indeed, note that the cubic defining $E_{\beta q}$
has coefficients in $R[[q]]$, and hence, $E_{\beta q}$ has a natural
projective (non-smooth) model ${\mathcal E}={\mathcal E}_{\beta
q}$ over $R[[q]]$ equipped with a ``zero'' section defined by
$T=-\frac{x}{y}$. The completion of ${\mathcal E}_{\beta q}$ along
this section is naturally isomorphic to
$Spf\ R[[q]][[T]]$, where $R[[q,T]]=R[[q]][[T]]$ is viewed with its $T$-adic
topology.
This isomorphism induces a natural embedding
$$\iota:qR[[q]] \ra {\mathcal E}(R[[q]]) \subset E(A)\, ,$$
that is explicitly given by sending any $u \in qR[[q]]$
into the point of ${\mathcal E}$ whose $(T,W)$-coordinates are
$$(u,u^3-\frac{1}{48}E_4(\beta q)u^7+\cdots) \, ,$$
cf. (\ref{dub}). We denote by $E_{1}(A)$ the image of this embedding, that is
to say, $E_{1}(A)=\iota(qR[[q]])$. When no confusion can arise, we will view
$\iota$ as an inclusion, and we will identify $\iota(qR[[q]])$ with $qR[[q]]$.

 The formal group law $\cF=\cF_E=\cF(T_1,T_2)$ of $E$ with
respect to $T=-\frac{x}{y}$ has coefficients in $R[[q]]$, and is
isomorphic over $R[[q]]$ to the formal group
$\cF_{\bG_m}=t_1+t_2+t_1t_2$ of $\bG_m$ via an isomorphism
$\sigma(t)=t+\cdots\in t+t^2R[[q,t]])$.

For convenience, we recall how $\sigma$ arises. We first notice that
for any variable $v$, the series
\begin{equation}
\label{uglys}
\begin{array}{rcl} X(q,v) & := &
\frac{v}{(1-v)^2}+ \sum_{n \geq 1} \left( \frac{\beta^n
    q^nv}{(1-\beta^n q^n v)^2}+
\frac{\beta^n
    q^nv^{-1}}{(1-\beta^n q^n v^{-1})^2}-2 \frac{\beta^n
    q^n}{(1-\beta^n q^n)^2}
    \right)\, , \vspace{1mm} \\
Y(q,v) & := & \frac{v^2}{(1-v)^3}+ \sum_{n \geq 1} \left(
 \frac{(\beta^n
q^nv)^2}{(1-\beta^n q^n v)^3} - \frac{\beta^n
q^nv^{-1}}{(1-\beta^n q^n v^{-1})^3}+ \frac{\beta^n
q^n}{(1-\beta^n q^n)^2} \right)\, ,
\end{array}
\end{equation}
make sense as elements of the ring $\bZ[v,v^{-1}][[q]][\frac{1}{1-v}]$.
Cf. \cite{sil2}, p. 425. So, if we specialize $v \mapsto 1+t \in
\bZ[[t]]$ where $t$ is a variable, then $X$ and $Y$ become
elements in $\bZ[[q]]((t))$, that is to say, Laurent series in $t$ with
coefficients in $\bZ[[q]]$. Since $p \geq 5$, the series
\[\begin{array}{rcl}
x & = & X + \frac{1}{12} \vspace{1mm} \\
y & = & -Y-\frac{X}{2}
\end{array}\]
belong to $\bZ_{(p)}[[q]]((t))$. Note that $X=t^{-2}+\cdots$ and
$Y=t^3+\cdots$, hence, $x=t^{-2}+\cdots$ $y=-t^3+\cdots$ and so
$T=-\frac{x}{y}=t+\cdots=:\sigma(t) \in t+t^2\bZ_{(p)}[[q,t]]$.
It turns out that $t \mapsto T=\sigma(t)$ defines an isomorphism
$\cF_{\bG_m} \simeq \cF_E$; cf. \cite{sil2}, p. 431.

\begin{definition}
 The  {\it characteristic polynomial} $\mu(\xi_p,\xi_q)$ of a character
$\psi^1_{pq}$ is the Fr\'{e}chet symbol of $\psi^1_{\di \ddi}$ with respect
to $\omega$:
$$
\mu(\xi_p,\xi_q):=\eta \xi_q - \gamma  \xi_p
+\gamma\, .
$$
Note that $\mu(0,0) \in R^{\times}$. The {\it
characteristic  integers} of $\psi^1_{pq}$ are the  integers
$\kappa$ such that $\mu(0,\kappa)=0$, that is to say, solutions of
$\eta \kappa =-\gamma$. The set $\cK$ of characteristic integers has
 at most one element and is given by
$$\cK=\{-\gamma/\eta\} \cap \bZ\, .$$
A {\it totally non-characteristic} integer $\kappa \in \bZ$ is an
integer such that $\kappa \not\equiv 0$ mod $p$, and $\mu(0,\kappa)
\not\equiv 0$ mod $p$, that is,
$$\eta \kappa^2+\gamma \kappa \not\equiv 0\quad \text{mod}\; p\, .$$
We denote by $\cK'$ the set of totally non-characteristic
integers.
 For any  $0 \neq \kappa \in \bZ$ and
  any $\alpha \in R$, the {\it basic  series} is
\begin{equation}
\label{filmm} u_{E,\kappa,\alpha}:=e_E \left( \int
u_{a,\kappa,\alpha}^{\mu} \frac{dq}{q} \right) \in qK[[q]]\, ,
\end{equation}
where $e_E(T) \in TK[[q]][[T]]$ is the exponential of the formal
group law $\cF_E$, and
$u_{a,\kappa,\alpha}^{\mu}$ is as in  Equation
\ref{ua}. If in addition $\kappa$ is characteristic (i.e.
$\frac{\eta}{\gamma}=-\frac{1}{\kappa}$) then
\begin{equation}
\label{nuilacasino}
u_{a,\kappa,\alpha}^{\mu}:=\sum_{j \geq 0} (-1)^j
\frac{\alpha^{\phi^j}}{F_j(p)} q^{\kappa p^j}\, .
\end{equation}
Cf.  Example \ref{rupdi}.
\end{definition}

\begin{lemma}
\label{tate4} Let us assume that $\kappa \in \bZ \backslash p\bZ$.
Then we
have $u_{E,\kappa,\alpha} \in qR[[q]]$. Moreover, the map $R \ra
R[[q]]$, $\alpha \mapsto u_{E,\kappa,\alpha}$ is a pseudo
$\d$-polynomial map.
\end{lemma}

{\it Proof}. In order to prove the first assertion, let
$$h:=l_E(u_{E,\kappa,\alpha})=\int u_{a,\kappa,\alpha}^{\mu}
\frac{dq}{q}\, ,
$$
where $l_E(T)=l(T)\in K[[q]]$ is the logarithm of the formal group law
$\cF_E$. The isomorphism $\sigma(T)=T+\cdots \in
R[[q]][[T]]$ between $\cF_E$ and $\cF_{\bG_m}$ clearly has the
property that
$$e_E(T)=\sigma(e_{\bG_m}(T))\, ,$$
where
$$e_{\bG_m}(T)=\exp{(T)}-1=T+\frac{T^2}{2!}+\frac{T^3}{3!}+\cdots $$
So in order to check that $e_E(h) \in qR[[q]]$, it is enough to
show that $e_{\bG_m}(h) \in qR[[q]]$, that is to say, show that $\exp{(h)} \in
1+qR[[q]]$. By Dwork's Lemma \ref{dwor}, it is enough to show that
$(\phi-p)h \in pqR[[q]]$. Now, by Lemma \ref{unidul} we have that
$$\begin{array}
{rcl} p(\eta \kappa+\gamma)\alpha q^{\kappa} & = & p(\eta
\dd-\gamma \phip+\gamma)u_{a,\kappa,\alpha}^{\mu}\\
 & = & p(\eta \dd-\gamma \phip+\gamma)\dd h\\
 & = & \dd(p \eta \dd -\gamma \phip +p \gamma)h \, .
\end{array}$$
Hence
\begin{equation}
\label{tarre} (p\eta \dd-\gamma \phip+p\gamma)h-p(\eta
\kappa+\gamma)\alpha \kappa^{-1}q^{\kappa}
 \in R \cap qK[[q]]=0\, ,
\end{equation}
which gives
$$\begin{array}{rcl}
\gamma(\phip-p)h & = & p \eta \dd h-p(\eta \kappa+\gamma)\alpha
\kappa^{-1} q^{\kappa}\\
  & = & p \eta u_{a,\kappa,\alpha}^{\mu}-p(\eta \kappa +\gamma)
\alpha \kappa^{-1}q^{\kappa}\\
  & \in & pqR[[q]]\, ,
\end{array}
$$
and we are done.

The second assertion is proved exactly as in Lemma
\ref{ol}. \qed

We have the following ``diagonalization'' result.

\begin{lemma}
\label{oashteptpei}
\label{gainuli}
 Let us assume that $\kappa \in \bZ \backslash p\bZ$,
 and that $\alpha \in R$. Then
$$
\begin{array}{rcl}
\psi^2_q u_{E,\kappa,\alpha} & = & \gamma
\kappa u^{\mu^{(p)}}_{a,\kappa,\alpha}\\
\psi^1_{pq}u_{E,\kappa,\alpha}& = & (\eta \kappa+\gamma)\alpha
\kappa^{-1} q^{\kappa}\, .
\end{array}
$$
\end{lemma}

{\it Proof}. We have that
$$
\begin{array}{rcl}
\psi^2_q u_{E,\kappa,\alpha} & = &
\gamma \dd^2
l_E \left (e_E \left (
\int u_{a,\kappa,\alpha}^{\mu} \frac{dq}{q} \right) \right) \vspace{1mm} \\
\  & = & \gamma \dd^2 \left(\int u_{a,\kappa,\alpha}^{\mu}
\frac{dq}{q} \right) \vspace{1mm} \\
\  & = & \gamma \dd u^{\mu}_{a,\kappa,\alpha} \vspace{1mm} \\
\  & = & \gamma \kappa u^{\mu^{(p)}}_{a,\kappa,\alpha}.
\end{array}
$$
Consequently we have
$$\begin{array}{rcl}
\dd^2 \psi^1_{\di \ddi} u_{E,\kappa,\alpha} & = & \frac{1}{\gamma}
(\eta \dd -p \gamma\phip +\gamma) \psi^2_{\ddi}(u_{E,\kappa,\alpha})
\vspace{1mm} \\
  & = & (\eta \dd -p \gamma\phip +\gamma)
\kappa u_{a,\kappa,\alpha}^{\mu^{(p)}} \vspace{1mm} \\
  & = & \kappa \mu^{(p)}(\phip,\dd)
u_{a,\kappa,\alpha}^{\mu^{(p)}} \vspace{1mm} \\
  & = & \kappa \mu^{(p)}(0,\kappa) \alpha q^{\kappa} \vspace{1mm} \\
  & = & \kappa(\eta \kappa+\gamma) \alpha q^{\kappa} \vspace{1mm} \\
  & = & \dd^2((\eta \kappa+\gamma) \alpha \kappa^{-1} q^{\kappa})\, .
\end{array}
$$
Hence \[\dd(\psi^1_{pq}u_{E,\kappa,\alpha}-(\eta\kappa+\gamma)\alpha
\kappa^{-1}q^{\kappa}) \in R \cap qR[[q]]=0,\] and, therefore,
\[\psi^1_{pq}u_{E,\kappa,\alpha}-(\eta\kappa+\gamma)\alpha
\kappa^{-1}q^{\kappa} \in R \cap qR[[q]]=0,\] completing the proof.
\qed

\begin{remark}
\begin{enumerate}
\item[1)] For $\kappa \in \bZ \backslash p\bZ$ the map
$$\begin{array}{rcl}
R & \ra &  {\mathcal E}(R[[q]]) \\ \alpha & \mapsto &
\iota(u_{E,\kappa,\alpha})
\end{array}
$$
is a group homomorphism.
\item[2)] In analogy with (\ref{vocea}), if $\kappa \in
\bZ \backslash p\bZ$, we
may define {\it boundary value operators}
$$\begin{array}{rcl}
B_{\kappa}^0:{\mathcal E}(qR[[q]]) & \ra & R\, ,\\
 B_{\kappa}^0 u & = &
\Gamma_{\kappa} \psi^2_q u\, .
\end{array}
$$
Then, for $\kappa_1,
  \kappa_2 \in \bZ \backslash p\bZ$, we have
$$\begin{array}{rcl}
B_{\kappa_1}^0 u_{E,\kappa_2,\alpha} & = & \Gamma_{\kappa_1}
\psi^2_q u_{E,\kappa_2,\alpha} \vspace{1mm} \\
  & = & \Gamma_{\kappa_1} \gamma \kappa_2
u^{\mu^{(p)}}_{a,\kappa_2,\alpha} \vspace{1mm} \\
  & = & \gamma \alpha \kappa_2 \delta_{\kappa_1,
\kappa_2}\, .
\end{array}
$$
\item[3)] For $\kappa \in \bZ \backslash p\bZ$ we have
\begin{equation}
\label{pit} u_{E,\kappa,\zeta^{\kappa} \alpha}(q)
=u_{E,\kappa,\alpha}(\zeta q)
\end{equation}
for all $\zeta \in \bmu(R)$. In particular, if
$\alpha=\sum_{i=0}^{\infty} m_i \zeta_i^{\kappa}$, $\zeta_i \in
\bmu(R)$, $m_i \in \bZ$, $v_p(m_i) \ra \infty$, then
$$u_{E,\kappa,\alpha}(q)=\left[ \sum_{i=0}^{\infty}
\right] [m_i](u_{E,1,1}(\zeta_i^{\kappa} q^{\kappa}))= \left[
\sum_{i=0}^{\infty} \right] [m_i](u_{E,\kappa,1}(\zeta_i q))\,  .
$$
Here $[ \sum ]$ is the sum taken in the formal group law, and
$[m_i](T)\in R[[q]][[T]]$ is the series induced by
``multiplication by $m_i$'' in the formal group. It is known that
the sequence $[m_i](T)$ converges to $0$ in $R[[q]][[T]]$ in the
$(p,T)$-adic topology, so the series in the right hand side of the
equality above converges in the $(p,q)$-adic topology of
$R[[q]]$. Morally, $u_{E,\kappa,\alpha}$ is obtained from
$u_{E,\kappa,1}$ via ``convolution.'' (We have not defined
``convolution'' for groups not defined over $R$. This is possible but
we will not do it here.)
\item[4)] We have the following ``rationality'' property: if $\alpha,
\beta \in \bZ_{(p)}$ then $u_{E,\kappa,\alpha} \in
\bZ_{(p)}[[q]]$ for $\kappa \in \bZ \backslash p\bZ$.
\end{enumerate}
\end{remark}

For the next Proposition, we denote by
\begin{equation}
\label{lalalay} \cU^1_{\ra}:=\{u \in E(A)\ |\ \psi^1_{\di
\ddi}u=0\}
\end{equation}
the group of solutions of $\psi^1_{\di \ddi}$ in $A=R((q))\h$, and
we set \begin{equation} \label{poponeata} \cU^1_1:=\cU^1_{\ra}
\cap E_1(A).\end{equation}
 Note that there is no natural analogue
of $\cU_{\la}$ in this case.

 We let $E_{0}(A)$ be the subgroup of all $u \in E(A)$ that are solutions to
all $\dd$-characters of $E$. By \cite{man} or \cite{ajm},
$E_{0}(A)$ is actually the set of all solutions of the Manin map
$\psi^2_{\ddi}$. We set \begin{equation} \label{tedrav}
\cU^1_{0}:=\cU^1_{\ra} \cap E_{0}(A).\end{equation} We should view
$\cU^1_{0}$ as a substitute for the groups of stationary solutions
in the cases of $\bG_a$ and $\bG_m$, respectively.

\begin{theorem}
\label{unde} Let $E_{\beta q}$ be the Tate curve with parameter
$\beta q$, where $\beta \in R^{\times}$.  Let $\cU^1_{\ra}$ be the
group of solutions in $R((q))\h$ of the $\D$-character $\psi_{\di
\ddi}^1$. Then the following hold:
\begin{enumerate}
\item[1)] If $\frac{\eta}{\gamma} =-\frac{1}{\kappa}$ for some integer
$\kappa \geq 1$ coprime to $p$, and if $u_{E,\kappa,\alpha}$ is the
basic series in {\rm (\ref{filmm})},  then
$$\begin{array}{rcl}
\cU^1_{\ra} & = & \cU^1_{0}+\cU^1_1\, ,\\
 \cU^1_1 & = & \{
u_{E,\kappa,\alpha}\ |\ \alpha \in R\}\, .
\end{array}
$$
\item[2)] If $\frac{\eta}{\gamma}$ is not of the form
$\frac{\eta}{\gamma}=-\frac{1}{\kappa}$ for some integer $\kappa
\geq 1$ coprime to $p$, then
$$\cU^1_1=0\, .$$
\end{enumerate}
\end{theorem}

{\it Proof}. It is sufficient to prove the following claims:
\begin{enumerate}
\item[1)] If $\frac{\eta}{\gamma}=-\frac{1}{\kappa}$ for some
integer $\kappa \geq 1$ coprime to $p$, then $\psi^1_{\di \ddi}
u_{E,\kappa,\alpha}=0$ for all $\alpha$.
\item[2)] Assume  $\psi^1_{\di \ddi} u=0$ for some $0 \neq u
\in qR[[q]]$. Then  $\frac{\eta}{\gamma}=-\frac{1}{\kappa}$ for
some integer $\kappa \geq 1$ coprime to $p$, and
$u=u_{E,\kappa,\alpha}$ for some $\alpha$.
\item[3)] If $\frac{\eta}{\gamma}=-\frac{1}{\kappa}$ for some
integer $\kappa \geq 1$ coprime to $p$, and if $u \in E(A)$ is such
that $\psi^1_{\di \ddi}u=0$, then there exists $\alpha_u$ such that
$\psi^2_{\ddi}(u-u_{E,\kappa,\alpha_u})=0$.
\end{enumerate}

The first of these claims follows directly by Lemma \ref{oashteptpei}.

For the second claim, we write
$$0=\dd^2 \psi^1_{\di
\ddi}u=\left(\eta \dd-p \gamma \phip+\gamma \right)
\dd^2(l_E(u))\, .$$
Now
$$0 \neq v:=\dd^2(l_E(u)) \in qR[[q]]\, ,$$
so, by Theorem \ref{addeq}, there exists an integer $\kappa$ such
that $\eta \kappa+\gamma=0$ and
$v=u_{a,\kappa,\alpha}^{\mu^{(p)}}$ for some $\alpha \in R$. By
(\ref{sasesi}), $\dd l_E(u)=u_{a,\kappa,\alpha}^{\mu}$ so
$$u=e_E\left( \int u_{a,\kappa,\alpha}^{\mu} \frac{dq}{q}
\right)=u_{E,\kappa,\alpha}\, ,
$$ and the result follows.

For the third claim note that, as before, we have
$$0=\dd^2 \psi^1_{\di \ddi} u=
\frac{1}{\gamma} (\eta \dd-p \gamma \phip+\gamma) \psi^2_q u\, .
$$
By Theorem \ref{addeq}, $\psi^2_{\ddi}
u=a_0+u_{a,\kappa,\alpha}^{\mu^{(p)}}$ for some $\alpha \in R$ and
some $a_0 \in R$, with
$(\eta \dd-p \gamma \phip+\gamma) a_0=0$. The latter simply says
$p  \phi_p(a_0)=a_0$, and since $v_p(\phi_p(a_0))=v_p(a_0)$, we must have
$a_0=0$. Hence $\psi^2_{\ddi} u=u_{a,\kappa,\alpha}^{\mu^{(p)}}$.
On the other hand,
$$\psi^2_{\ddi} u_{E,\kappa,\alpha/\gamma \kappa}=
\gamma \kappa u^{\mu^{(p)}}_{a,\kappa,\alpha/\gamma \kappa}=
u^{\mu^{(p)}}_{a,\kappa,\alpha}.$$
and so $\psi^2_{\ddi}(u-u_{E,\kappa,\alpha/\gamma \kappa})=0$, completing the
proof.
\qed

\begin{remark}
The condition in Theorem \ref{unde} that the quotient
of the Kodaira-Spencer classes has the form
$$\frac{\eta}{\gamma} =-\frac{1}{\kappa}$$
for some integer
$\kappa \geq 1$ coprime to $p$ has a nice interpretation in terms of wave
lengths (in the sense of Remark \ref{funnyterm}). The condition is equivalent
to  saying that
the wave length of the parameter $\beta q$ of the Tate curve
belongs to
$$\left\{ -\frac{c}{\gamma},-\frac{2 c}{\gamma},-\frac{3
c}{\gamma},\ldots\right\} \cap \bZ_p^{\times}\, .
$$ This is, again, a
``quantization'' condition.
\end{remark}

\begin{corollary}
Under the hypotheses of Theorem {\rm \ref{unde}} and the
assumption that $\frac{\eta}{\gamma}=-\frac{1}{\kappa}$ for some
integer $\kappa$, let $u \in \cU_1^1$. Then  we have that the
series $\overline{\psi_q^2 u} \in k[[q]]$ is integral over $k[q]$,
and the field extension $k(q) \subset k(q,\overline{\psi_q^2 u})$
is Abelian with Galois group killed by $p$.
\end{corollary}

{\it Proof}. This follows immediately by Theorem
\ref{unde} and Lemmas \ref{irinusescoala} and \ref{gainuli}, respectively.
\qed

The next Corollary says that the ``boundary value problem at
$q=0$" is well posed.

\begin{corollary}
\label{vainiciel}
Under the hypotheses of Theorem {\rm \ref{unde}} and the assumption that
$\frac{\eta}{\gamma}=-\frac{1}{\kappa}$ for some integer $\kappa$, then
for any $\alpha \in R$ there exists a unique $u \in \cU_{1}^1$ such
that $B_\kappa^0 u=\alpha$.
\end{corollary}

Next we study $\D$-characters of order $2$ of the form
$\psi_{\ddi}^2+\lambda \psi_{\di}^2$.

\begin{theorem}
Let $E_{\beta q}$ be the Tate curve with parameter $\beta q$ for
$\beta$ not a root of unity. Let $\cU^2_{\ra}$ be the set
of solutions of the $\D$-character $\psi_{\ddi}^2+\lambda
\psi_{\di}^2$ in $A=R((q))\h$, where $\lambda \in R$. We set
$$
\begin{array}{rcl}
\cU^2_{0} & := & \cU^2_{\ra} \cap E_{0}(A)\, , \\
\cU^2_{1} & := & \cU^2_{\ra} \cap E_{1}(A)\, .
\end{array}
$$
Then the following hold:
\begin{enumerate}
\item[1)] If $\frac{\eta}{\gamma}=-\frac{1}{\kappa}$ and $\lambda=\kappa^3$
for some integer $\kappa \geq 1$ coprime to $p$, and if
$u_{E,\kappa,\alpha}$ is the basic series of 
$\psi^1_{pq}$ (cf. (\ref{filmm})), then
$$
\begin{array}{rcl}
  \cU^2_{\ra} & = & \cU^2_{0}+\cU^2_{1}\, , \\
\cU^2_{1} & = & \{u_{E,\kappa,\alpha}\, | \; \alpha \in R\} \, .
\end{array}
$$
\item[2)] If  $\frac{\lambda \eta^{\phi}}{\gamma} \not\in \{-n^2\, | \; n \in
\bZ\}$, then $\cU^2_{1}=0$.
\end{enumerate}
\end{theorem}

{\it Proof}. We prove assertion 1). The fact that $\cU^2_{1}$ consists of
exactly the $u_{E,\kappa,\alpha}$'s follows directly by Theorem \ref{unde}
and Proposition \ref{pomu}. We just need to observe that,
by Theorem \ref{addeq}, the map
$$(\gamma \eta^2\dd+\gamma^2 \eta\phip-\gamma^2\eta):qR[[q]]\ra qR[[q]]
$$
is injective. Now if $u \in \cU^2_{\ra}$ then, by Proposition
\ref{pomu}, we have
$$(\gamma \eta^2 \dd+\gamma^2 \eta \phip-\gamma^2 \eta) \psi^1_{\di \ddi}u=0
\, .
$$
By Theorem  \ref{addeq}, $\psi^1_{\di \ddi}u \in R$. Again, by
Proposition \ref{pomu}, we get
$$0=\gamma \dd^2 \psi^1_{\di \ddi} u=(\eta \dd-p\gamma
\phip+\gamma) \psi^2_{\ddi}u\, .
$$
By Theorem \ref{addeq} there is an
$\alpha \in R$ such that
$$\psi^2_{\ddi}u =u_{a,\kappa,\alpha}^{\mu^{(p)}}\, .
$$
On the other hand, as in the proof of Theorem \ref{unde}, we have that
$u_{a,\kappa,\alpha}^{\mu^{(p)}}=\psi^2_{\ddi}
u_{E,\kappa,\alpha/\gamma \kappa}$. Hence
$\psi^2_{\ddi}(u-u_{E,\kappa,\alpha/\gamma \kappa})=0$, and the assertion
is proved.

In order to prove assertion 2), let us assume $u \in qR[[q]]$ is such that
$(\psi^2_{\ddi}+\lambda \psi^2_{\di})u=0$. Then
$$
\begin{array}{rcl}
0 = \dd(\psi^2_{\ddi}+\lambda \psi^2_{\di})u & = & \dd\left[ \gamma
\dd^2+\frac{\lambda}{p} \left( \eta
\phip^2-(\eta^{\phi}+p\eta)\phip+p\eta^{\phi} \right) \right]l(u)\\
 & = & [\gamma \dd^2+p \lambda \eta
\phip^2-\lambda(\eta^{\phi}+p\eta)\phip+\lambda\eta^{\phi}]\dd
l(u)\, .
\end{array}
$$
By our assumption, $\frac{\lambda
\eta^{\phi}}{\gamma}+n^2\neq 0$ for all $n \in \bZ$. By Theorem
\ref{addeq},
$$\gamma \dd^2+p \lambda \eta
\phip^2-\lambda(\eta^{\phi}+p\eta)\phip+\lambda\eta^{\phi}
$$
has no non-zero solution in $qR[[q]]$. Since $\dd l(u) \in qR[[q]]$, we
get $\dd l(u)=0$ so $l(u) \in K \cap qK[[q]]=0$ so $u=0$.
\qed

Next we address the boundary value problem at $q \neq 0$ for
$\psi^1_{\di \ddi}$. For $q_0 \in pR$ and $\beta \in R^{\times}$
not a root of unity, we let ${\mathcal E}_{\beta q_0}$ be the
(non-smooth) plane projective curve that in the affine plane is given
by the equation
$$y^2=x^3-\frac{1}{48}E_4(\beta q_0)x+\frac{1}{864}E_6(\beta q_0)\, .$$
There is a commutative diagram of groups
\begin{equation}
\label{pichiu}
\begin{array}{ccc}
qR[[q]] & \stackrel{\iota}{\ra} & {\mathcal E}_{\beta q}(R[[q]])\\
\da &  & \da \pi_{q_0}\\
pR & \stackrel{\iota_{q_0}}{\ra} & {\mathcal E}_{\beta q_0}(R)
\end{array}\, ,
\end{equation}
where the vertical arrows are induced by the ring homomorphism
$$\begin{array}{rcl}
R[[q]] & \ra & R \\ u(q) & \mapsto & u(q_0)
\end{array}\, ,
$$
$\iota_{q_0}$ sends $pa$ into $(pa,(pa)^3-\frac{1}{48}E_4(\beta q_0)(pa)^7+
\cdots)$ for $a \in R$, and the image of $\iota_{q_0}$ is ${\mathcal
E}_{\beta q_0}(pR)$, the preimage of $0$ via the reduction modulo
$p$ mapping
\begin{equation}
\label{petro} {\mathcal E}_{\beta q_0}(R) \ra {\mathcal E}_{\beta
q_0}(k)\, .
\end{equation}
So we have $\iota_{q_0}(pR)={\mathcal E}_{\beta q_0}(pR)$.
We note that we have usually identified $\iota(qR[[q]])$ with
$qR[[q]]$, but in order to avoid the potential confusion produced by the choice
of $q_0$, we will {\bf not} identify $\iota_{q_0}(pR)$ with $pR$.

Let ${\mathcal E}_{\beta q_0}'(R)$ be the pull-back by the map
(\ref{petro}) of the locus ${\mathcal E}^{reg}_{\beta q_0}(k)$ of
regular points on the cubic ${\mathcal E}_{\beta q_0}(k)$. So we
have inclusions of groups
$${\mathcal E}_{\beta q_0}(pR) \subset {\mathcal E}'_{\beta q_0}(R)
\subset {\mathcal E}_{\beta q_0}(R)\, .
$$
Let us recall some facts about torsion points on Tate curves. There is a
natural injective homomorphism
\begin{equation} \label{paru}
\begin{array}{rcl}
\tau:\bmu(R) & \ra & {\mathcal E}_{\beta q}(R[[q]])\\
v & \mapsto & (X(q,v)+\frac{1}{12},-Y(q,v)-\frac{1}{2}X(q,v))\, ,
\quad  v \neq 1\, ,\\
1 & \mapsto & \infty\, ,
\end{array}
\end{equation}
given by the series in (\ref{uglys}). (The formula makes
sense because if $1 \neq v \in \bmu(R)$ then $1-v \in
R^{\times}$.) The composition of the homomorphism in (\ref{paru}) with
the specialization map $\pi_{q_0}$ in (\ref{pichiu}) gives a homomorphism
$$\tau_{q_0}:\bmu(R) \ra {\mathcal E}_{\beta q_0}(R)\, $$
The composition of $\tau_{q_0}$ with the reduction mod $p$ mapping
${\mathcal E}_{\beta q_0}(R) \ra {\mathcal E}_{\beta q_0}(k)$
is the map
$$\begin{array}{rcl}
\bmu(R) & \ra & {\mathcal E}_{\beta q_0}(k)\\
\zeta & \mapsto & \left( \frac{\zeta}{(1-\zeta)^2}+\frac{1}{12}\
\text{mod}\ p, -\frac{\zeta^2}{(1-\zeta)^3}-
\frac{\zeta}{2(1-\zeta)^2}
\ \text{mod}\ p \right)\, , \; \zeta \neq 1\, , \\
1 & \mapsto & \infty\, ,
\end{array}
$$
which is an isomorphism of $\bmu(R) \simeq k^{\times}$ onto
${\mathcal E}_{\beta q_0}^{reg}(k)$. We conclude that the mapping $\tau_{q_0}$
above is an injective map, and any point $P \in {\mathcal E}'_{\beta q_0}(R)$
can be written uniquely as
$$P=\iota_{q_0}(pa)+\tau_{q_0}(\zeta)\, ,$$
where $a \in R$ and $\zeta \in \bmu(R)$.

\begin{corollary}
\label{titoo}
Let $\cU^1_{\ra}$ be the set of solutions of $\psi^1_{\di \ddi}$ in
$R((q))\wh$, and let
$$\cU'_{\ra}:=\tau(\bmu(R)) \cdot \cU^1_{1} \subset \cU^1_{\ra}\, .$$
 Assume
$\frac{\eta}{\gamma}=-\frac{1}{\kappa}$ with $\kappa \geq 1$ an
integer coprime to $p$. Then the following hold:
\begin{enumerate}
\item[1)]  For any  $q_0 \in p^{\nu}R^{\times}$ with $\nu \geq 1$ and any
$g \in \iota_{q_0}(p^{\kappa \nu} R) \subset
\iota_{q_0}(pR)={\mathcal E}_{\beta q_0}(pR)$, there exists a
unique $u \in \cU^1_{1}=\iota(qR[[q]])=qR[[q]]$ such that
$u(q_0)=g$.
\item[2)] Assume $\kappa=1$. Then for any $q_0 \in pR^{\times}$ and any
$g \in {\mathcal E}'_{\beta q_0}(R)$, there exists a unique $u \in
\cU'_{\ra}$ such that $u(q_0)=g$.
\end{enumerate}
\end{corollary}

{\it Proof}. Let us  prove assertion 1).  By Theorem \ref{unde} and
Equations
(\ref{filmm}) and (\ref{nuilacasino}),
 it is enough to show that the mapping $R \ra p^{\kappa
\nu}R$ defined by
$$\alpha \mapsto e_E \left( \sum_{j \geq 0} (-1)^j
\frac{\alpha^{\phi^j}}{\kappa p^j F_j(p)} q_0^{\kappa p^j}
\right)$$
is a bijection. Let us recall that $e_E:p^N R \ra
\iota_{q_0}(p^N R)$ is an isomorphism for all $N$. So it is enough
to show that the map $R \ra p^{\kappa \nu}R$ given by
$$ \alpha \mapsto  \sum_{j \geq 0} (-1)^j
\frac{\alpha^{\phi^j}}{\kappa p^j F_j(p)} q_0^{\kappa p^j}\, ,$$
is a bijection. But this is clear by Lemma \ref{ajutator}.

For the proof of assertion 2), we let $g \in {\mathcal E}_{\beta
q_0}'(R)$, and write $g=\tau_{q_0}(\zeta_0)+w_0$, $w_0 \in
{\mathcal E}_{\beta q_0}(pR)$. By assertion 1), there exists $w \in
\cU_{1}^1$ such that $w(q_0)=w_0$. Let us set $u:=\tau(\zeta_0)+w$. Then
$u(q_0)=\tau_{q_0}(\zeta_0)+w_0=g$, which completes the proof of the existence
part of the assertion. In order to prove uniqueness, we let
$u_1= \tau(\zeta_1)+w_1$ with $w_1 \in \cU^1_{1}$ such that
$u_1(q_0)=g$. Then
$$\tau_{q_0}(\zeta_1)+w_1(q_0)=\tau_{q_0}(\zeta_0)+w_0\, .
$$
This implies that $\zeta_1=\zeta_0$ and $w_1(q_0)=w_0$. By the
uniqueness in assertion 1), we get $w_1=w$. We conclude that
$u_1=u$, and we are done.
\qed

The following Corollary is concerned with the inhomogeneous
equation $\psi^1_{pq} u=\varphi$, and it is an immediate consequence
of Lemmas \ref{oashteptpei}, \ref{vainiciel}, and \ref{irinusescoala}

\begin{corollary}
Under the hypotheses of Theorem {\rm \ref{unde}},  assume that
$\frac{\eta}{\gamma}=-\frac{1}{\kappa}$ for some integer $\kappa$,
and that $\varphi \in qR[[q]]$ is a series whose support is
contained in the set $\cK'$ of totally non-characteristic integers
of $\psi^1_{pq}$. Then the equation $\psi^1_{pq} u=\varphi$ has a
unique solution $u \in E_1(A)$ such that the support of $\psi^2_q
u$ does not contain $\kappa$. Moreover, if $\bar{\varphi} \in
k[q]$ the series $\overline{\psi^2_q u} \in k[[q]]$ is integral
over $k[q]$ and the field extension $k(q) \subset
k(q,\overline{\psi^2_q u})$ is Abelian with Galois group killed by
$p$.
\end{corollary}

\begin{remark}
Let us assume that $\frac{\eta}{\gamma}=-1$.
By  Corollary \ref{titoo}, for any $q_0 \in pR^{\times}$
the group homomorphism
$$\begin{array}{rcl}
S_{q_0}:\bmu(R) \times R & \ra & {\mathcal E}'_{\beta q_0}(R) \\
(\xi,\alpha) & \mapsto &\tau_{q_0}(\xi)+\iota_{q_0}(u_{E,1,\alpha}(q_0))
\end{array}
$$
is an isomorphism. So, for any $q_1,q_2 \in pR^{\times}$ we have a
group isomorphism
$$
S_{q_1, q_2}:=S_{q_2} \circ S_{q_1}^{-1}:{\mathcal E}'_{\beta q_1}(R)
\ra {\mathcal E}'_{\beta q_2}(R)\, .
$$
The latter mapping can be viewed as the ``propagator'' attached to
$\psi^1_{\di \ddi}$.

It is natural to make the propagator act as an endomorphism of a given group.
We accomplish this by defining a new {\it propagator}
$$S_{q_1,q_2}^{q_0}:=\Gamma_{q_0} \circ \Gamma_{q_2}^{-1} \circ
S_{q_1,q_2} \circ \Gamma_{q_1} \circ \Gamma_{q_0}^{-1}: {\mathcal
E}'_{\beta q_0}(R) \ra {\mathcal E}'_{\beta q_0}(R)\, ,$$
where $\Gamma_{q_i}$ is the group isomorphism
$$\begin{array}{rcl}
\bmu(R) \times R & \stackrel{\Gamma_{q_i}}{\ra } &
{\mathcal E}'_{\beta q_i}(R) \\
(\xi,a) & \mapsto & \iota_{q_i}(pa)+\tau_{q_i}(\xi)
\end{array}\, .
$$
Then it is easy to see that for $\zeta_1,\zeta_2 \in \bmu(R)$ we
have that
$$S_{q_0,\zeta_1 \zeta_2 q_0}^{q_0}=S_{q_0,\zeta_2 q_0}^{q_0}
\circ S_{q_0,\zeta_1  q_0}^{q_0}\, ,
$$
which, once again, we view as a (weak) form of the ``Huygens principle.''
\end{remark}

\subsection{Elliptic curves over $R$}
In this section we consider an elliptic curve $E$ over
$A:=R$ defined by
$$f=y^2-(x^3+a_4x+a_6)\, ,$$
where $a_4,a_6 \in R$, and equip it with the $1$-form
$$\omega=\frac{dx}{y}\, .$$
We use the notation and discussion in Example
\ref{cucurigu} that applies to $E$ over $R$.

Clearly, $f^1_{\ddi}=0$. Let us fix in what follows the \'{e}tale
coordinate $T=-\frac{x}{y}$, and denote by $l(T)=l_E(T) \in K[[T]]$
the logarithm of the formal group $\cF_E$ of $E$ with respect to $T$.

\begin{proposition}
\label{logder} The following holds:
\begin{enumerate}
\item[1)] The image of $\psi^1_{\ddi}$ in
$$\cO\left(J^1_{\di \ddi}\left(\frac{A[x,y]}{(f)} \right) \right)=
\frac{A[x,y,\d x, \d y, \dd x, \dd y]\h}{(f,\d f,\dd f)}$$
is equal to ${\displaystyle \frac{\dd x}{y}}$.
\item[2)] The image of $\psi^1_{\ddi}$ in $A[[T]][\d T, \dd T]\h$ is equal
to ${\displaystyle \dd l(T)=\frac{dl}{dT}(T) \cdot \dd T}$.
\end{enumerate}
\end{proposition}

The first assertion above says that the $\D$-character $\psi^1_{\ddi}$
attached to $(E,\omega)$ coincides with the {\it Kolchin's logarithmic
derivative} \cite{kolchin}.

{\it Proof}. We first check that $\frac{\dd x}{y}$ comes from a $\D$-character.
For that, it is enough to show that its image
$\left( \frac{\dd x}{y}\right) (T,\dd T)$ in $A[[T]][\d T,\dd T]\h$ defines a
homomorphism from the formal group of $J^1_{\di \ddi}(E)$ to the
formal group of $\bG_a$. This is so because ``a partially defined map,
which generically is a homomorphism, is an everywhere defined
homomorphism.'' We now recall that $\omega$ and $l$ are related by the
equation
$$\omega=\frac{dl}{dT} (T) \cdot dT\, .$$
Consequently, if we denote by $x(T)$ and $y(T)$ the images of $x$ and
$y$ in $R((T))$, respectively,  we have that
$$\frac{dl}{dT} (T) \cdot dT=\left(\frac{dx}{y} \right) (T)=
(y(T))^{-1} \frac{dx}{dT} (T)\cdot  dT\, .
$$
Hence
\begin{equation}
\label{stelutza} \left( \frac{\dd x}{y} \right) (T,\dd
T)=(y(T))^{-1} \frac{dx}{dT}(T) \cdot \dd T=\frac{dl}{dT}(T) \cdot
\dd T=\dd l(T)\, .
\end{equation}
But clearly $\dd l(T)$ defines a homomorphism at the level of
formal groups. Hence $\frac{\dd x}{y}$ is a $\D$-character of $E$.
Now, (\ref{stelutza}) also shows that
$$\psi^1_{\ddi}-\frac{\dd x}{y} \in A[[T]]\, ,$$
and thus, $\psi^1_{\ddi}-\frac{\dd x}{y}$ defines a homomorphism
$\hat{E} \ra \hat{\bG}_a$. Therefore, $\psi^1_{\ddi}-\frac{\dd x}{y}=0$
and assertion 1) is proved.

The second assertion follows from the first in combination with
(\ref{stelutza}).
\qed

\begin{remark}
As in the previous subsection, we identify $q^{\pm
1}R[[q^{\pm 1}]]$ with its image $E(q^{\pm 1}R[[q^{\pm 1}]])$
in $E(R[[q^{\pm 1}]])$ under the embedding
$$\iota:q^{\pm 1} R[[q^{\pm 1}]] \ra E(q^{\pm 1}R[[q^{\pm 1}]])$$
given in $(T,W)$- coordinates by
$$u \mapsto (u,u^3+a_4 u^7+\cdots)\, .$$
\end{remark}

In the sequel, we fix an elliptic curve $E/R$. We define next
the {\it characteristic polynomial} of a $\D$-character,  {\it
non-degenerate} $\D$-characters, their {\it characteristic  integers},
and {\it basic series}. We distinguish the two cases
$f^1_{\di} \neq 0$ and $f^1_{\di} = 0$, respectively.

If $f^1_{\di} \neq 0$, by \cite{book}, p. 197,
$v_p(f^2_{\di}) \geq v_p(f^1_{\di})$. We may therefore
consider the $\D$-character
$$\psi^2_{\di,nor} :=\frac{1}{f^1_{\di}} \psi^2_{\di} \in
\bX^2_{\di \ddi}(E)\, .$$
Its image in $R[[T]][\d T, \d^2 T]\h$ is
$$\psi^2_{\di,nor}=\frac{1}{p}(\phip^2+\gamma_1 \phip+p
\gamma_0)l(T),$$ where $$\gamma_1:=-\frac{f^2_{\di}}{f^1_{\di}}\in
R,\ \ \gamma_0:= \frac{(f^1_{\di})^{\phi}}{f^1_{\di}}\in
R^{\times}.$$ By Proposition \ref{ileana1}, any $\D$-character of
$E$ is a $K$-multiple of a $\D$-character of the form
\begin{equation}
\label{mau} \psi_E:=\nu(\phip,\dd)  \psi^1_{\ddi}+\lambda(\phip)
\psi^2_{\di,nor}\, ,
\end{equation}
 where $\nu(\xi_p,\xi_q) \in R[\xi_p,\xi_q]$ and $\lambda(\xi_p) \in R[\xi_p]$.
 The Picard-Fuchs symbol of $\psi_E$ with respect to $T$
 is easily easily seen to be
$$\sigma(\xi_p,\xi_q)=p \nu(\xi_p,\xi_q)\xi_q
 +\lambda(\xi_p)(\xi_p^2+\gamma_1 \xi_p+p\gamma_0)\, .
$$
Hence, the Fr\'{e}chet symbol of $\psi_E$ with respect to $\omega$
is
$$\theta(\xi_p,\xi_q)=\frac{\sigma(p\xi_p,\xi_q)}{p}=
\nu(p\xi_p,\xi_q)\xi_q+\lambda(p\xi_p)(p\xi_p^2+\gamma_1
\xi_p+\gamma_0)\, .
$$

\begin{definition}
Let us assume that $E/R$ has $f^1_{\di} \neq 0$, and let $\psi_E$ be
a  $\D$-character of $E$ of the form in (\ref{mau}). We define the {\it
characteristic polynomial} $\mu(\xi_p,\xi_q)$  of $\psi_E$ to be the
Fr\'{e}chet symbol $\theta(\xi_p,\xi_q)$ of $\psi_E$ with respect to $\omega$.
We say that the $\D$-character $\psi_E$ is {\it non-degenerate} if
$\mu(0,0)\in R^{\times}$, or equivalently, if $\lambda(0)\in
 R^{\times}$. For a non-degenerate character $\psi_E$, we define
its {\it characteristic integers} to be the integers $\kappa$
such that $\mu(0,\kappa) =0$ (so any such $\kappa$ is coprime to
$p$), and denote by $\cK$ the set of all such. We call $\kappa$ a
{\it totally non-characteristic} integer if $\kappa \not\equiv 0$ mod $p$ and
$\mu(0,\kappa) \not\equiv 0$ mod $p$. The set of all totally
non-characteristic integers is denoted  by $\cK'$.
For $0 \neq \kappa \in \bZ$ and $\alpha \in R$, we define the
{\it basic series} by
\begin{equation}
\label{caino} u_{E,\kappa,\alpha}=e_E\left( \int
u_{a,\kappa,\alpha} \frac{dq}{q} \right) \in q^{\pm 1}K[[q^{\pm
1}]]\, ,
\end{equation}
where
\begin{equation}
\label{cain} u_{a,\kappa,\alpha}:=u_{a,\kappa,\alpha}^{\mu};
\end{equation}
Cf. (\ref{ua}).
\end{definition}

Let us now consider now the case where $f^1_{\di} =0$. We may then look at the
$\D$-character
$\psi^1_{\di} \in \bX^1_{\di}(E)$, cf. (\ref{psi1p}). Its image in
$R[[T]][\d T]\h$ is
$$\psi^1_{\di}=\frac{1}{p}(\phip +p \gamma_0) l(T)\, ,$$
where $\gamma_0 \in R$; cf. \cite{book}, Remark 7.21. (We
actually have $\gamma_0 \in R^{\times}$ if the cubic defining $E$
has coefficients in $\bZ_p$; cf. \cite{frob}, Theorem 1.10.) By
Proposition \ref{ileana2}, any $\D$-character of $E$ is a
$K$-multiple of a $\D$-character of the form
\begin{equation}
\label{drrr}
\psi_E:=\nu(\phip,\dd) \psi^1_{\ddi}+\lambda(\phip)
\psi^1_{\di}\, ,
\end{equation}
where $\nu(\xi_p,\xi_q) \in R[\xi_p,\xi_q]$ and $\lambda(\xi_p)
\in R[\xi_p]$. We easily see that the Picard-Fuchs symbol of $\psi_E$ is
given by
$$\sigma(\xi_p,\xi_q)=p\nu(\xi_p,\xi_q)\xi_q+\lambda(\xi_p)
(\xi_p+p\gamma_0)\, .
$$
Thus, the Fr\'{e}chet symbol at the origin (with respect to $dT$) is
$$\theta(\xi_p,\xi_q)=
\frac{\sigma(p\xi_p,\xi_q)}{p}=\nu(p\xi_p,\xi_q)\xi_q
+\lambda(p\xi_p)(\xi_p+\gamma_0)\, .
$$

\begin{definition}
Let us assume that $E/R$ has $f^1_{\di}=0$. (Recall that then
$E/R$ has ordinary reduction mod $p$; cf. \cite{book}, Corollary
8.89.), and let us fix a $\D$-character of $E$ of the form
$\psi_E$ as in {\rm (\ref{drrr})}. We define the {\it
characteristic polynomial} $\mu(\xi_p,\xi_q)$ of $\psi_E$ to be
the Fr\'{e}chet symbol $\theta(\xi_p,\xi_q)$ of $\psi_E$ with
respect to $\omega$. We say that the $\D$-character $\psi_E$ is
{\it non-degenerate} if $\mu(0,0) \in R^{\times}$, or
equivalently, if $\lambda(0)\in R^{\times}$ and $\gamma_0 \in
R^{\times}$. For a non-degenerate character $\psi_E$, we define
the {\it characteristic integers} to be the integers $\kappa$ such
that $\mu(0,\kappa) =0$ (so any such $\kappa$ is coprime to $p$),
and denote by $\cK$ the set of all such. We say that $\kappa$ is a
{\it totally non-characteristic} integers is $\kappa \not\equiv 0$
mod $p$
 and $\mu(0,\kappa) \not\equiv 0$ mod $p$. The set of all
totally non-characteristic integers is denoted by $\cK'$. For $0
\neq \kappa \in \bZ$ and $\alpha \in R$, we define the {\it basic
series}
\begin{equation}
\label{errr} u_{E,\kappa,\alpha}=e_E\left( \int
u_{a,\kappa,\alpha} \frac{dq}{q} \right) \in q^{\pm 1}K[[q^{\pm 1}]]\, ,
\end{equation}
 where
\begin{equation}
\label{rrrr} u_{a,\kappa,\alpha}:=u_{a,\kappa,\alpha}^{\mu};
\end{equation} cf. (\ref{ua}).
\end{definition}

In the sequel, we consder an elliptic curve $E/R$ without imposing
any a priori restriction on the vanishing of $f^1_{\di}$.

\begin{example}
 Let us assume that $\nu=1$, $\lambda \in R^{\times}$, that is to say,
$\psi_E$ is either $\psi^1_q+\lambda \psi^2_{p,nor}$ or
$\psi^1_q+\lambda \psi^1_p$ according as  $f^1_p \neq 0$ or $f^1_p
=0$, respectively.
 Then the characteristic polynomial is unmixed (provided that
 $\gamma_1 \in R^{\times}$ in case $f^1_p \neq 0$),
$$
\cK=\{-\lambda \gamma_0\} \cap \bZ\, ,
$$
and $\cK \neq \emptyset$ if, and only if, $\lambda \gamma_0 \in
\bZ$. In this case, $\psi_E$ should be viewed as an analogue of
either the heat equation or the convection equation according as
$f^1_{\di} \neq 0$ or $f^1_{\di}=0$, respectively.
\end{example}

\begin{example}
Let us now assume that $\nu=\xi_q$ and $\lambda \in R^{\times}$,
that is to say, $\psi_E$ is either $\dd \psi^1_q+\lambda
\psi^2_{p,nor}$ or $\dd \psi^1_q+\lambda \psi^1_p$ according as
$f^1_p \neq 0$ or $f^1_p =0$, respectively. Then the
characteristic polynomial is unmixed (provided that
 $\gamma_1 \in R^{\times}$ in case $f^1_p \neq 0$),
$$
\cK=\{\pm \sqrt{-\lambda \gamma_0}\} \cap \bZ\, ,$$ and $\cK \neq
\emptyset$ if, and only if, $-\lambda \gamma_0$ is  a perfect
square. If this the case, $\psi_E$ should be viewed as analogue of
either the wave equation or the sideways heat equation according
as $f^1_{\di} \neq 0$ or $f^1_{\di}=0$, respectively.
\end{example}

\begin{example}
Let us consider the ``simplest'' of the energy functions in the case
where $f^1_p=0$:
$$
H=a(\psi^1_q)^2+2b \psi^1_p \psi^1_q +c(\psi_p^1)^2 \, .
$$
Here, $a,b,c \in R$.

A computation essentially identical to the one in Example \ref{disom} leads
to the following formula for the Euler-Lagrange equation
$\epsilon^1_{H,\partial}$ attached to $H$
and the vector field $\partial$ dual to $\omega$:
$$
\epsilon^1_{H,\partial}=(-2a^{\phi}\phip \dd-2b^{\phi}
\phip^2+2b)\psi^1_q+ (2c^{\phi}\gamma_0^{\phi} \phip
+2c)\psi^1_p\, .
$$
Hence, the characteristic polynomial of
$\epsilon^1_{H,\partial}$ is
$$\mu(\xi_p,\xi_q)=(-2a^{\phi}p \xi_p \xi_q-2b^{\phi}p^2
\xi_p^2+2b)\xi_q+ (2c^{\phi}\gamma_0^{\phi}p
\xi_p+2c)(\xi_p+\gamma_0)\, ,
$$
so the $\D$-character $\epsilon^1_{H,\partial}$ is non-degenerate if,
and only if, $c \in R^{\times}$ and $\gamma_0 \in R^{\times}$. Moreover,
 the set of characteristic  integers is
$$
\cK=\{-\frac{c \gamma_0}{b} \}\cap \bZ\, .
$$

On the other hand, when $f^1_p \neq 0$ we consider the
``simplest'' energy function
$$
H=a(\psi^1_q)^2+2b \psi^1_q \psi^2_{p,nor} +c(\psi_{p,nor}^2)^2\, ,
$$
for $a,b,c \in R$. Then we have the following values for the
Fr\'{e}chet symbols:
$$
\begin{array}{rcl}
\theta_{\psi^1_q,\omega} & = & \xi_q\\
\theta_{\psi^2_{p,nor},\omega} & = & p \xi_p^2+\gamma_1 \xi_p
+\gamma_0
\end{array}\, ,
$$
and we get the following formula for the Euler-Lagrange equation
attached to $H$ and the vector field $\partial$ dual to $\omega$:
$$
\begin{array}{rcl}
\epsilon^2_{H,\partial} & = & (-2b^{\phi^2}p \phip^3 -2a^{\phi^2}
\phip^2 \dd)\psi^1_q \\
  & \ & + [2b^{\phi^2}(\gamma_0^{\phi^2}-\gamma_1^{\phi})
\phip^2+ (2b^{\phi}\gamma_1^{\phi}-2b^{\phi^2}
\gamma_0^{\phi})\phip+ 2bp] \psi^1_q\\
  &  & + (2c^{\phi^2} \gamma_0^{\phi^2} \phip^2
+2c^{\phi}\gamma_1^{\phi} \phip +2cp) \psi^2_{p,nor}\, .
\end{array}
$$
In particular $\epsilon^2_{H,\partial}$ is degenerate for all
values of $a,b,c$. More is true, actually:
$\epsilon^2_{H,\partial}$ is not a $K$-multiple of a
non-degenerate $\D$-character.
\end{example}

\begin{lemma}
Let $\kappa \in \bZ\backslash p\bZ$. Then we have
$u_{E,\kappa,\alpha} \in q^{\pm 1}R[[q^{\pm 1}]]$ for all $\alpha
\in R$, and
$$\begin{array}{rcl}
R & \ra & R[[q^{\pm 1}]] \\ \alpha & \mapsto & u_{E,\kappa,\alpha}
\end{array}
$$
is a pseudo $\d$-polynomial mapping. Consequently,
$$\begin{array}{rcl}
R & \ra & E(R[[q^{\pm 1}]]) \\
\alpha & \mapsto & \iota(u_{E,\kappa,\alpha})
\end{array}
$$
is also a pseudo $\d$-polynomial mapping.
\end{lemma}

{\it Proof}. As an element of $R[[T]][\d T, \dd T,\d^2 T,\d \dd T,\dd^2T]\wh$,
the $\D$-character $\psi_E$ coincides with
$$\psi_E=\left( \nu(\phip,\dd)\dd + \frac{\lambda(\phip)}{p}
 \phip^2 + \frac{\lambda(\phip)}{p}
\gamma_1 \phip+\lambda(\phip) \gamma_0 \right)l(T)
$$
if $f^1_{\di} \neq 0$, and it coincides with
$$\psi_E=(\nu(\phip,\dd)\dd+\frac{\lambda(\phip)}{p}
\phip+\lambda(\phip) \gamma_0)l(T)
$$
if $f^1_{\di}=0$. In the first case, by
Proposition \ref{logder} we have that
\begin{equation} \label{casc}
\begin{array}{rcl} \dd \psi_E & = & (
\nu(p\phip,\dd)\dd
 + \lambda(p\phip) p
\phip^2 +\lambda(p\phip) \gamma_1 \phip
+ \lambda(p\phip) \gamma_0) \dd l(T)\\
\  & = & (\nu(p\phip,\dd)\dd + \lambda(p\phip) p \phip^2
+\lambda(p\phip) \gamma_1 \phip + \lambda(p\phip) \gamma_0)
\psi^1_{\ddi}\, .
\end{array}
\end{equation}
Similarly, if $f^1_{\di}=0$ we have
\begin{equation}
\label{mz23}
\begin{array}{rcl} \dd \psi_E & = & (\nu(p\phip,\dd)\dd +
\lambda(p\phip) \phip
+ \lambda(p\phip) \gamma_0) \dd l(T)\\
 & = & (\nu(p\phip,\dd)\dd + \lambda(p\phip)  \phip +
\lambda(p\phip) \gamma_0) \psi^1_{\ddi}\, .
\end{array}
\end{equation}
Thus, if  $\psi_E u=0$, we have $\psi^1_{\ddi}u$ a solution of
$$
\nu(p\phip,\dd)\dd+\lambda(p\phip)
 p \phip^2+ \lambda(p\phip) \gamma_1 \phip +\lambda(p\phip) \gamma_0
$$
when $f^1_{\di}\neq 0$, or a solution of
$$\nu(p\phip,\dd)\dd+ \lambda(p\phip)  \phip +\lambda(p\phip) \gamma_0
$$
when $f^1_{\di}= 0$.

Let $\alpha \in R$, and set
$$h:=\int u_{a,\kappa,\alpha} \frac{dq}{q} \in qK[[q]]\, ,$$
with $u_{a,\kappa,\alpha}$ as in (\ref{cain}) or (\ref{rrrr}),
respectively. Note that
$$
\dd h=u_{a,\kappa,\alpha}\, .
$$
By Theorem  \ref{addeq}, when $f^1_{\di} \neq 0$ we get that
$$
\begin{array}{rcl}
\dd(\mu(0,\kappa)\alpha \kappa^{-1}q^{\kappa}) \! \! \! & = & \! \! \!
\mu(0,\kappa)\alpha q^{\kappa} \vspace{1mm} \\
 & = & \! \! \!
(\nu(p\phip,\dd)\dd \! + \! \lambda(p\phip) p \phip^2+ \lambda(p\phip)
\gamma_1 \phip +\lambda(p\phip)
\gamma_0) u_{a,\kappa,\alpha} \vspace{1mm} \\
 & = & \! \! \!  (\nu(p\phip,\dd)\dd+\lambda(p\phip) p \phip^2+
\lambda(p\phip) \gamma_1 \phip +\lambda(p\phip)
\gamma_0) \dd h \vspace{1mm} \\
 & = &\! \! \!  \dd \left(\nu(\phip,\dd)\dd+ \frac{\lambda(\phip)}{p}
\phip^2+ \frac{\lambda(\phip)}{p} \gamma_1 \phip +\lambda(\phip)
\gamma_0 \right)h
\end{array}\, .
$$
Hence
\begin{equation}
\label{davidstar}  \begin{array}{c} \left(\nu(\phip,\dd)\dd+
\frac{\lambda(\phip)}{p} \phip^2+ \frac{\lambda(\phip)}{p}
\gamma_1 \phip + \lambda(\phip)\gamma_0 \right)h
-\mu(0,\kappa)\alpha \kappa^{-1}q^{\kappa}\\
 \in R \cap q^{\pm 1}K[[q^{\pm 1}]]=0,\end{array}
\end{equation}
and we have that
\begin{equation}
\begin{array}{rcl}
 \lambda(\phip) \left(\frac{1}{p} \phip^2+ \frac{1}{p} \gamma_1
\phip + \gamma_0 \right)h  & = & -\nu(\phip,\dd)\dd h
+\mu(0,\kappa)\alpha \kappa^{-1}q^{\kappa}
\\ & = & -\nu(\phip,\dd)u_{a,\kappa,\alpha}
+\mu(0,\kappa)\alpha \kappa^{-1}q^{\kappa}\\
  &  \in & q^{\pm 1}R[[q^{\pm 1}]]\, .
\end{array}
\end{equation}
By Lemma \ref{floricica}, we get
\begin{equation}
\label{onestar} \left( \frac{1}{p} \phip^2+\frac{1}{p} \gamma_1
\phip+\gamma_0 \right) h \in q^{\pm 1}R[[q^{\pm 1}]] \subset
R[[q^{\pm 1}]]\, .
\end{equation}

Similarly, when $f^1_{\di}=0$, we get
\begin{equation}
\label{oozz} \left( \frac{1}{p} \phip+ \gamma_0 \right) h \in
R[[q^{\pm 1}]]\, .
\end{equation}

We claim that that if $f^1_{\di} \neq 0$ we have
\begin{equation}
\label{twostars} \left(\frac{1}{p} \phip^2+ \frac{1}{p} \gamma_1
\phip +\gamma_0 \right)l_E(q) \in R[[q^{\pm 1}]]\, ,
\end{equation}
where $l_E(q) \in K[[q^{\pm 1}]]$ is obtained from $l(T)=l_E(T)
\in K[[T]]$ by substitution of $q^{\pm 1}$ for $T$, and
$\phip:R[[q^{\pm 1}]] \ra R[[q^{\pm 1}]]$ is defined by
$\phip(q^{\pm 1})=q^{\pm p}$. In order to check that this
holds, it is sufficient to check that
\begin{equation}
\label{fourstars} (\phip^2+ \gamma_1 \phip + p \gamma_0)l_E(q^{\pm
1}) \in pR[[q^{\pm 1}]]\, .
\end{equation}
Now recall that
\begin{equation}
\label{threestars} (\phip^2+ \gamma_1 \phip + p \gamma_0)l_E(T)
\in pR[[T]][\d T, \d^2 T]\h\, ,
\end{equation}
where,  as usual, the mappings $R[[T]] \stackrel{\phip}{\ra}
R[[T]][\d T]\h \stackrel{\phip}{\ra}R[[T]][\d T, \d^2 T]\h$ are
defined by $\phip(T)=T^p+p \d T$, $\phip(\d T)=(\d T)^p+p\d^2 T$.
Taking the image (\ref{threestars}) under the unique $\d$-ring
homomorphism $R[[T]][\d T, \d^2 T]\h \ra R[[q^{\pm 1}]]$ that
sends $T$ into $q^{\pm 1}$, we get equality (\ref{fourstars}),
completing the verification that (\ref{twostars}) holds.

Similarly, when $f^1_{\di}=0$ we have that
\begin{equation}
\label{zzzt} \left( \frac{1}{p} \phip + \gamma_0 \right) l_E(q)
\in R[[q^{\pm 1}]]\, .
\end{equation}

If $f^1_{\di} \neq 0$, then by (\ref{onestar}), (\ref{twostars}), and
Hazewinkel's Functional Equation Lemma
\ref{haze}, it follows that $e(h) \in R[[q^{\pm 1}]]$. Since $e(h)
\in q^{\pm 1}K[[q^{\pm 1}]]$, we get
$$
u_{E,\kappa,\alpha}=e_E(h) \in q^{\pm 1}R[[q^{\pm 1}]]\, .
$$

Similarly, if $f^1_{\di}=0$, by (\ref{oozz}) and (\ref{zzzt}), we get
that
$$
u_{E,\kappa,\alpha} \in q^{\pm 1}R[[q^{\pm 1}]]\, .
$$

As in the proof of Lemma \ref{ol}, we see that
$$\begin{array}{rcl}
R & \ra & R[[q^{\pm 1}]] \\
\alpha & \mapsto & u_{E,\kappa,\alpha}
\end{array}
$$ is pseudo $\d$-polynomial mapping.
\qed

We have the following diagonalization result.

\begin{lemma}
\label{gainuseli}
\label{viyne} Let $\kappa \in \bZ \backslash p\bZ$, and $\alpha
\in R$. Then
$$\begin{array}{rcl}
\psi^1_q u_{E,\kappa,\alpha} & = & u_{a,\kappa,\alpha}\\
\psi_E u_{E,\kappa,\alpha} & = & \mu(0,\kappa)\alpha \kappa^{-1}
q^{\kappa}.\end{array}
$$
\end{lemma}

{\it Proof}. The first equality is clear.
Now, by (\ref{davidstar}), if $f^1_{\di}\neq 0$ we have that
$$
\begin{array}{rcl}
\psi_E u_{E,\kappa,\alpha} & = &  \left( \nu(\phip,\dd)\dd+
\frac{\lambda(\phip)}{p} \phip^2+\frac{\lambda(\phip)}{p} \gamma_1
\phip +\lambda(\phip) \gamma_0 \right) l_E(e_E(h))\\
 & = & \mu(0,\kappa) \alpha \kappa^{-1} q^{\kappa}\, .
\end{array}
$$
A similar argument can be given in the case $f^1_p=0$.
\qed

\begin{remark}
\begin{enumerate}
\item[1)] For  $\kappa \in \bZ\backslash p\bZ$ we have
\begin{equation}
\label{pie}
u_{E,\kappa,\zeta^{\kappa} \alpha}(q)=u_{E,\kappa,\alpha}(\zeta q)
\end{equation}
for all $\zeta \in \bmu(R)$. So if $\alpha=\sum_{i=0}^{\infty} m_i
\zeta_i^{\kappa}$, $\zeta_i \in \bmu(R)$, $m_i \in \bZ$, $v_p(m_i)
\ra \infty$, then
$$
u_{E,\kappa,\alpha}(q)
= \left[ \sum_{i=0}^{\infty} \right] [m_i](u_{E,\kappa,1}(\zeta_i
q))\, .
$$
If $f \in \bZ \bmu(R)\h$ is such that
$(f^{[\kappa]})^{\sh}=\alpha \in R$, then
$$
u_{E,\kappa,\alpha}=f \star u_{E,\kappa,1}\, .
$$
Note that $\{u_{E,\kappa,\alpha}\ |\ \alpha \in R\}$ is a $\bZ
\bmu(R)\h$-module (under convolution). This module structure comes
from  an $R$-module structure (still denoted by $\star$) by base
change via the morphism
$$\bZ \bmu(R)\h \stackrel{[\kappa]}{\ra}
\bZ \bmu(R)\h \stackrel{\sh}{\ra} R
$$
(cf. (\ref{notinjjj})), and the $R$-module $\{u_{E,\kappa,\alpha}\, | \;
\alpha \in R\}$ is free, with basis $u_{E,\kappa,1}$. Thus, for
$g \in \bZ \bmu(R)\h$, $\beta=(g^{[\kappa]})^{\sh}$, we have that
$$
\beta \star u_{E,\kappa,\alpha}=g \star u_{E,\kappa,\alpha}\, .
$$
In particular
$$
u_{E,\kappa,\alpha}  = \alpha \star u_{E,\kappa,1}\, .
$$
\item[2)] We recall (see (\ref{vocea})) the natural group homomorphisms
attached to $\psi_q^1$,
$$\begin{array}{ccc}
B_{\kappa}^0:E(R[[q^{\pm 1}]]) & \ra & R \\
B_{\kappa}^0 u & = & \Gamma_{\kappa} \psi^1_q u
\end{array}\, .$$
For $\kappa_1,\kappa_2 \in \bZ \backslash p\bZ$, we obtain that
\begin{equation}
\label{blleah} B_{\kappa_1}^0 u_{E,\kappa_2,\alpha}
=\Gamma_{\kappa_1}
 u^{\mu}_{a,\kappa_2,\alpha}=
\alpha \cdot \delta_{\kappa_1 \kappa_2}\, .
\end{equation}
\item[3)] Assume that $E$ is defined over $\bZ_p$ and $f^1_{\di}\neq 0$. Then
$f^1_{\di} \in \bZ_p$, so $\gamma_0=1$. Also, by the Introduction of
\cite{frob}, $\gamma_1 \in \bZ$. Thus, if in addition to $\alpha,
\lambda \in \bZ_{(p)}$ we have that the cubic defining $E$ has coefficients
in $\bZ_{(p)}$, then $u_{E,\kappa,\alpha} \in \bZ_{(p)}[[q^{\pm
1}]]$  for $\kappa \in \bZ \backslash p\bZ$.
\item[4)] Assume $E$ is defined over $\bZ_p$ and that $f^1_{\di}= 0$. Then
the Introduction in \cite{frob}, $\gamma_0$ is an integer in a quadratic
extension $F$ of $\bQ$. We view $F$ as embedded into $\bQ_p$, and set
$\cO_{(p)}:=F \cap \bZ_p$. Then, if in addition to $\alpha, \lambda
\in \bZ_{(p)}$ we have that the cubic defining $E$ has coefficients in
$\bZ_{(p)}$, we obtain that $u_{E,\kappa,\alpha} \in \cO_{(p)}[[q^{\pm
1}]]$ for $\kappa \in \bZ \backslash p\bZ$.
\end{enumerate}
\end{remark}

\begin{theorem}
\label{bloo} Let $E$ be an elliptic curve over $R$, $\psi_E$
be a non-degenerate $\D$-character of $E$, and $\cU_*$ be the
corresponding groups of solutions. If $\cK$ is the set of
characteristic integers of $\psi_E$, and $u_{E,\kappa,\alpha}$
be the basic series, then we have that
$$
\cU_{\pm 1}=\oplus_{\kappa \in \cK_{\pm}} \{u_{E,\kappa,\alpha}\ |\
\alpha \in R\}\, ,
$$
where $\oplus$ stands for the internal direct sum. In
particular, $\cU_{\pm 1}$ are free $R$-modules under convolution,
with bases
$
\{u_{E,\kappa,1}\ |\ \kappa \in \cK_{\pm}\}\, ,
$
respectively.
\end{theorem}

{\it Proof}.
 By Lemma \ref{viyne}, $u_{E,\kappa,\alpha} \in
\cU_{1}$. Conversely, if $u \in \cU_{1}$, that is to say, if $u \in
E(qR[[q]])$ and $\psi_E u=0$, we have that $\psi^1_{\ddi} u \in
qR[[q]]$ is a solution of
$$
\nu(p\phip,\dd)\dd+\lambda(p\phip) p
\phip^2+\lambda(p\phip) \gamma_1 \phip+\lambda(p\phip) \gamma_0\, ,
$$
or a solution of
$$
\nu(p\phip,\dd)\dd+\lambda(p\phip)
\phip+\lambda(p\phip) \gamma_0 \, ,
$$
if $f^1_{\di}\neq 0$ or $f^1_{\di}=0$, respectively. Cf. (\ref{casc})
 and (\ref{mz23}),  respectively. By Theorem
\ref{addeq}, we have
$$
\psi^1_{\ddi} u=\sum_{i=1}^s u^{\mu}_{a,\kappa_i,\alpha_i}
$$
for some $\alpha_i \in R$, where $\cK_+=\{\kappa_1,\ldots,\kappa_s\}$.
Therefore,
$$
u=e_E \left( \int \sum_{i=1}^s u^{\mu}_{a,\kappa_i,\alpha_i}
\frac{dq}{q} \right)=\sum_{i=1}^s u_{E,\kappa_i,\alpha_i}\, .
$$
This representation is unique because of (\ref{blleah}).

A similar argument works for $u \in \cU_{-1}$.
\qed

\begin{corollary}
Under the hypotheses of Theorem {\rm \ref{bloo}}, let $u \in
\cU_{\pm}$. Then the following hold:
\begin{enumerate} 
\item The series $\overline{\psi^1_q u} \in k[[q^{\pm1}]]$ is integral over
$k[q^{\pm 1}]$, and the field extension $k(q) \subset
k(q,\overline{\psi^1_q u})$ is Abelian with Galois group killed by
$p$. 
\item If the characteristic polynomial of $\psi_E$ is unmixed
and $\cK_{\pm}$ is short then  $u$ is transcendental over $K(q)$.
\end{enumerate}
\end{corollary}

{\it Proof}. Assertion 1 follows immediately by Theorem \ref{bloo}
and Lemmas \ref{irinusescoala} and \ref{gainuseli}, respectively.
In order to check assertion 2, note that by the above results, $\psi^1_q u$
is transcendental over $K(q)$. If $u$ were algebraic over $K(q)$,
the point 
$$(T(q),W(q))=(u,u^3+a_4u^7+\cdots )
$$
would have algebraic coordinates over $K(q)$. Hence, the same would be true
about the point
$$(x(q),y(q))=\left( \frac{T(q)}{W(q)}, -\frac{1}{W(q)} \right)\, ,$$
and therefore, $\psi^1_q u=\dd x(q)/y(q)$ would be algebraic over $K(q)$, a
contradiction. Thus, $u$ is transcendental over $K(q)$, and
assertion 2 is proved.
 \qed

\begin{corollary}
\label{dacaniciel}
Under the hypotheses of Theorem {\rm \ref{bloo}}, the maps
$B_{\pm}^0:\cU_{\pm 1} \ra R^{\rho_{\pm}}$ are $R$-module
isomorphisms. Furthermore, for any $u \in \cU_{\pm 1}$, we have
$$u=\sum_{\kappa \in {\mathcal K}_{\pm}} (B_{\kappa}^0 u) \star
u_{E,\kappa,1}\, .$$
\end{corollary}

In particular the ``boundary value problem at $q^{\pm 1}=0$'' is
well posed.

Next we address the ``boundary value problem at $q \neq 0$'' for
$\psi_E$.

If $E$ is an elliptic curve over $R$, we denote by $E(pR)$ the
kernel of the reduction modulo $p$ map $E(R) \ra E(k)$. As usual,
the group $E(pR)$ will be identified with $pR$ via the bijection
$$\begin{array}{rcl}
pR & \ra & E(pR) \\ pa & \mapsto &
(pa, (pa)^3+a_4(pa)^7+\cdots)
\end{array}\, .
$$

We denote by $E'(k)$ the group of all points in $E(k)$ of order prime to $p$.
And we denote by $E'(R)$ the subgroup of all points in $E(R)$ whose image
in $E(k)$ lie in $E'(k)$. There is a split exact sequence
$0 \ra E(pR) \ra E'(R) \ra E'(k) \ra 0$, hence
$E'(R)=E(pR) \oplus E'(R)_{tors}$.

\begin{corollary}
\label{haicasiel}
Under the hypotheses of Theorem {\rm \ref{bloo}}, and letting
$$
\cU'_+:=E'(R)_{tors}  \cdot \cU_{1} \subset \cU_+\, ,
$$
the following hold:
\begin{enumerate}
\item[1)] Assume $\cK_+=\{\kappa\}$. For any $q_0 \in p^{\nu}R^{\times}$
with $\nu \geq 1$, and any $g \in p^{\nu \kappa}R \subset
pR=E(pR)$ there exists a unique $u \in \cU_{1}$ such that
$u(q_0)=g$.
\item[2)] Assume $\cK_+=\{1\}$. Then for any $q_0 \in pR^{\times}$ and
any $g \in E'(R)$ there exists a unique $u \in  \cU'_+$ such that
$u(q_0)=g$.
\end{enumerate}
\end{corollary}

{\it Proof}. The first assertion follows exactly as in the case of $\bG_m$; cf.
Corollary \ref{maisus}.

In order to check the second assertion, note that
we can write $g$ uniquely as $g=g_{1}+g_{2}$, where $g_{1} \in
E'(R)_{tors}$ and $g_{2} \in E(pR)$. By the first part, there exists
$u_2 \in \cU_{1}$ such that $u_2(q_0)=g_{2}$. We set $u=g_{1}+u_2$.
Then, $u(q_0)=g$.

The uniqueness of $u$ is also checked easily.
\qed

The following Corollary is concerned with the inhomogeneous
equation $\psi_E u=\varphi$, and it is an immediate consequence of
Lemmas \ref{viyne}, \ref{dacaniciel}, and \ref{irinusescoala}.

\begin{corollary}
Let $\psi_E$ be a non-degenerate $\D$-character of $E$, and let
$\varphi \in q^{\pm 1}R[[q^{\pm 1}]]$ be a series whose support is
contained in the set $\cK'$ of totally non-characteristic integers
of $\psi_E$. Then the following hold:
\begin{enumerate} 
\item The equation $\psi_E u=\varphi$ has a unique solution $u \in E(q^{\pm
1}R[[q^{\pm 1}]])$ such that $\psi^1_q u$ has support disjoint
from the set $\cK$ of characteristic integers. 
\item If
$\bar{\varphi} \in k[q^{\pm 1}]$ the series $\overline{\psi^1_q u}
\in k[[q^{\pm 1}]]$ is integral over $k[q^{\pm 1}]$ and the field
extension $k(q) \subset k(q,\overline{\psi^1_q u})$ is Abelian
with Galois group killed by $p$. 
\item If the characteristic
polynomial of $\psi_E$ is unmixed and the support of $\varphi$ is
short then $u$ is transcendental over $K(q)$.
\end{enumerate}
\end{corollary}

\begin{remark}
 Corollary \ref{haicasiel} implies that if $\cK_+=\{1\}$ and $q_0 \in
pR^{\times}$, the group homomorphism
$$
\begin{array}{rcl}
S_{q_0}:E'(R)_{tors} \times R & \ra & E'(R) \\
(P,\alpha) & \mapsto & P+u_{E,1,\alpha}(q_0)
\end{array}
$$
is an isomorphism. So for any $q_1, q_2
\in pR^{\times}$, we have an isomorphism
$$
S_{q_1, q_2}:=S_{q_2} \circ S_{q_1}^{-1}:E'(R) \ra E'(R)\, .
$$
The latter map should be viewed as the ``propagator'' attached to
$\psi_E$.

As in the case of $\bG_a$ and $\bG_m$, respectively, if $\zeta \in
\bmu(R)$ and $q_0 \in p R^{\times}$, then, by (\ref{pie}),
$$
S_{\zeta q_0}=S_{q_0} \circ M_{\zeta} \\
$$
where
$$\begin{array}{rcl}
M_{\zeta}:E'(R)_{tors} \times R & \ra & E'(R)_{tors} \times  R \\
M_{\zeta}(P,\alpha) & := & (P,\zeta \alpha)
\end{array}
$$
Hence, for $\zeta_1, \zeta_2 \in \bmu(R)$, we obtain that
$$S_{\zeta_1 q_0,\zeta_2 q_0}=S_{q_0} \circ M_{\zeta_2/\zeta_1}
\circ S_{q_0}^{-1}\, .
$$
In particular,
$$
S_{q_0,\zeta_1 \zeta_2 q_0}=S_{q_0,\zeta_2 q_0} \circ
S_{q_0,\zeta_1 q_0}\, ,
$$
which can be interpreted  as a (weak) ``Huygens principle.''
\end{remark}

\bibliographystyle{amsalpha}

\end{document}